\documentclass[twoside, 11pt]{article}

\usepackage{amssymb}
\usepackage{amsfonts}
\usepackage{amsthm}
\usepackage{amsmath}
\usepackage{amscd}

\begin{document}

\pagestyle{myheadings}
\markboth{BASARAB}
{coGalois theory}

\title{A more general framework for coGalois theory\/}

\author{\c{S}ERBAN A. BASARAB \\
``Simion Stoilow'' Institute of Mathematics of the
Romanian Academy \\
P.O. Box 1--764\\
014700 Bucharest, ROMANIA\\
\texttt{{\itshape e-mail}: Serban.Basarab@imar.ro}}
\date{}  
\maketitle

\newtheorem{te}{Theorem}[section]
\newtheorem{pr}[te]{Proposition}
\newtheorem{rem}[te]{Remark}
\newtheorem{rems}[te]{Remarks}
\newtheorem{co}[te]{Corollary}
\newtheorem{lem}[te]{Lemma}
\newtheorem{prob}[te]{Problem}
\newtheorem{exs}[te]{Examples}
\newtheorem{defs}[te]{Definitions}
\newtheorem{de}[te]{Definition}
\newtheorem{que}[te]{Question}
\newtheorem{re}[te]{Remark}
\newtheorem{res}[te]{Remarks}
\newtheorem{ex}[te]{Example}
\newtheorem{aus}[te]{}

\newenvironment{aussage}%
{\renewcommand{\theequation}{\alph{equation}}\setcounter{equation}{0}%
\begin{aussage*}}{\end{aussage*}}

\def\R{\mathbb{R}}
\def\bbbm{{\rm I\!M}}
\def\N{\mathbb{N}}
\def\F{\mathbb{F}}
\def\H{\mathbb{H}}
\def\I{\mathbb{I}}
\def\K{\mathbb{K}}
\def\P{\mathbb{P}}
\def\D{\mathbb{D}}
\def\Q{\mathbb{Q}}
\def\Z{\mathbb{Z}}
\def\C{\mathbb{C}}
\def\S{\mathbb{S}}
\def\T{\mathbb{T}}
\def\L{\mathbb{L}}
\def\U{\mathbb{U}}

\def\cA{\mathcal{A}}
\def\cB{\mathcal{B}}
\def\cC{\mathcal{C}}
\def\cD{\mathcal{D}}
\def\cE{\mathcal{E}}
\def\cF{\mathcal{F}}
\def\cG{\mathcal{G}}
\def\cH{\mathcal{H}}
\def\cI{\mathcal{I}}
\def\cJ{\mathcal{J}}
\def\cK{\mathcal{K}}
\def\cL{\mathcal{L}}
\def\cM{\mathcal{M}}
\def\cN{\mathcal{N}}
\def\cO{\mathcal{O}}
\def\cP{\mathcal{P}}
\def\cQ{\mathcal{Q}}
\def\cR{\mathcal{R}}
\def\cS{\mathcal{S}}
\def\cT{\mathcal{T}}
\def\cU{\mathcal{U}}
\def\cV{\mathcal{V}}
\def\cW{\mathcal{W}}
\def\cX{\mathcal{X}}
\def\cY{\mathcal{Y}}
\def\cZ{\mathcal{Z}}


\newcommand{\notdiv}{{\not{|}\,}}
\newcommand{\ctg}{{\rm ctg\,}}
\newcommand{\sh}{{\rm sh\,}}
\newcommand{\ch}{{\rm ch\,}}
\newcommand{\Spin}{{\rm Spin\,}}
\newcommand{\spin}{{\rm spin\,}}
\newcommand{\Spek}{{\rm Spek\,}}
\newcommand{\Ker}{{\rm Ker\,}}
\newcommand{\di}{\mbox{\rm d\,}}
\newcommand{\cont}{\mbox{\rm cont\,}}
\newcommand{\rank}{\mbox{\rm rank\,}}
\newcommand{\codim}{\mbox{\rm codim\,}}
\newcommand{\ssqrt}[2]{\sqrt[\scriptstyle{#1}]{#2}}
\newcommand{\dbigcup}{\displaystyle\bigcup}
\newcommand{\dbigcap}{\displaystyle\bigcap}
\newcommand{\dinf}{\displaystyle\inf}
\newcommand{\dint}{\displaystyle\int}
\newcommand{\dsum}{\displaystyle\sum}
\newcommand{\dprod}{\displaystyle\prod}
\newcommand{\comb}[2]{\left(\begin{array}{c}#1\\#2\end{array}\right)}
\newcommand{\dlim}{\displaystyle\lim_{\stackrel{\longleftarrow}{}}}
\def\be{\begin{equation}}
\def\ee{\end{equation}}
\newcommand{\rs}{respectively }
\newcommand{\fZ}{Z^1(\Gam,A)}
\newcommand{\G}{\fZ}
\def\bp{\begin{proof}}
\def\ep{\end{proof}}
\def\hf{\hfill $\square$}
\def\ben{\begin{enumerate}}
\def\een{\end{enumerate}}
\def\ba{\begin{eqnarray*}}
\def\ea{\end{eqnarray*}}

\def\ls{\leqslant}
\def\gs{\geqslant}
\def\lra{\longrightarrow}
\def\hra{\hookrightarrow}
\def\Llra{\Longleftrightarrow}
\def\Lra{\Longrightarrow}
\def\p{\perp}

\newcommand{\wh}{\,\widehat}
\newcommand{\la}{\langle}
\newcommand{\ra}{\rangle}
\newcommand{\wt}{\widetilde}
\newcommand{\sm}{\setminus}
\newcommand{\sse}{\subseteq}
\newcommand{\es}{\varnothing}
\newcommand{\ov}{\overline}

\newcommand{\f}{\frac}
\newcommand{\q}{\quad}
\newcommand{\n}{\vartriangleleft}

\newcommand{\al}{h }
\newcommand{\De}{\Delta}
\newcommand{\del}{\delta}
\newcommand{\eps}{\varepsilon}
\newcommand{\gam}{\gamma}
\newcommand{\Gam}{\Gamma}
\newcommand{\Lam}{\Lambda}
\newcommand{\lam}{\lambda}
\newcommand{\om}{\omega}
\newcommand{\Om}{\Omega}
\newcommand{\si}{\sigma}


\def\fa{\mathfrak{a}}
\def\fb{\mathfrak{b}}
\def\fc{\mathfrak{c}}
\def\fd{\mathfrak{d}}
\def\fe{\mathfrak{e}}
\def\ff{\mathfrak{f}}
\def\fg{\mathfrak{g}}
\def\fh{\mathfrak{h}}
\def\fp{\mathfrak{p}}
\def\fq{\mathfrak{q}}
\def\fm{\mathfrak{m}}
\def\fn{\mathfrak{n}}
\def\fs{\mathfrak{s}}
\def\ft{\mathfrak{t}}
\def\fu{\mathfrak{u}}
\def\fA{\mathfrak{A}}
\def\fE{\mathfrak{E}}
\def\fF{\mathfrak{F}}
\def\fM{\mathfrak{M}}
\def\fG{\mathfrak{G}}
\def\fH{\mathfrak{H}}
\def\fT{\mathfrak{T}}
\def\fX{\mathfrak{X}}
\def\fY{\mathfrak{Y}}
\def\fP{\mathfrak{P}}
\def\fL{\mathfrak{L}}

\begin{abstract}
The paper is an extended version of a talk given at the {\em International Conference: Experimental and 
Theoretical Methods in Algebra, Geometry and Topology}, which took place June 21--24, 2013 in Eforie Nord (Romania). 
Its purpose is to present a more general framework for a fairly new theory in Field Theory, called 
{\em coGalois Theory}, which is somewhat dual to the very classical {\em Galois Theory} and is more general 
than the {\em Kummer Theory}. The main object of investigation in this more general framework is the {\em coGalois
connexion} naturally associated to any triple $\,(\Gam, \fG, \eta)\,$, where $\,\Gam\,$ is a profinite group, $\,\fG\,$
is a profinite operator $\,\Gam$-group, and $\,\eta: \Gam \lra \fG\,$ is a continuous $\,1$-cocycle with the 
property that the profinite group $\,\fG\,$ is topologically generated by the image $\,\eta(\Gam)$.
\smallskip

\smallskip

\end{abstract}

\section{Introduction}

The so called {\em coGalois theory} is a more or less recent development \cite{Orozco, G-H, B-R-V, A-N1, A-N2, 
B-L, TA1} of the study of finite radical extensions, carried out, in chronological order, by H. Hasse (1930), 
A. Besicovitch (1940) \cite{Bes}, L.J. Mordell (1953) \cite{Mordell}, C.L. Siegel (1972) \cite{Siegel}, 
M. Kneser (1975) \cite{K} and A. Schinzel (1975) \cite{Schinzel}, among others. Roughly speaking, 
the coGalois theory is an extension of the {\em Kummer theory}, being somewhat dual to the very classical 
{\em Galois theory}. For details and more references the reader may consult the monograph \cite{CT}.

More precisely, given a radical field extension $L/K$, the main goal of the coGalois theory is to investigate 
the relation between the intermediate fields of the extension $L/K$ and the subgroups of the torsion subgroup 
of the multiplicative factor group $L^\times/K^\times$, called the {\em coGalois group} of the extension $L/K$, 
denoted by $\text{coG}(L/K)$. Assuming that the extension $L/K$ is Galois with $\Gam := \text{Gal}(L/K)$, 
$\text{coG}(L/K)$ is canonically isomorphic, via Hilbert's Theorem 90, to the group
$\,Z^1(\Gam,\mu(L))\,$ of all continues 1-cocycles (crossed homomorphisms) of the profinite group $\Gam$ with 
coefficients in the group $\,\mu(L)\,$ of all roots of unity in $L$ (\cite{B-R-V}). Notice that the multiplicative group 
$\,\mu(L)\,$ is quasi-cyclic, so it is isomorphic (in a noncanonical way) to a subgroup of the additive group $\,\Q/\Z$.

By analogy with Neukirch's {\em Abstract Galois Theory} within his {\em Abstract Class Field Theory} \cite{N},
an abstract, group theoretic framework of the coGalois theory is developed in the papers \cite{A-B, Serdica, JPAA}, 
where basic concepts of the field theoretic coGalois theory, as well as their main properties, are generalized 
to arbitrary continuous actions of profinite groups on discrete quasi-cyclic groups.

Such a continuous action $\Gam \times A \lra A$, where $\Gam$ is a profinite group and the discrete quasi-cyclic 
group $A$ is identified with a subgroup of $\Q/\Z$, establishes through the evaluation map 
$\,\Gam \times \G \lra A,\, (\gam, \alpha) \mapsto \alpha(\gam),\,$ a Galois connexion between the
lattice $\,\L(\Gam)\,$ of all closed subgroups of $\Gam$ and the lattice $\,\L(\G)\,$ of all subgroups of 
the torsion abelian group $\,\G\,$ of continuous $1$-cocycles.
 
On the other hand, the continuous action of $\,\Gam\,$ on $A$ endows the dual group 
$\,\G^\vee = {\rm Hom}(\G,\Q/\Z) = {\rm Hom}(\G, A)\,$ with a natural structure of profinite $\,\Gam$-module, related 
to $\,\Gam\,$ through a continuous cocycle $\,\eta : \Gam \lra \G^\vee\,$ with the property that the abelian 
profinite group $\,\G^\vee\,$ is topologically generated by its closed subset $\eta(\Gam)$. The continuous cocycle 
$\eta$ plays a key role in \cite{A-B, Serdica, JPAA} for the study of several interesting classes of subgroups of 
$\Gam$ and $\G$ induced by the Galois connexion above ({\em radical, hereditarily radical, Kneser, hereditarily Kneser, 
coGalois} and {\em strongly coGalois}).

The major role played by the continuous cocycle $\eta$ above in what we may call {\em cyclotomic abstract 
coGalois theory} is the motivation for a more general approach of coGalois theory with potential new applications.
The present paper, based on the unpublished preprint \cite{Preprint}, is organized in 6 sections (2--7). Some basic 
notions used throughout the paper (profinite and spectral spaces, profinite groups and operator groups, subgroup 
lattices, Galois and coGalois connexions) are briefly explained in Section 2. Section 3 introduces a general abstract 
framework for coGalois theory having as main object of investigation the triples $\,(\Gam, \fG, \eta)\,$, where 
$\,\Gam\,$ is a profinite group, $\,\fG\,$ is a profinite operator $\,\Gam$-group, and $\,\eta: \Gam \lra \fG\,$ is a 
continuous $\,1$-cocycle with the property that $\,\eta(\Gam)\,$ topologically generates the profinite group $\,\fG\,$;
such a cocycle $\,\eta\,$ is called {\em generating cocycle}. 
To any such triple one assigns a natural coGalois connexion between the lattice $\,\L(\Gam)\,$ of all closed subgroups
of $\,\Gam\,$ and the modular lattice $\,\L(\fG)\,$ of all closed $\,\Gam$-invariant normal subgroups of $\,\fG\,$,
called {\em ideals} of the profinite $\,\Gam$-group $\,\fG\,$. The main properties of this coGalois connexion are 
collected in Proposition~\ref{Pr:MainProperties}.

The general framework presented in Section 3 is applied in Section 4 to an abelian context which extends as most as 
possible the framework of cyclotomic coGalois theory by considering continuous actions of profinite groups $\,\Gam\,$ on 
arbitrary discrete torsion abelian groups $\,A\,$, not necessarily quasi-cyclic. To any such $\,\Gam$-module $\,A\,$ 
one assigns a profinite $\,\Gam$-module $\,\fG\,$ and a generating cocycle $\,\eta: \Gam \lra \fG\,$ such that the 
Galois connexion between the lattice $\,\L(\Gam)\,$ and the lattice $\,\L(Z^1(\Gam, A))\,$ of all subgroups of the
torsion abelian group $\,Z^1(\Gam, A)\,$ is obtained by composing the coGalois connexion between $\,\L(\Gam)\,$ and 
$\,\L(\fG)\,$, introduced in Section 3, with a natural Galois connexion between the lattices $\,\L(\fG)\,$ and 
$\,\L(Z^1(\Gam, A))\,$ (Propositions~\ref{Pr:GalCon1}, \ref{Pr:coGCon}, \ref{Pr:GalCon2}, Corollary~\ref{Cor:Composing}).
Section 4 ends with 4 relevant examples: the first two examples are concerned with the coGalois theory of separable 
radical extensions and its abstract cyclotomic version; the third example is devoted to an additive analogue of 
the coGalois theory of separable radical extensions, based on Witt calculus and higher Artin-Schreier theory, 
while the fourth example is an extension of the cyclotomic context to Galois algebras. 

Some of the main notions and results of the cyclotomic abstract coGalois theory are extended in Sections 5 and 6 to
the general framework introduced in Section 3. Here two remarkable types of triples $\,(\Gam, \fG, \eta)\,$ are 
investigated: the {\em Kneser triples}, where the cocycle $\,\eta: \Gam \lra \fG\,$ is surjective, and the 
{\em coGalois triples}, i.e,, Kneser triples for which the associated coGalois connexion is perfect. Given a triple
$\,(\Gam, \fG, \eta)\,$, where $\,\eta\,$ is a generating cocycle, the main properties of the space $\,\K(\fG)\,$
($\,\C\mathbb{G}(\fG)\,$), consisting of the ideals $\,{\bf a}\,$ of the profinite $\,\Gam$-group $\,\fG\,$ for
which the induced triple $\,(\Gam, \fG/{\bf a}, \eta_{\bf a}: \Gam \lra \fG/{\bf a})\,$ is Kneser (coGalois), are
collected in Propositions~\ref{Pr:KneserSpace}, \ref{Pr:coGaloisSpace}. A special attention is paid in \ref{Subsubsec:BijectCocycles}
to a procedure for obtaining bijective cocycles by deformation of a profinite group via a continuous action on itself.
General {\em Kneser} and {\em coGalois criteria} are provided by Propositions~\ref{Pr:AbsKneserCrit}, 
\ref{Pr:AbscoGaloisCrit}, and two remarkable classes of finite structures ({\em minimal non-Kneser} and 
{\em Kneser minimal non-coGalois triples}) arising naturally from these general criteria are introduced. The open 
Problems~\ref{Prob:mnK}, \ref{Prob:mncG} are concerned with the classification of these finite structures. 

Partial answers to Problem~\ref{Prob:mnK} are given in Section 7. In particular, the special case when $\,\Gam,\,\fG\,$ 
are abelian is completely solved under an additional assumption on the local subring of the endomorphism ring of the
abelian $\,p$-group $\,\fG\,$, generated by $\,\Gam\,$ (Propositions~\ref{Pr:Ab1}, \ref{Pr:Ab2}). As an immediate consequence, 
we find again \cite[Lemma 1.18, Theorem 1.20]{A-B}, the abstract version of the classical {\em Kneser criterion for 
separable radical extensions} \cite{K}, \cite[Theorem 11.1.5]{CT}.

\section{Preliminaries}

\subsection{Profinite and spectral spaces} 
A {\em profinite space}, also called {\em Stone} or {\em boolean space}, is a compact Hausdorff and totally
disconnected topological space. A {\em spectral space}, also called {\em coherent} or {\em quasi-boolean space} is  
a compact $\,T_0$ topological space which admits a base of compact sets for its topology; equivalently, a topological
space $\,X\,$ is spectral if the family of compact open sets is closed under finite intersections (in particular, $\,X\,$
itself is compact) and forms a base for the topology of $\,X\,$, and $\,X\,$ is {\em sober}, i.e., every irreducible 
closed subset of $\,X\,$ is the closure of a unique point of $\,X$.

The spectral spaces form a category $\,{\bf SPECS}\,$ having as morphisms the so called {\em coherent maps},
i.e., the maps $\,f : X \lra Y\,$ for which $\,f^{- 1}(V)\,$ is a compact open subset of the spectral space 
$\,X\,$ provided $\,V\,$ is a compact open subset of the spectral space $\,Y\,$. In particular, the coherent 
maps are continuous, so $\,{\bf SPECS}\,$ is a non-full subcategory of the category $\,{\bf TOP}\,$ of
topological spaces. We denote by $\,{\bf PFS}\,$ the full subcategory of $\,{\bf TOP}\,$ whose objects are
the profinite spaces. Obviously, $\,{\bf PFS}\,$ is also a full subcategory of $\,{\bf SPECS}\,$.

By {\em Stone Representation Theorem}, $\,{\bf PFS}\,$ and $\,{\bf SPECS}\,$ are duals to the categories
of boolean algebras and (bounded) distributive lattices, respectively. For more details concerning the topological 
spaces above and their dual structures the reader may consult \cite{Hochster, J, Pop, Dual}. 
\subsection{Profinite groups}
\label{Subsec:ProfGps}
A {\em profinite group} is a group object in $\,{\bf PFS}\,$, i.e., a compact totally disconnected topological group;
equivalently, a topological group $\,\Gam$ is profinite if the identity element $\,1$ of $\,\Gam$ admits a 
fundamental system $\,\cal{U}$ of open neighborhoods $\,U$ such that $\,U$ is a normal subgroup of $\,\Gam$, and 
$\,\Gam = \displaystyle\lim_{\stackrel{\longleftarrow}{U \in \cal{U}}}\,\Gam/U$, the inverse limit of the inverse 
system of discrete finite groups $\,\{\Gam/U\,|\,U \in \cal{U}\}\,$. For details on profinite groups see 
\cite{SER, NSW, RibesZ}.

Let $\,\cC\,$ be a {\em variety} of finite groups, i.e., a nonempty class of finite groups, closed under 
subgroups, quotients, and finite direct products. The profinite groups $\,\Gam\,$, with $\,\Gam/U \in \cC\,$ for all
open normal subgroups $\,U \sse \Gam\,$, called {\em pro-$\,\cC\,$ groups}, form a full subcategory 
$\,{\bf P}\,\cC\,{\bf G}\,$ of the category $\,{\bf PFG}\,$ of profinite groups, with continuous homomorphisms. 
In particular, for any subgroup $\,A \sse \Q/\Z\,$, let $\,{\bf FAb}_A\,$ denote the variety consisting of those 
finite abelian groups $\,G\,$ for which $\,\frac{1}{{\rm exp}(G)} \Z/\Z \sse A\,$. For instance, taking 
$\,A = \Q/\Z,\,(\Q/\Z)(p) \cong \Q_p/\Z_p,\,\frac{1}{n} \Z/\Z\,$, we obtain the varieties of finite abelian groups, of 
finite abelian $\,p$-groups for any prime number $\,p\,$, and of finite abelian groups of exponent dividing $\,n \geq 1$.  
Moreover any variety of finite abelian groups is of the form above. Setting $\,A^\vee:= {\rm Hom}(A, \Q/\Z)$, the 
Pontryagin dual of the discrete quasi-cyclic group $\,A\,$, with the canonical structure of profinite ring, the 
pro-$\,{\bf FAb}_A\,$ groups are identified with the profinite $\,A^\vee$-modules, the duals of discrete (torsion) 
$\,A^\vee$-modules.  

\subsection{The topological lattice of closed subgroups of a profinite group}
\label{Subsec:LatticeSgps}
For any profinite group $\,\Gam\,$, we denote by $\,\L(\Gam)\,$ the poset with respect to inclusion of all closed
subgroups of $\,\Gam\,$; moreover $\,\L(\Gam)\,$ is a bounded lattice with obviously defined operations 
$\,\wedge = \cap\,$ and $\,\vee\,$. In addition, $\,\L(\Gam)\,$ becomes a spectral space as a closed subspace 
of the spectral space of all closed subsets of the underlying profinite space of $\,\Gam\,$. The {\em spectral topology}
$\,\tau_s\,$ on $\,\L(\Gam)\,$ is defined by the base of compact open sets $\,\L(\De)\,$ for $\,\De\,$ ranging over 
all open subgroups of $\,\Gam$. Note that for any $\,\Lam \in \L(\Gam)$, the closure of the one-point set
$\,\{ \Lam \}\,$ is $\,\overline{\{\Lam\}} = \L(\,\Gam \,|\, \Lam)$, the set of all closed subgroups of 
$\,\Gam$ lying over $\,\Lam$. Thus the spectral space $\,\L(\Gam)\,$ is irreducible with the generic point 
$\,\{ 1 \}\,$, while $\,\Gam\,$ is its unique closed point. Since the poset $\,\L(\Gam)\,$ is the inverse limit of 
the inverse system of finite posets $\,\L(\Gam/\De)\,$ for $\,\De\,$ ranging over $\,\cN(\Gam)\,$, the set of all 
open normal subgroups of $\,\Gam$, with natural order-preserving connecting maps, the topology $\,\tau_s\,$ is 
exactly the inverse limit of the $T_0$ topologies induced by the partial order given by inclusion on the finite 
sets $\,\L(\Gam/\De)$.

The {\em boolean (profinite) completion\/} $\,\tau_{b}\,$ of the spectral topology $\,\tau_{s}\,$ on $\,\L(\Gam)$, 
also called the {\em patch topology\/}, is the topology with the base of clopen sets
${\cal V}_{\De,\De'}=\{\,\Lam \in \L(\Gam)\,|\, \Lam \De = \De' \}$
for all pairs $\,(\De, \De')$, with $\,\De \in \cN(\Gam)\,$ and $\,\De'\in \L(\Gam | \De)$. The profinite space 
above is the inverse limit of the discrete finite spaces $\L(\Gam/\De)$ for $\De$ ranging over $\,\cN(\Gam)\,$. 
A subset $\,\cal U\,$ of $\,\L(\Gam)\,$ is $\,\tau_{s}$-open if and only if  $\,\cal U\,$ is both $\,\tau_{b}$-open 
and  a {\em lower} subset of $\,\L(\Gam)$; the later condition means that $\Lam\in {\cal U}\,$ and 
$\,\Lam'\in \L(\Lam)$ imply $\,\Lam' \in {\cal U}$.

\begin{rem}
\label{Rem:Coherence} \em
One checks easily that the canonical action of the profinite group $\,\Gam$ on the spectral space $\,\L(\Gam)\,$,
$\,(\gam, \Lam) \mapsto \gam \Lam \gam^{- 1}\,$, and the join operation $\,(\Lam_1, \Lam_2) \mapsto \Lam_1 \vee
\Lam_2$ are coherent maps, in particular, continuous, while the meet operation $\,(\Lam_1, \Lam_2) \mapsto \Lam_1 \cap 
\Lam_2\,$ is continuous, not necessarily coherent.
\end{rem}
\subsection{Profinite operator groups}
\label{Subsec:OperatorGroups}

Let $\,\Gam\,$ be a profinite group. By a {\em profinite (operator) $\Gam$-group} we understand a profinite group 
$\,\fG\,$ together with a continuous action by automorphisms $\Gam \times \fG \lra \fG,\, (\gam, g) \mapsto \gam g\,$;
equivalently, the profinite group $\,\fG\,$ possesses a system of neighbourhoods of the identity consisting of 
open $\Gam$-invariant normal subgroups. Denote by $\,{\rm Fix}_\Gam(\fG)\,$ the kernel of the action, a closed 
normal subgroup of $\,\Gam\,$. The abelian profinite $\,\Gam$-groups are usually called {\em profinite 
$\,\Gam$-modules}. For any nontrivial subgroup $\,A \sse \Q/\Z\,$, the pro-$\,{\bf FAb}_A\,$ $\,\Gam$-groups are
identified with the profinite (left) $\,A^\vee[[\Gam]]$-modules.

We denote by $\,{\bf PFOG}$ the category having as objects the pairs $\,(\Gam, \fG)\,$, where $\,\Gam\,$ is a
profinite group and $\,\fG\,$ is a profinite $\,\Gam$-group. The morphisms $(\Gam, \fG) \lra (\Gam', \fG')$ are
pairs $\,(\varphi : \Gam \lra \Gam', \psi : \fG \lra \fG')$ of continuous homomorphisms such that 
$\psi(\gam g) = \varphi(\gam) \psi(g)$ for all $\gam \in \Gam, g \in \fG$. The composition law in $\,{\bf PFOG}\,$ 
is naturally defined. On the other hand, we consider the category $\,{\bf SEPI}$ of {\em splitting epimorphisms} 
having as objects the tuples $\,(\Gam, \fE, p, s)$ consisting of profinite groups $\,\Gam, \fE$, an epimorphism 
$\,p : \fE \lra \Gam$, and a continuous homomorphic section $\,s : \Gam \lra \fE$ of $\,p$. As morphisms 
$\,(\Gam, \fE, p, s) \lra (\Gam',\fE', p', s')$ in $\,{\bf SEPI}$, we take the pairs $\,(\varphi : \Gam \lra \Gam',
\psi : \fE \lra \fE')$ of profinite group morphisms satisfying $\,\varphi \circ p = p' \circ \psi$ and 
$\,\psi \circ s = s' \circ \varphi$, with the natural composition law. 
\begin{lem}
\label{Lem:EqPFOG-SEPI}
The categories $\,{\bf PFOG}$ and $\,{\bf SEPI}$ are equivalent.
\end{lem}
\bp
Let $\,F : {\bf PFOG} \lra {\bf SEPI}$ denote the covariant functor induced by the map assigning to any profinite
operator group $\,(\Gam, \fG)$ the object $\,(\Gam, \fE, p, s)$ of $\,{\bf SEPI}$, where $\,\fE\,$ is the semidirect 
product $\,\fG \rtimes \Gam\,$ induced by the continuous action of $\,\Gam\,$ on $\,\fG,\,p : \fE \lra \Gam\,$ is 
the natural projection with $\,\Ker(p) = \fG$, and $\,s : \Gam \lra \fE\,$ is the canonical homomorphic section of 
$\,p\,$ which identifies $\,\Gam\,$ with a closed subgroup of $\,\fE\,$ satisfying $\fG \cdot \Gam = \fE,\, 
\fG \cap \Gam = \{ 1 \}\,$. One checks easily that the functor $\,F\,$ is faithfully full and essentially surjective, 
and hence it yields an equivalence of categories. Note that the inverse equivalence $\,{\bf SEPI} \lra {\bf PFOG}\,$
is induced by the map assigning to an object $\,(\Gam, \fE, p, s)\,$ of $\,{\bf SEPI}$ the profinite $\,\Gam$-group 
$\,\fG = \Ker(p)\,$ with the action of $\,\Gam\,$ defined by $\,\gam g := s(\gam)\,g\,s(\gam)^{- 1}\,$ for
$\,\gam \in \Gam, g \in \fG\,$.
\ep

Note that the equivalent categories $\,{\bf PFOG}\,$ and $\,{\bf SEPI}\,$ have inverse limits. Moreover the equivalence
above is naturally extended to an equivalence of suitable categories of bundles containing $\,{\bf PFOG}\,$ and
$\,{\bf SEPI}\,$ as reflective subcategories \cite{Preprint}.

\subsubsection{The lattice of ideals of a profinite operator group}
\label{Subsubsec:LatticeIdeals}
For any profinite $\,\Gam$-group $\,\fG\,$, we denote by $\,\L(\fG)\,$ the poset with respect to inclusion of all
$\,\Gam\,$-invariant closed normal subgroups of $\,\fG\,$, called {\em ideals} of the profinite $\,\Gam$-group
$\,\fG\,$. The ideals of $\,\fG\,$ are exactly those closed normal subgroups of the semidirect product 
$\,\fG \rtimes \Gam\,$ which are contained in $\,\fG\,$. $\,\L(\fG)\,$ is a bounded {\em modular} lattice, dual to the 
lattice of all quotients of the profinite $\,\Gam$-group $\,\fG\,$; an ideal $\,{\bf a} \in \L(\fG)\,$ is open if
and only if the quotient $\,\Gam$-group $\,\fG/{\bf a}$ is finite. 

In addition, $\,\L(\fG)\,$ is equipped with a spectral topology defined by the base of compact open sets
$\L({\bf a}) = \{\,{\bf b} \in \L(\fG)\,|\,{\bf b} \sse {\bf a}\,\}$
for $\,{\bf a}\,$ ranging over the open ideals of $\,\fG\,$. The join operation $\,({\bf a},\,{\bf b}) \mapsto
{\bf a} \vee {\bf b}:= {\bf a} \cdot {\bf b}\,$ is coherent, while the meet operation $\,({\bf a},\,{\bf b}) \mapsto
{\bf a} \cap {\bf b}\,$ is continuous, not necessarily coherent.

\subsection{Galois and coGalois connexions}
The notions of Galois and coGalois connexions are remarkable special cases of the more general concept of 
{\em adjunction} from category theory.

According to Ore \cite{Ore}, a {\em Galois connexion} is a system $\,(X, Y, \varphi, \psi)\,$, where $\,X, Y\,$ are
posets, and $\,\varphi:X \lra Y$, $\,\psi:Y \lra X\,$ are {\em order-reversing} maps satisfying  
$\,x \leq \psi(\varphi(x))\,$, $\,y \leq \varphi(\psi(y))\,$ for all $\,x \in X, y \in Y$. It follows that the maps
$\,\varphi\,$ and $\,\psi\,$ are {\em quasi-inverse} to one another, i.e., $\,\varphi \circ \psi \circ \varphi = 
\varphi,\,\psi \circ \varphi \circ \psi = \psi\,$, and the order-preserving endomaps $\,\psi \circ \varphi: X \lra X$,
$\,\varphi \circ \psi:Y \lra Y\,$ are {\em closure operators} with $\,X_{\rm c}:= \{x \in X\,|\,\psi(\varphi(x)) = x\} = 
\psi(Y)\,$, $\,Y_{\rm c}:= \{y \in Y\,|\,\varphi(\psi(y)) = y\} = \varphi(X)\,$ as sets of closed points. The Galois
connexion $\,(X, Y, \varphi, \psi)\,$ is {\em perfect} if $\,X = X_{\rm c}, Y = Y_{\rm c}\,$, i.e., the maps
$\,\varphi: X \lra Y,\,\psi:Y \lra X\,$ are {\rm anti-isomorphisms} inverse to one another.

By duality, a {\em coGalois connexion} is a system $\,(X, Y, \varphi, \psi)\,$, where $\,X, Y\,$ are posets, and
$\,\varphi:X \lra Y, \psi:Y \lra X\,$ are {\em order-preserving} maps satisfying $\,x \leq \psi(\varphi(x))\,$,
$\,\varphi(\psi(y)) \leq y\,$ for all $\,x \in X, y \in Y\,$. The coGalois connexion $\,(X, Y, \varphi, \psi)\,$ is
{\em perfect} if $\,\varphi: X \lra Y, \psi: Y \lra X\,$ are {\em isomorphisms} inverse to one another.

\subsubsection{Standard Galois connexions} 
Let $\,L/K\,$ be an arbitrary field extension. According to \cite{Krull1}, \cite[6.3]{Shimura}, 
$\,\Gam:= {\rm Gal}(L/K) = \{\gam \in {\rm Aut}(L)\,:\,\si|_K = 1_K\}\,$ has a natural structure of 
{\em totally disconnected (Hausdorff) topological group} with respect to the weakest topology for which the action 
of $\,\Gam\,$ on the discrete field $\,L\,$ is continuous; the subgroups $\,{\rm Gal}(L/F)\,$, where 
$\,K \sse F \sse L\,$ ranges over the finitely generated field extensions of $\,K\,$, form a system of 
open neighbourhoods of the identity $\,1_L\,$. For algebraic extensions $\,L/K\,$, $\,\Gam\,$ is compact, whence 
profinite, and the topology above is usually called {\em Krull topology}. Note that $\,\Gam\,$ is locally compact 
if and only if $\,{\rm tr.deg}\,L/K < \infty\,$.

Consider the bounded lattices (with respect to inclusion) $\,\cL(L/K)\,$ of all subfields of $\,L\,$ containing 
$\,K\,$, $\,\cL(\Gam)\,$ of all subgroups of $\,\Gam\,$, and $\,\L(\Gam)\,$ of all closed  
subgroups of $\,\Gam\,$; note that $\,\L(\Gam)\,$ is a sub-semilattice of $\,\cL(\Gam)\,$ with respect to 
the meet operation (intersection), while its retract $\,\cL(\Gam) \lra \L(\Gam),\,\Lam \mapsto \overline{\Lam}\,$,
is a morphism of semilattices with respect to the join operations. These lattices are naturally related through 
the order-reversing maps  
$\,\cD: \cL(\Gam) \lra \cL(L/K),\,\Lam \mapsto L^\Lam:=$ $\,\{x \in L\,|\,\forall \gam \in \Lam, \gam(x) = x\}\,$ 
({\em Dedekind connexion}), and $\,\cK: \cL(L/K) \lra \cL(\Gam), F \mapsto {\rm Gal}(L/F)\,$ ({\em Krull 
connexion}). Note that $\,\cD\,$ factorizes trough the canonical surjection $\,\cL(\Gam) \lra \L(\Gam),\,\Lam 
\mapsto \overline{\Lam}\,$, while $\,\cK\,$ factorizes through the embedding $\,\L(\Gam) \hra \cL(\Gam)\,$.  

The systems $\,(\cL(L/K), \cL(\Gam), \cK, \cD)\,$ and $\,(\cL(L/K), \L(\Gam), \cK, \cD)\,$ are both Galois connexions.
The first one is perfect if and only if $\,L/K\,$ is a {\em finite Galois extension} (solution of 
{\em Steinitz's problem} \cite{Steinitz}), while the latter one is perfect if and only if $\,L/K\,$ is an 
{\em algebraic} (not necessarily finite) {\em Galois extension} \cite{Krull1}. In the first case, the Galois group 
$\,\Gam\,$ is finite of order $\,[L : K]\,$, while in the latter case, $\,\Gam\,$ is a profinite group isomorphic to 
$\,\displaystyle\lim_{\stackrel{\longleftarrow}{}} {\rm Gal}(F/K)\,$, $\,F\,$ ranging over the finite normal 
extensions of $\,K\,$ contained in $\,L$.

\begin{res} \em
(1) According to Barbilian \cite{DB}, a field extension $\,L/K\,$ is {\em Dedekindian} (i.e., 
$\,\forall F \in \cL(L/K)\,$, $\cD(\cK(F)) = F^{{\rm Gal}(L/F)} = F\,$) if and only if for all $\,F \in \cL(L/K)\,$, the
relative algebraic closure of $\,F\,$ in $\,L\,$ is Galois over $\,K\,$. In particular, assuming $\,{\rm char}\,K > 0\,$,
$\,L/K\,$ is Dedekindian if and only if $\,L/K\,$ is an algebraic Galois extension.

(2) Extending Barbilian's paper \cite{DB}, Krull investigates in \cite{Krull2} the field extensions $\,L/K\,$, called 
locally normal by Barbilian, and simply {\em normal} by Krull, having the property that for all 
$\,F \in \cL(L/K)\,$, the relative algebraic closure of $\,F\,$ in $\,L\,$ is normal over 
the base field $\,K\,$. He shows that an arbitrary field extension $\,L/K\,$ is normal if and only if for every 
{\em Steinitz decomposition} $\,F \in \cL(L/K)\,$ (i.e., $\,F\,$ purely transcendental over $\,K\,$, and $\,L/F\,$ 
algebraic), the algebraic extension $\,L/F\,$ is normal (necessarily infinite), therefore the purely transcendental
extensions are not normal. Using this general notion of normality, Barbilian's main result reads as follows: 
$\,L/K\,$ {\em is Dedekindian if and only if $\,L/K\,$ is normal and separable}. The following open problem is raised 
by Krull in \cite{Krull2}: {\em Do there exist transcendental normal extensions which are not algebraically closed ?} 
\footnote[1]{Unfortunately, I couldn't find in literature some references to this fundamental problem stated by Krull 
in 1953. I found only a paper in Arch. Math. {\bf 61} (1993), 238--240, with the title {\em Gibt es nichttriviale 
absolut-normale K\" orpererweiterungen ?}, but, strangely enough, this paper contains only some of Barbilian and 
Krull's notions and results without to mention their papers in bibliography ! }.

(3) Though the Galois algebraic extensions play a key role in the most applications of the field arithmetic, there
exist also significant contexts in which the Galois groups of transcendental extensions are essential tools of 
investigation. As a relevant example, we mention the automorphism group of the field of modular functions, rational
over the maximal abelian extension of $\,\Q\,$ \cite[Ch. 6]{Shimura}. Another situation, treated in \cite{Rov}, concerns
finitely generated extensions $\,F/K\,$ of transcendence degree $\,1\,$ over an algebraically closed field $\,K\,$ of
characteristic $\,0\,$. In this case, the {\em absolute Galois group} $\,{\rm Gal}(\wt{F}/F)\,$ is a free profinite
group, and hence, its structure tells nothing about the field $\,F\,$. However, considering two finitely generated 
extensions $\,F_1, F_2\,$ of $\,K\,$, contained into an algebraically closed field $\,L\,$ of transcendence degree 
$\,1\,$ over $\,K\,$, the author proves that $\,F_1\,$ and $\,F_2\,$ are isomorphic over $\,K\,$, provided there is 
a continuous and open automorphism of $\,{\rm Gal}(L/K)\,$ inducing by restriction an isomorphism 
$\,{\rm Gal}(L/F_1) \cong {\rm Gal}(L/F_2)\,$.  
\end{res}

\section{An abstract framework for coGalois Theory}
\label{Sec:GenFrame}

\subsection{Generating cocycles}
\label{Subsec:GenCocycles}
For any profinite $\,\Gam$-group $\fG$, let $\,Z^1(\Gam, \fG)\,$ denote the set of all continuous $1$-cocycles 
(crossed homomorphisms) of $\,\Gam\,$ with coefficients in $\fG$, i.e., the continuous maps $\eta: \Gam \lra \fG$ 
satisfying $\eta(\si \tau) = \eta(\si) \cdot \si \eta(\tau)\,$ for all $\,\si, \tau \in \Gam\,$;
in particular, $\eta(\gam^{- 1}) = \gam^{- 1} \eta(\gam)^{- 1}$ for all $\gam \in \Gam$, and $\eta(1) = 1$.
The set $Z^1(\Gam, \fG)$ contains the trivial cocycle $\gam \mapsto 1$ as a privileged element, so it is an object
of the category $\,{\bf PS}\,$ of pointed sets with naturally defined morphisms. $\,Z^1(\Gam, \fG)\,$ becomes an
abelian group whenever $\fG$ is a profinite $\,\Gam$-module.

Note that Ker$\,(\eta) := \eta^{- 1}(1)\,$ is a closed subgroup of $\,\Gam\,$ for every $\,\eta \in Z^1(\Gam, \fG)\,$.
Set $\,Z^1(\Gam\,|\,\Lam, \fG):= \{\eta \in Z^1(\Gam, \fG)\,|\,\Lam \sse \Ker\,(\eta)\}\,$, where $\,\Lam\,$ is a closed
subgroup of $\,\Gam$. 

We denote by $\,\cZ^1\,$ the {\em category of $\,1$-cocycles} whose objects are the triples $\,(\Gam, \fG, \eta)\,$ 
consisting of a profinite group $\,\Gam\,$, a profinite $\,\Gam$-group $\,\fG\,$ and a cocycle 
$\,\eta \in Z^1(\Gam, \fG)\,$. As morphisms $\,(\Gam, \fG, \eta) \lra (\Gam', \fG', \eta')\,$ we take those morphisms 
$\,(\varphi : \Gam \lra \Gam',\,\psi : \fG \lra \fG')\,$ in the category $\,{\bf PFOG}$ of profinite operator 
groups which in addition are compatible with the cocycles $\,\eta$ and $\,\eta'$, i.e., 
$\,\psi \circ \eta = \eta' \circ \varphi$. The composition law in $\,\cZ^1\,$ is induced from the category
$\,{\bf PFOG}$. 
\begin{de}
\label{Def:GenerCocycle}
$\,\eta \in Z^1(\Gam, \fG)\,$ is said to be a {\em generating cocycle} (for short {\em g-cocycle}) if the profinite 
group $\,\fG\,$ is topologically generated by its closed subset $\,\eta(\Gam)$.
\end{de}
For any g-cocycle $\,\eta \in Z^1(\Gam, \fG)$, we obtain 
\[
\text{Fix}\,_\Gam(\fG) = \text{Fix}\,_\Gam(\eta(\Gam)) = \{ \gam \in \Gam\,|\,\forall \si \in \Gam,\,\eta(\gam \si) =
\eta(\gam) \eta(\si)\,\},
\] 
while the core of $\De:= \text{Ker}(\eta)$ is 
$\,\wt{\De} := \bigcap_{\gam \in \Gam} \gam \De \gam^{- 1} = \text{Fix}\,_\Gam(\fG) \cap \De$.

We denote by $\,\cG\cZ^1\,$ the full subcategory of $\,\cZ^1\,$ whose objects are the g-cocycles.

On the other hand, consider the category of {\em pairs of homomorphic sections}, denoted $\,\cP\cH\cS\,$, having as 
objects the tuples $\,(\Gam, \fE, p, s_1, s_2)\,$ consisting of two profinite groups $\,\Gam\,$ and $\,\fE\,$, 
a splitting epimorphism $\,p : \fE \lra \Gam\,$ and a pair $\,(s_1, s_2)\,$ of homomorphic continuous sections of 
$\,p\,$. A morphism $\,(\Gam, \fE, p, s_1, s_2) \lra (\Gam', \fE', p', s_1', s_2')\,$ is a pair of continuous
homomorphisms $\,(\varphi : \Gam \lra \Gam',\, \psi : \fE \lra \fE')\,$ satisfying 
$\,\varphi \circ p = p' \circ \psi,\, \psi \circ s_i = s_i' \circ \varphi,\,i = 1, 2$. The composition law in 
$\,\cP\cH\cS\,$ is naturally defined. We denote by $\,\cG\cP\cH\cS\,$ the full subcategory of those objects 
$\,(\Gam, \fE, p, s_1, s_2)\,$ of $\,\cP\cH\cS\,$ for which the profinite group $\,\fE\,$ is topologically
generated by the union $\,s_1(\Gam) \cup s_2(\Gam)$. 
\begin{lem}
\label{Lem:EqZ-PHS}
The categories $\,\cZ^1\,$ and $\cP\cH\cS\,$, as well as  
the corresponding full subcategories $\,\cG\cZ^1\,$ and
$\,\cG\cP\cH\cS\,$, are equivalent.
\end{lem}
\bp
We have to extend the equivalence of the categories $\,{\bf PFOG}\,$ and $\,{\bf SEPI}\,$ provided by Lemma 
~\ref{Lem:EqPFOG-SEPI}. To a cocycle $\,\eta \in Z^1(\Gam, \fG)\,$ we first assign the quadruple 
$\,(\Gam, \fE, p, s_1)\,$ associated to the profinite $\,\Gam$-group $\,\fG$, where $\,\fE\,$ is the semidirect 
product $\,\fG \rtimes \Gam\,$ induced by the action of $\,\Gam\,$ on $\,\fG,\,p : \fE \lra \Gam\,$ is the 
natural projection and $\,s_1 : \Gam \lra \fE\,$ is the canonical homomorphic section to $\,p\,$ identifying 
$\,\Gam\,$ with a closed subgroup of $\,\fE\,$ satisfying $\fG \cdot \Gam = \fE,\, \fG \cap \Gam = \{ 1 \}\,$. 
Next we extend the quadruple above by adding the homomorphic section $\,s_2 : \Gam \lra \fE$ to $\,p$ induced 
by the cocycle $\,\eta : s_2(\gam) = \eta(\gam) s_1(\gam)\,$ for all $\gam \in \Gam$. 

To obtain the inverse equivalence, we assign to an object $(\Gam, \fE, p, s_1, s_2)$ of $\cP\cH\cS$ 
the profinite $\,\Gam$-group  $\,\fG = \text{Ker}(p)\,$, with the action of $\,\Gam\,$ defined by 
$\,\gam g := s_1(\gam)\,g\,s_1(\gam)^{- 1}\,$ for $\,\gam \in \Gam, g \in \fG\,$, and the cocycle
$\,\eta \in Z^1(\Gam, \fG)\,$ induced by the homomorphic section $\,s_2 : \eta(\gam) = s_2(\gam) s_1(\gam)^{- 1}\,$ 
for $\,\gam \in \Gam$. 
\ep
Note that the equivalent categories $\,\cZ^1\,$ and $\,\cP\cH\cS\,$, as well as their full subcategories $\,\cG\cZ^1\,$
and $\,\cG\cP\cH\cS\,$, have inverse limits and free products of bundles \cite{Preprint}.

\subsubsection{Universal cocycles}
Given a profinite group $\,\Gam\,$ and a closed subgroup $\,\De \sse \Gam\,$, 
there exists uniquely (up to isomorphism) a pair $\,(\Om_{\Gam, \De},\,\om_{\Gam, \De})\,$ consisting of a profinite
$\,\Gam$-group $\,\Om_{\Gam, \De}\,$ and a g-cocycle $\,\om_{\Gam, \De}: \Gam \lra \Om_{\Gam, \De}\,$, with 
$\,{\rm Fix}_\Gam(\Om_{\Gam, \De}) = \wt{\De} = \bigcap_{\gam \in \Gam} \gam \De \gam^{- 1}$,
$\,\Ker(\om\,_{\Gam, \De}) = \De\,$, such that for all profinite $\,\Gam$-groups $\,\fG\,$, the map
$\,{\rm Hom}_\Gam(\Om_{\Gam, \De}, \fG) \lra Z^1(\Gam\,|\,\De, \fG)\,$, $\,\varphi \mapsto \varphi \circ \om_{\Gam, \De}\,$,
is a functorial bijection \cite{Preprint}. To construct the {\em universal pair}
$\,(\Om_{\Gam, \De},\, \om_{\Gam, \De})\,$, we consider the (generalized) free profinite group   
$\,(\Om_X,\, \om_X : X \lra \Om_X)\,$ generated by the pointed profinite space $\,X:= (\Gam/\De;\,\De)\,$.
Since for any $\,\si \in \Gam\,$, the continuous injective map $\,X \lra \Om_X,\,\tau \De \mapsto \om_X(\si \De)^{- 1}\,
\om_X(\si \tau \De)\,$ extends uniquely to an automorphism, $\,\Om_X\,$ becomes a profinite $\,\Gam$-group denoted
$\,\Om_{\Gam, \De}\,$, with $\,{\rm Fix}_\Gam(\Om_{\Gam, \De}) = \wt{\De}\,$,
while the map $\,\om_{\Gam, \De}: \Gam \lra \Om_{\Gam, \De},\,\gam \mapsto \om_X(\gam \De)\,$,
is a g-cocycle with kernel $\,\De\,$ as desired. In particular, for $\,\De = 1\,$, the action of $\,\Gam\,$ on
$\,\Om_\Gam:= \Om_{\Gam, 1}\,$ is faithful, and the g-cocycle $\,\om_\Gam := \om_{\Gam, 1}\,$ is injective.
As a profinite group, $\,\Om_\Gam\,$ is free of rank
\[
\kappa = \begin{cases}
          |\Gam| - 1 & \mbox{if}\,\, \Gam\, \mbox{is finite},\\
\mbox{max}(\aleph_0, \mbox{rank}\,(\Gam)) & \mbox{if}\,\, \Gam\,\mbox{is\,infinite}.
         \end{cases}
\]

\begin{res} 
\label{Rems:UnivCocycles} \em
(1) The pairs $\,(\Om_{\Gam, \De},\,\om_{\Gam, \De})\,$ are also universal objects for certain {\em embedding problems}
for profinite operator groups \cite{Preprint}.

(2) Let $\,\cC\,$ be a variety of finite groups containing nontrivial groups. For any profinite group $\,\Gam\,$ and 
any closed subgroup $\,\De \sse \Gam\,$, let $\,\Om_{\Gam, \De}^{\,\cC}\,$ be the maximal pro-$\cC$ quotient of 
$\,\Om_{\Gam, \De}\,$, with the induced action of $\,\Gam\,$, and denote by 
$\,\om_{\Gam, \De}^{\,\cC}: \Gam \lra \Om_{\Gam, \De}^{\,\cC}\,$ the g-cocycle induced by $\,\om_{\Gam, \De}\,$. 
Since the variety $\,\cC\,$ contains nontrivial groups, it follows from the construction of $\,\Om_{\Gam, \De}\,$ 
\cite{Preprint} that $\,\Ker(\om_{\Gam, \De}^{\,\cC}) = \De\,$, in particular, 
$\,\om_\Gam^{\,\cC}: \Gam \lra \Om_\Gam^{\,\cC}\,$ 
is injective. The g-cocycles $\,\eta: \Gam \lra \fG\,$, with $\,\De \sse \Ker(\eta)\,$ and $\,\fG\,$ a 
pro-$\cC$ $\,\Gam$-group, are, up to isomorphism, in 1--1 correspondence with the quotients of the
pro-$\cC\,$ $\,\Gam$-group $\,\Om_{\Gam, \De}^{\,\cC}\,$. 

In particular, taking $\,\cC = {\bf FAb}_A\,$, where $\,A\,$ 
is an arbitrary nontrivial subgroup of $\,\Q/\Z\,$, the pro-$\cC\,$ $\,\Gam$-groups are identified with the profinite
(left) $\,A^{\vee}[[\Gam]]\,$-modules, and for any such module $\,\fG\,$, every continuous 
$1$-cocycle $\,\eta : \Gam \lra \fG\,$ extends uniquely to a derivation of the complete group algebra 
$\,A^{\vee}[[\Gam]]\,$ into $\,\fG\,$, i.e., to a continuous $\,A^{\vee}\,$-linear map 
$\,D : A^{\vee}[[\Gam]] \lra \fG\,$ satisfying $\,D(fg) = D(f) \cdot \eps(g) + f \cdot D(g)\,$ for all 
$\,f, g \in A^{\vee}[[\Gam]]\,$, where $\,\eps : A^{\vee}[[\Gam]] \lra A^{\vee}, \gam \in \Gam \mapsto 1\,$, is the 
augmentation map. $\,\Omega_\Gam^{\cC}\,$ becomes the module $\,\Om_{A^{\vee}[[\Gam]]}\,$ of 
(noncommutative) differential forms of $\,A^{\vee}[[\Gam]]\,$, the injective universal cocycle
$\,\om_\Gam^{\cC}: \Gam \lra \Om_\Gam^{\cC}\,$ extends to the universal derivation 
$\,d : A^{\vee}[[\Gam]] \lra \Om_{A^{\vee}[[\Gam]]}, f \mapsto df\,$, and the map 
$\,d\gam \mapsto \gam - 1, \gam \in \Gam\,$, extends to a canonical isomorphism of profinite $\,A^{\vee}[[\Gam]]\,$-
modules from $\Om_{A^{\vee}[[\Gam]]}$ onto the augmentation ideal $\,I\,$ of $\,A^{\vee}[[\Gam]]\,$, the kernel
of the augmentation map $\,\eps\,$. For any closed subgroup $\,\De\,$ of $\,\Gam\,$, the canonical isomorphism 
above induces an isomorphism of profinite $\,A^{\vee}[[\Gam]]$-modules from $\,\Om_{\Gam, \De}^{\cC}\,$ onto the
quotient of the augmentation ideal $\,I\,$ by the left closed ideal $\,J_\De\,$ of $\,A^{\vee}[[\Gam]]\,$ generated by
$\,\{ \del - 1\,|\,\del \in \De \}\,$.  
\end{res}

\subsection{The coGalois connexion associated to a generating cocycle}
\label{Subsec:GeneralFrame}

Fix an object $(\Gam, \fG, \eta)$ of the category $\,\cG\cZ^1$; thus $\,\Gam$ is a profinite group acting 
continuously on the profinite group $\,\fG$, and $\,\eta: \Gam \lra \fG$ is a g-cocycle, so the profinite
group $\,\fG$ is topologically generated by $\,\eta(\Gam)$. Set $\,\De = \text{Ker}(\eta):= \eta^{- 1}(1)$,
$\De' = \text{Fix}_\Gam(\fG)$, and $\wt{\De} = \bigcap_{\gam \in \Gam} \gam \De \gam^{- 1} = \De \cap \De'$.
\begin{de}
\label{Def:Normalization}
The triple $\,(\Gam, \fG, \eta)$ is {\em normalized} (we say also that the g-cocycle $\,\eta$ is 
{\em normalized}) if the closed normal subgroup $\,\wt{\De}$ is trivial. In particular, if $\,\Gam\,$ is abelian
then $\,\eta\,$ is normalized if and only if $\,\eta\,$ is injective.
\end{de}

Given a g-cocycle $\,\eta: \Gam \lra \fG\,$, its {\em normalization} is obtained by replacing the
profinite group $\Gam$ and the cocycle $\eta: \Gam \lra \fG$ with the quotient $\,\Gam':= \Gam/\wt{\De}$ and the
cocycle $\,\eta': \Gam' \lra \fG$ induced by $\eta$, with $\text{Ker}(\eta') = \De/\wt{\De}$, respectively. 

Thus we may assume from the beginning that the triple $\,(\Gam, \fG, \eta)$ is normalized. We associate to 
$\,(\Gam, \fG, \eta)$ two bounded lattices related through natural maps induced by the g-cocycle 
$\,\eta: \Gam \lra \fG$: the lattice $\,\L(\Gam\,|\,\De)\,$ of all closed subgroups of $\,\Gam\,$ lying over $\,\De\,$ 
and the modular lattice $\,\L(\fG)\,$ of all ideals of the profinite $\,\Gam$-group $\,\fG\,$, dual to the lattice of all 
quotients of $\,\fG\,$. According to \ref{Subsec:LatticeSgps} and \ref{Subsubsec:LatticeIdeals}, $\,\L(\Gam\,|\,\De)\,$
and $\,\L(\fG)\,$ are also irreducible spectral spaces with generic points $\,\De\,$ and $\,\{1\}\,$ respectively, 
coherent join operations and continuous meet operations.

The posets $\,\L(\Gam\,|\,\De)\,$ and $\,\L(\fG)\,$ are naturally related through the following canonical 
order-preserving maps induced by the cocycle $\,\eta: \Gam \lra \fG\,$ 
\[
\cJ: \L(\Gam\,|\,\De) \lra \L(\fG), \Lam \mapsto \cJ(\Lam) := \text{the ideal generated by}\,\eta(\Lam),
\]
and
\[
\cS: \L(\fG) \lra \L(\Gam\,|\,\De), {\bf a} \mapsto \cS({\bf a}) := \eta^{- 1}({\bf a}).
\]  
For any subset $\,X \sse \Gam$, let $\,\Lam \in \L(\Gam\,|\,\De)$ be the closed subgroup generated by the union 
$\,X \cup \De$; it follows that $\,\cJ(X) = \cJ(\Lam)$, where $\,\cJ(X)$ is the ideal of the profinite 
$\,\Gam$-group $\,\fG$ generated by $\,\eta(X)$. Note that $\,\fG/\cJ(\Lam)\,$ is the maximal quotient $\,\fH\,$
of the profinite $\,\Gam$-group $\,G\,$ for which the cocycle obtained by composing $\,\eta : \Gam \lra \fG\,$ with
the natural projection $\,\fG \lra \fH\,$ vanishes on $\,\Lam\,$, while, for any $\,{\bf a} \in \L(\fG)$, 
$\,\cS({\bf a}) = \text{Ker}(\eta_{\bf a})\,$, where $\,\eta_{\bf a}: \Gam \lra \fG/{\bf a}$ is the cocycle 
obtained by composing $\,\eta : \Gam \lra \fG\,$ with the natural projection $\,\fG \lra \fG/{\bf a}\,$.
 
\begin{rem} 
\label{Rem:Alternative} \em
To give an alternative description of the operators $\,\cJ\,$ and $\,\cS\,$ above, set $\,\fE := \fG \rtimes \Gam\,$, 
and let $\,\Gam_i \cong \Gam, i = 1, 2,\,$ be the complements of the closed normal subgroup $\,\fG\,$ of $\,\fE\,$ 
induced by the canonical section $\,s_1 : \Gam \lra \fE\,$ and the homomorphic section $\,s_2 : \Gam \lra \fE$ 
determined by the cocycle $\,\eta : \Gam \lra \fG\,$, with $\,\De = \text{Ker}(\eta)$, respectively. As $\,\eta\,$ 
is by assumption a g-cocycle, the profinite group $\,\fE\,$ is topologically generated by the union 
$\,\Gam_1 \cup \Gam_2\,$. For any $\,\Lam \in \L(\Gam\,|\,\De)\,$, let $\,\Lam_i \sse \Gam_i, i = 1, 2,\,$ denote 
the image of $\,\Lam\,$ through the homomorphic sections above. As $\,\De \sse \Lam\,$, we obtain 
$\,\Lam_1 \cap \Lam_2 = \Gam_1 \cap \Gam_2 = s_1(\De) = s_2(\De) \cong \De\,$. Denote by $\,\wt{\fG}\,$ the 
intersection of $\,\fG\,$ with the closed subgroup $\,\wt{\Lam}\,$ of $\,\fE\,$ generated by the union 
$\,\Lam_1 \cup \Lam_2\,$, and note that $\,\wt{\fG}\,$ is topologically generated by $\,\eta(\Lam)\,$, while 
$\,\Lam_1\,$ and $\,\Lam_2\,$ are complements of $\,\wt{\fG}\,$ in $\,\wt{\Lam}\,$. It follows that $\,\cJ(\Lam)\,$ 
is the smallest closed normal subgroup of $\,\fE\,$ containing $\,\wt{\fG}\,$. On the other hand, given 
$\,{\bf a} \in \L(\fG)\,$, let $\,\fH := \fG/{\bf a},\,\fE' := \fH \rtimes \Gam\,$ and $\,\Gam'_i \cong \Gam, 
i = 1, 2,\,$ be the complements of $\,\fH\,$ in $\,\fE'\,$ induced by the canonical section $\,s'_1 : \Gam \lra \fE'\,$ 
and the homomorphic section $\,s'_2\,$ determined by the cocycle $\,\eta': \Gam \lra \fH\,$, obtained by composing 
the cocycle $\,\eta : \Gam \lra \fG\,$ with the projection $\,\fG \lra \fH\,$, respectively. It follows that
$\,\cS({\bf a}) \cong \Gam'_1 \cap \Gam'_2\,$ is the image of $\,\Gam'_1 \cap \Gam'_2\,$ through the projection 
$\,\fE' \lra \Gam$.
\end{rem}

The next result collects together the main properties of the operators $\,\cJ\,$ and $\,\cS\,$ as defined above.
\begin{pr}
\label{Pr:MainProperties} Let $\,\eta \in Z^1(\Gam, \fG)\,$ be a g-cocycle and $\,\De := \Ker(\eta)\,$. 
The following assertions hold. 

{\rm (1)} For all $\,\Lam \in \L(\Gam\,|\,\De), {\bf a} \in \L(\fG)\,$ one has
$\,\Lam \sse \cS(\cJ(\Lam))\,$, $\,\cJ(\cS({\bf a})) \sse {\bf a}\,$, so the pair of operators $\,(\cJ,\,\cS)\,$ 
establishes a coGalois connexion between the posets $\,\L(\Gam\,|\,\De)\,$ and $\,\L(\fG)\,\,$.
 
{\rm (2)} For arbitrary families $\,({\bf a}_i)_{i \in I}\,$ and $\,(\Lam_i)_{i \in I}\,$, with $\,{\bf a}_i \in
\L(\fG),\,\Lam_i \in \L(\Gam\,|\,\De)\,$, one has
\[
\cS(\bigcap_{i \in I} {\bf a}_i) = \bigcap_{i \in I} \cS({\bf a}_i)
\,\,\text{ and }\,\,\cJ(\bigvee_{i \in I} \Lam_i) = \bigvee_{i \in I} \cJ(\Lam_i),
\]
i.e., $\,\cS\,$ and $\,\cJ\,$ are complete semi-lattice morphisms with respect to $\,\cap\,$ and $\,\vee$ respectively. 

{\rm (3)} The map $\,\cJ : \L(\Gam\,|\,\De) \lra \L(\fG)\,$ is coherent, in particular continuous.

{\rm (4)} The map $\,\cS : \L(\fG) \lra \L(\Gam\,|\,\De)$ is continuous.
\end{pr}
\bp
The assertions (1) and (2) are obvious.

(3) Let $\,{\bf b}\,$ be an open ideal of $\,\fG\,$. As the quotient $\,\Gam$-group $\,\fG/{\bf b}\,$ is finite 
and the map $\,\Gam/\cS({\bf b}) \lra \fG/{\bf b}\,$ induced by the cocycle $\,\eta : \Gam \lra \fG\,$ is injective, 
it follows that $\,\cS({\bf b})\,$ is an open subgroup of $\,\Gam\,$ lying over $\,\De\,$. Consequently, the inverse 
image
\[
\{\Lam \in \L(\Gam\,|\,\De)\,|\,\cJ(\Lam) \sse {\bf b}\} =
\{\Lam \in \L(\Gam\,|\,\De)\,|\,\Lam \sse \cS({\bf b})\} = \L(\cS({\bf b})\,|\,\De)
\]
of the basic compact open set $\,\L({\bf b})\,$ of the spectral space $\,\L(\fG)\,$ through the map $\,\cJ\,$ 
is open and compact as desired.

(4) For any open subgroup $\,\Lam \in \L(\Gam\,|\,\De)\,$, let $\,\cW := \{\,{\bf b} \in \L(\fG)\,|\,\cS({\bf b})
\sse \Lam\,\}\,$ denote the inverse image through the map $\,\cS\,$ of the basic open set $\,\L(\Lam\,|\,\De)\,$ of 
the spectral space $\,\L(\Gam\,|\,\De)\,$. As $\,\cS(\{ 1 \}) = \De \sse \Lam,\,\cW\,$ is nonempty. 
For any $\,{\bf b} \in \cW\,$, denote by $\,\L(\fG\,|\,{\bf b})_{\text{o}}\,$ the poset of all {\em open} ideals of 
$\,\fG\,$ lying over $\,{\bf b}\,$, so $\,{\bf b} = \bigcap_{{\bf a} \in \L(\fG\,|\,{\bf b})_{\text{o}}} {\bf a}\,$. 
By (2) and by the compactness of $\,\Gam\,$ it follows that there exists $\,{\bf a} \in \L(\fG\,|\,{\bf b})_{\text{o}}$ 
such that $\,\cS({\bf b}) \sse \cS({\bf a}) \sse \Lam\,$. Consequently, the nonempty set $\,\cW_{\rm max}\,$ of all 
maximal members of $\,\cW\,$ with respect to inclusion consists of open ideals of the profinite $\,\Gam$-group $\,\fG\,$, 
and hence $\,\cW = \bigcup_{{\bf c} \in \cW_{\rm max}}\L({\bf c})\,$ is open as a union of basic open sets of 
the spectral space $\,\L(\fG)\,$.
\ep

\begin{re} 
\label{Re:CoherentS} \em
By contrast with the map $\,\cJ\,$ which is always coherent, the continuous map $\,\cS\,$ is not coherent in general. 
For instance, let $\,\Gam = (\wh{\Z},\,+)\,$ and $\,\fG = \dprod_{p \in \cP'} (\Z/p \Z,\,+)\,$, where $\cP'$ is the 
set of the odd prime numbers $\,p\,$ for which the order $\,f_p\,|\,(p - 1)\,$ of $\,2\,\text{mod}\,p\,$ is even. 
Consider the continuous action $\,\Gam \times \fG \lra \fG,\,(\gam, g) \mapsto 2^\gam g$ and the coboundary 
$\,\eta: \Gam \lra \fG,\,\gam \mapsto 2^\gam - 1\,$; since $\eta$ sends the topological generator $\,1\,$ of
$\,\Gam\,$ to a topological generator of $\,\fG\,$, $\,\eta$ is a g-cocycle. Note that 
$\,\De:= \text{Fix}_\Gam(\fG) = \Ker(\eta) = \dbigcap_{p \in \cP'} f_p \wh{\Z}$. For the open subgroup 
$\,\Lam:= 2 \wh{\Z} \in \L(\Gam\,|\,\De)\,$, we obtain, with the notation from Proposition~\ref{Pr:MainProperties}, (4),
$\,\cW_{\text{max}} = \{\dprod_{p \in \cP' \sm \{l\}} \Z/p \Z\,|\,l \in \cP'\}\,$, the set of all maximal open
subgroups of $\,\fG\,$, whence the open set $\,\cW\,$ is not compact since the set $\cP'\,$ is infinite (any odd prime
number $\,p \not \equiv \pm 1\,\text{mod}\,8$ belongs to $\,\cP'\,$).  

Thus, given g-cocycles $\,\eta: \Gam \lra \fG\,$, it is natural to look for suitable closed subspaces 
$\,X\,$ of the spectral space $\,\L(\fG)\,$ for which the restriction $\,\cS|_X: X \lra \L(\Gam\,|\,\De)\,\,$ 
becomes coherent. 
\end{re}
Other useful properties of g-cocycles are collected in the next lemma whose proof is straightforward.
\begin{lem}
\label{Lem:PropGenCocycles}
Let $\,\eta \in Z^1(\Gam, \fG)\,$ be a g-cocycle with $\,\De := \Ker(\eta),\,\De':= {\rm Fix}\,_\Gam(\fG)\,$,
$\,\wt{\De}:= \cap_{\gam \in \Gam} \gam \De \gam^{- 1} = \De \cap \De'\,$, $\,\De'':= 
\{\si \in \Gam\,|\,\forall \gam \in \Gam,\,\gam\,\eta(\si) = \eta(\gam \si \gam^{- 1})\},\,\overline{\De}:= 
\De' \cap \De''\,$. The following assertions hold.

{\rm (1)} $\,\De''\,$ is the maximal closed normal subgroup $\,\Lam \sse \Gam\,$ for which the restrictiom map
$\,\eta\,|_\Lam: \Lam \lra \fG\,$ is $\,\Gam$-equivariant; in particular, $\,\overline{\De}\,$ is a closed normal 
subgroup of $\,\Gam$ containing $\,\wt{\De}$. 

{\rm (2)} If $\,\fG\,$ is abelian, then $\,\De' \sse \De''\,$, so $\,\overline{\De} = \De'$.

{\rm (3)} $\,\eta(\overline{\De}) = \cJ(\overline{\De})\,$ is contained in the center $\,C(\fG)\,$ of $\,\fG$.

{\rm (4)} $\,\overline{\De} = \De' \cap \cS(C(\fG))$.

{\rm (5)} $\,\eta\,$ induces by restriction an isomorphism of profinite $\,\Gam$-modules 
$\,\overline{\De}/\wt{\De} \cong \eta(\overline{\De})$.
\end{lem}

\section{Continuous actions on discrete abelian groups}
\label{Sec:Modules}

In this section we apply the general framework provided by \ref{Subsec:GeneralFrame} to the case of continuous
actions of profinite groups on discrete abelian groups, including as particular cases among others the framework of the 
classical coGalois theory of separable radical extensions as well as  its abstract version developed in 
\cite{A-B, Serdica, JPAA}.
\subsection{The coGalois group of a discrete module}
\label{Subsec:coG}
Let $\,\Gam\,$ be a profinite group, $\,E$ a discrete $\,\Gam$-module, and $\,A:= t(E)\,$ the torsion group of the 
abelian group $\,E\,$ with the induced action of $\,\Gam\,$. In the following we extend to arbitrary discrete
$\,\Gam$-modules $\,E\,$ the notion of {\em coGalois group} $\,\text{coG}(L/K)\,$ of a field extension $\,L/K\,$, 
introduced by Greither and Harrison \cite{G-H}, which plays a major role in the {\em coGalois theory of radical
extensions}.  

Setting $\,H^0(\Gam, E) = E^\Gam\ = \{x \in E\,|\,\gam x = x\,\text{for all}\,\gam \in \Gam\}\,$, we denote 
\[
\text{Rad}(E\,|\,E^\Gam):= \{x \in E\,|\,n x \in E^\Gam\,\text{for some}\,n \in \N \sm \{0\}\,\}. 
\]
$\text{Rad}(E\,|\,E^\Gam)$ is a $\,\Gam$-submodule of $\,E\,$, $\,\text{Rad}(E\,|\,E^\Gam)^\Gam = E^\Gam\,$, and 
$\,t(\text{Rad}(E\,|\,E^\Gam)) = t(E) = A\,$. By the {\em coGalois group} of the discrete $\,\Gam$-module $\,E\,$, 
denoted by $\,\text{coG}(E)$, we understand the quotient group $\,\text{Rad}(E\,|\,E^\Gam)/E^\Gam\,$, which is 
nothing else than the torsion group $\,t(E/E^\Gam)\,$ of the quotient group $\,E/E^\Gam\,$. Note that 
$\,\text{coG}(E)\,$ is a discrete $\,\Gam$-module, and $\,\text{coG}(A) = A/A^\Gam$ is a $\,\Gam$-submodule of 
$\,\text{coG}(E)\,$. Since $\,\text{coG}(E) = \text{coG}(\text{Rad}(E\,|\,E^\Gam))$, in order to study the coGalois
group of a discrete $\,\Gam$-module $\,E\,$ we may assume without loss that $\,E = \text{Rad}(E\,|\,E^\Gam)$, so
$\,\text{coG}(E) = E/E^\Gam\,$.

The profinite group $\,\Gam\,$ and the discrete abelian group $\,E = \text{Rad}(E\,|\,E^\Gam)\,$ are naturally 
related through the map $\Gam \times E \lra A,\,(\gam, x) \mapsto \gam x - x\,$, which induces a map \\
$\,\Gam/\text{Fix}_\Gam(E) \times \text{coG}(E) \lra A$ relating the profinite quotient
group $\,\Gam/\text{Fix}_\Gam(E)\,$ and the discrete abelian torsion group $\,\text{coG}(E)\,$, so we may also
assume that $\,\text{Fix}_\Gam(E)\,$ is trivial, i.e., the action of $\,\Gam\,$ on $\,E\,$ is faithful. 

Consider the homomorphism $\,\theta: E \lra Z^1(\Gam, A)$ defined by $\,\theta(x)(\gam) = \gam x - x\,$ for
$x \in E, \gam \in \Gam\,$, where $\,Z^1(\Gam, A)\,$ is the discrete torsion abelian group of all continuous 
$\,1$-cocycles of $\,\Gam\,$ with coefficients in $A$. As $\,\Ker(\theta) = E^\Gam\,$, $\,\text{coG}(E) \cong
\theta(E)\,$ is identified with a subgroup of $\,Z^1(\Gam, A)\,$, while $\,\text{coG}(A) \cong \theta(A) = 
B^1(\Gam, A)\,$, the subgroup of $\,1$-coboundaries \\$\,f_a: \Gam \lra A,\,\gam \mapsto \gam a - a\,$, for $\,a \in A\,$.
Consequently, the quotient group $\,\text{coG}(E)/\text{coG}(A)$ is identified with a subgroup of 
$\,H^1(\Gam, A):= Z^1(\Gam, A)/B^1(\Gam, A)\,$.

Using the exact sequence of cohomology groups in low dimensions associated to the short exact sequence of discrete 
$\,\Gam$-modules $\,0 \lra A \lra E \lra E/A \lra 0\,$, we obtain 
\begin{lem}
\label{Lem:H1}
Let $\,E\,$ be an arbitrary discrete $\,\Gam$-module, with $\,{\rm Rad}(E\,|\,E^\Gam)\,$ may be properly
contained in $\,E\,$, and $\,A:= t(E)\,$. Then the following assertions are equivalent.

{\rm (1)} $\,{\rm coG}(E) \cong Z^1(\Gam, A)$.

{\rm (2)} $\,{\rm coG}(E)/{\rm coG}(A) \cong H^1(\Gam, A)$.

{\rm (3)} $\,H^1(\Gam, E) = 0$.

{\rm (4)} $\,H^1(\Gam, {\rm Rad}(E\,|\,E^\Gam)) = 0$.
\end{lem}

\begin{rem}
\label{Rem:ActOnZ} \em
Consider the canonical continuous action of $\,\Gam\,$ on $\,Z^1(\Gam, A)\,$ defined by
\[ 
(\si \alpha)(\gam) = \si \alpha(\si^{- 1} \gam \si) = \alpha(\gam) + (\gam \alpha(\si) - \alpha(\si))\,\text{for}
\,\alpha \in Z^1(\Gam, A), \si, \gam \in \Gam\,, 
\]
i.e., $\,\si \alpha = \alpha + f_{\alpha(\si)}\,$. Thus the homomorphism $\,\theta: E \lra Z^1(\Gam, A)\,$ is 
$\,\Gam$-equivariant, and hence $\,\text{coG}(E)$ is identified with a $\,\Gam$-submodule of $\,Z^1(\Gam, A)\,$. 
In particular, if $\,\text{coG}(A) = 0\,$, i.e., $A^\Gam = A\,$, then $\,\theta(\si x) = \theta(x)\,$ for all 
$\,\si \in \Gam, x \in E\,$, and $\,\text{coG}(E)\,$ is identified with a subgroup of the torsion abelian group 
$\,\text{Hom}(\Gam, A)\,$ of continuous homomorphisms from $\,\Gam\,$ to $\,A\,$.
\end{rem}
   
\subsection{The Galois connexion associated to a discrete torsion module}
\label{Subsec:GalCon}

Let $\,\Gam$ be a profinite group acting continuously on a discrete torsion abelian group $\,A\,$.
Consider the {\em evaluation map} $\,\Gam \times Z^1(\Gam, A) \lra A,\,(\gam, \alpha) \mapsto \alpha(\gam)\,$,
relating the profinite group $\,\Gam\,$ and the discrete torsion $\,\Gam$-module $\,Z^1(\Gam, A)\,$, with the
action of $\,\Gam\,$ defined as in Remark~\ref{Rem:ActOnZ}.

Consider also the lattice $\,\L(\Gam)\,$, its modular sublattice $\,\L_{\text{n}}(\Gam)\,$ of all closed normal subgroups
of $\,\Gam\,$, the modular lattice $\,\L(\fZ)\,$ of all subgroups of $\,\fZ\,$, and its sublattice $\,\L_\Gam(\fZ)\,$ of
all $\,\Gam$-submodules of $\,\fZ\,$. Note that $\,\L(\fZ)\,$ is also an irreducible spectral space with the basic 
compact sets $\,\L(\fZ\,|\,F):= \{G \in \L(\fZ)\,|\,F \sse G\}\,$ for $\,F\,$ ranging over all finite subgroups 
of $\,\fZ$, generic point $\,\fZ$, and the unique closed point $\,\{0\}\,$. The action 
$\,\Gam \times \L(\fZ) \lra \L(\fZ),\,(\gam, G) \mapsto \gam G = \{\gam \alpha\,|\,\alpha \in G\}\,$
and the meet operation $\,\L(\fZ) \times \L(\fZ) \lra \L(\fZ),\,(G_1, G_2) \mapsto G_1 \cap G_2\,$ are coherent maps,
while the join operation $\,\L(\fZ) \times \L(\fZ) \lra \L(\fZ),\,(G_1, G_2) \mapsto G_1 + G_2\,$ is a continuous map,
not necessarily coherent. 
 
The posets $\,\L(\Gam)\,$ and $\,\L(\fZ)\,$ are related through the canonical order-reversing maps
\[
\L(\Gam) \lra \L(\fZ),\,\Lam \mapsto \Lam^\p:= Z^1(\Gam\,|\,\Lam, A) = \{\alpha \in \fZ\,|\,\alpha|_\Lam = 0\}
\]
and
\[
\L(\fZ) \lra \L(\Gam),\,G \mapsto G^\p:= \bigcap_{\alpha \in G} \alpha^\p,  
\]
where $\,\alpha^\p:= \Ker(\alpha) = \{\gam \in \Gam\,|\,\alpha(\gam) = 0\}\,$ is an open subgroup of $\,\Gam\,$ for
all $\alpha \in \fZ\,$. For any $\,\Lam \in \L(\Gam)\,$ one denotes by $\,\text{res}^\Gam_\Lam: \fZ \lra Z^1(\Lam, A),
\,\alpha \mapsto \alpha|_\Lam\,$, the restriction homomorphism. It follows that $\,\Lam^\p = \Ker(\text{res}^\Gam_\Lam)\,$
and $\,(\text{res}^\Gam_\Lam(G))^\p = G^\p \cap \Lam\,$ for all $\,\Lam \in \L(\Gam), G \in \L(\fZ)\,$.
                                                                           
In the following we assume without loss that the closed normal subgroup $\,\fZ^\p\,$ of $\,\Gam\,$ is trivial, 
and hence the closed normal subgroup $\,B^1(\Gam, A)^\p = \text{Fix}_\Gam(A) \sse \Gam\,$ is abelian since
$\,\alpha(\si \tau) = \alpha(\si) + \alpha(\tau) = \alpha(\tau \si)\,$ for 
$\alpha \in \fZ, \si, \tau \in \text{Fix}_\Gam(A)$, while the profinite quotient group $\,\Gam/\text{Fix}_\Gam(A)\,$
is identified with a closed subgroup of the totally disconnected topological group $\,\text{Aut}(A)$ for which
the subgroups $\,\text{Aut}(A\,|\,F):= \{\varphi \in \text{Aut}(A)\,|\,\varphi|\,_F = 1_F\}$ for $\,F\,$ ranging over 
the finite subgroups of $\,A\,$ form a fundamental system of open neighborhoods of $\,1_A\,$.
  
The next result is an analogue of Proposition~\ref{Pr:MainProperties}, and the proof is similar.
\begin{pr}
\label{Pr:GalCon1} The following assertions hold.

{\rm (1)} The pair of order-reversing maps
\[
\L(\Gam) \lra \L(\fZ),\,\Lam \mapsto \Lam^\p,\,\L(\fZ) \lra \L(\Gam),\,G \mapsto G^\p,
\]
establishes a Galois connexion between the posets $\,\L(\Gam)\,$ and $\,\L(\fZ)\,$, i.e., $\,X \sse X^{\p\,\p}\,$
for any element $\,X\,$ of $\,\L(\Gam)\,$ or $\,\L(\fZ)\,$.
 
{\rm (2)} The map $\,\Lam \mapsto \Lam^\p\,$ is a coherent complete-semi-lattice morphism \\
$\,(\L(\Gam), \vee) \lra (\L(\fZ), \wedge)\,$ satisfying $\,(\si \Lam \si^{- 1})^\p = \si \cdot \Lam^\p\,$ for all
$\,\si \in \Gam\,$, $\,\Lam \in \L(\Gam)$.
 
{\rm (3)} The map $\,G \mapsto G^\p\,$ is a continuous, not necessarily coherent, complete-semi-lattice morphism 
$\,(\L(\fZ), \vee) \lra (\L(\Gam), \wedge)\,$ satisfying $\,(\si \cdot G)^\p = \si \cdot G^\p \cdot \si^{- 1}\,$ 
for all $\,\si \in \Gam\,$, $\,G \in \L(\fZ)$.

{\rm (4)} The maps above induce by restriction a Galois connexion between the posets $\,\L_n(\Gam)\,$ and
$\,\L_\Gam(\fZ)\,$.
\end{pr}

\begin{rem} 
\label{Rem:GalConCogE} \em
Let $\,E\,$ be a discrete $\,\Gam$-module such that $\,t(E) = A\,$, and assume without loss that 
$\,\text{coG}(E):= t(E/E^\Gam) = E/E^\Gam\,$ and $\,\text{Fix}_\Gam(E)\,$ is trivial. According to 
Remark~\ref{Rem:ActOnZ}, $\,\text{coG}(E)\,$ is identified with a $\,\Gam$-submodule of $\,\fZ\,$, and hence 
$\,\fZ^\p \sse \text{coG}(E)^\p = \text{Fix}_\Gam(E) = \{1\}\,$, so $\,\fZ^\p$ is trivial. The assertions of 
Propsition~\ref{Pr:GalCon1} remain valid for the order-reversing maps 
\[
\L(\Gam) \lra \L(\,\text{coG}(E)),\,\Lam \mapsto \Lam^\p \cap \text{coG}(E),\,\,\L(\text{coG}(E)) \lra 
\L(\Gam),\,G \mapsto G^\p.
\] 
\end{rem}

\subsection{The associated profinite module and generating cocycle}
\label{Subsec:AssProfMod}
Consider the same data as in \ref{Subsec:GalCon}. We construct a profinite $\,\Gam$-module $\,\fG\,$, a 
sort of dual of the discrete torsion $\,\Gam$-module $\,\fZ\,$ with respect to the discrete torsion $\,\Gam$-module 
$\,A\,$, and a natural continuous g-cocycle $\,\eta: \Gam \lra \fG\,$. To this end, we consider the abelian
group $\fH:= \text{Hom}(\fZ, A)\,$, and for any subgroup $\,G \sse \fZ\,$, we denote by $\,r_G\,$ the restriction 
homomorphism $\,\fH \lra \text{Hom}(G, A),\,\varphi \mapsto \varphi|_G\,$, with kernel $\,\fH_G$. $\,\fH\,$ becomes 
a totally disconnected (Hausdorff) topological group with the compact-open topology, for which the subgroups 
$\,\fH_F\,$, with $\,F\,$ ranging over the {\em finite} subgroups of $\,\fZ\,$, serve as a fundamental system of 
open neighborhoods of the null homomorphism. The profinite group $\,\Gam\,$ acts continuously on $\,\fH\,$ 
according to the rule $\,(\gam \varphi)(\alpha):= \gam \varphi(\alpha)\,$ for 
$\,\gam \in \Gam\,$, $\varphi \in \fH\,$, $\alpha \in \fZ\,$, and the canonical map 
$\,\eta: \Gam \lra \fH,\,\gam \mapsto \eta_\gam\,$, defined by $\,\eta_\gam(\alpha):= \alpha(\gam)\,$ for 
$\,\gam \in \Gam, \alpha \in \fZ\,$, is a continuous $\,1$-cocycle; note that, with the notation from 
\ref{Subsec:GalCon}, $\,\eta^{- 1}(\fH_F) = F^\p = \dbigcap_{\alpha \in F} \alpha^\p\,$ is an open subgroup of 
$\,\Gam\,$ provided $\,F\,$ is a finite subgroup of $\,\fZ\,$. As $\,\Ker(\eta) = \fZ^\p$ is a closed normal 
subgroup of $\,\Gam\,$, we may assume without loss that the cocycle $\,\eta\,$ is injective.

Next let us consider the subgroup $\,\la \eta(\Gam) \ra \sse \fH\,$ generated by $\,\eta(\Gam)\,$; for any finite
subgroup $\,F \sse \fZ\,$, let $\,\fG^F:= r_F(\la \eta(\Gam) \ra) = \la \eta_\gam|_F\,:\,\gam \in \Gam \ra \cong 
\la \eta(\Gam \ra/(\la \eta(\Gam) \ra \cap \fH_F)\,$. Since any $\,\alpha \in \fZ\,$ is a continuous map, and hence
locally constant, with values in the torsion abelian group $\,A\,$, we deduce that $\,\fG^F \sse \text{Hom}(F, A)\,$ 
is a finite abelian group for all finite subgroups $\,F \sse \fZ\,$. Since the induced topology on the subgroup 
$\,\la \eta(\Gam) \ra \sse \fH\,$ is determined by the collection of the subgroups $\,\la \eta(\Gam) \ra \cap \fH_F\,$
of finite index filtered from below, it follows that the closure $\,\fG:= \overline{\la \eta(\Gam) \ra}$ is a
profinite abelian group isomorphic with $\displaystyle\lim_{\stackrel{\longleftarrow}{F}}\,\fG^F\,$, topologically
generated by $\,\eta(\Gam)\,$. Moreover $\,\fG\,$ is a profinite $\,\Gam$-module with the induced action of $\,\Gam\,$,
and $\,\eta: \Gam \lra \fG\,$ is an injective g-cocycle.

\begin{res}
\label{Rem:2ndActOnG} \em
(1) Using \ref{Rems:UnivCocycles} (2), the profinite $\,\Gam$-module $\,\fG\,$ is the quotient of 
$\,\Om_\Gam^{{\rm ab}}\,$, identified with the augmentation ideal of $\,\wh{\Z}[[\Gam]]\,$, by the closed submodule 
$\,\cap_{\varphi \in {\rm Hom}_\Gam(\Om_\Gam^{{\rm ab}},\,A)} \Ker(\varphi)\,$, while the g-cocycle 
$\eta: \Gam \lra \fG$ is induced by $\,\om_\Gam^{{\rm ab}}: \Gam \lra \Om_\Gam^{{\rm ab}}$. See also 
Proposition~\ref{Pr:Pairing}. 

(2) Aside from the continuous action $\,\Gam \times \fG \lra \fG,\,(\si, \varphi) \mapsto \si \varphi\,$, considered 
above, there is another continuous action of the profinite group $\,\Gam\,$ on the profinite abelian group 
$\,\fG \sse \text{Hom}(\fZ, A)\,$ induced by the continuous actions of $\,\Gam\,$ on the discrete torsion 
abelian groups $\,\fZ\,$ and $\,A\,$, defined by 
\[
({}^\si \varphi)(\alpha):= \si \cdot \varphi(\si^{- 1} \alpha) = (\si \varphi)(\alpha) + (\si \varphi)(f_{\alpha(\si^{- 1})})
\] 
for $\,\si \in \Gam, \varphi \in \fG, \alpha \in \fZ\,$.
In particular, taking $\,\varphi = \eta_\gam\,$ for $\,\gam \in \Gam\,$, we obtain $\,{}^\si \eta_\gam = 
\eta_{\si \gam \si^{- 1}}\,$ for all $\,\si \in \Gam\,$, whence the injective continuous map $\,\eta: \Gam \lra \fG,
\gam \mapsto \eta_\gam\,$, is $\,\Gam$-equivariant with respect to the action of $\,\Gam\,$ on itself by inner
automorphisms and the afore defined action $\,\Gam \times \fG \lra \fG,\,(\si, \varphi) \mapsto {}^\si \varphi\,$. 
\end{res}
We denote by $\,\L(\fG)\,$ the modular lattice of all closed $\,\Gam$-submodules of $\,\fG\,$ (with respect to the 
action $\,\Gam \times \fG \lra \fG,\,(\si, \varphi) \mapsto \si \varphi\,$), and by $\,\L_\Gam(\fG)\,$ its sublattice
consisting of those $\,\Gam$-submodules which are also invariant under the action 
$\,\Gam \times \fG \lra \fG,\,(\si, \varphi) \mapsto {}^\si \varphi\,$. 

According to \ref{Subsubsec:LatticeIdeals},
$\,\L(\fG)\,$ is also an irreducible spectral space with basic compact open sets 
$\,\L({\bf a}):= \{{\bf b} \in \L(\fG)\,|\,{\bf b} \sse {\bf a}\}\,$ for all open $\,{\bf a} \in \L(\fG)\,$, generic 
point $\,\{0\}\,$, unique closed point $\,\fG\,$, coherent join operation ($\,+\,$), and continuous meet operation 
($\,\cap\,$).

Applying Proposition~\ref{Pr:MainProperties} to the triple $\,(\Gam, \fG, \eta)\,$ above, we obtain
\begin{pr}
\label{Pr:coGCon}
The following assertions hold.

{\rm (1)} The order-preserving maps
\[
\cJ: \L(\Gam) \lra \L(\fG),\,\Lam \mapsto\,\mbox{the\,closed}\,\,\Gam-\mbox{submodule\, generated\, by}\,\,\eta(\Lam),
\]
\[
\cS: \L(\fG) \lra \L(\Gam),\,{\bf a} \mapsto \eta^{- 1}({\bf a}),
\]
establishes a coGalois connexion between the posets $\,\L(\Gam)\,$ and $\,\L(\fG)\,$, i.e., 
$\,\Lam \sse \cS(\cJ(\Lam))\,$ and $\,\cJ(\cS({\bf a})) \sse {\bf a}\,$ for all 
$\,\Lam \in \L(\Gam),\,{\bf a} \in \L(\fG)$.

{\rm (2)} The coherent map $\,\Lam \mapsto \cJ(\Lam)\,$ is a morphism of complete semilattices \\
$\,(\L(\Gam),\,\vee) \lra (\L(\fG),\,+)\,$, satisfying $\,\cJ(\si \Lam \si^{- 1}) = {}^\si \cJ(\Lam)$ for all 
$\,\si \in \Gam,\,\Lam \in \L(\Gam)$.

{\rm (3)} The continuous (not necessarily coherent) map $\,{\bf a} \mapsto \cS({\bf a})\,$ is a morphism of 
complete semilattices $\,(\L(\fG),\,\cap) \lra (\L(\Gam),\,\cap)\,$, satisfying $\,\cS({}^\si {\bf a}) =
\si \cS({\bf a}) \si^{- 1}\,$ for all $\,\si \in \Gam, {\bf a} \in \L(\fG)$.

{\rm (4)} The order-preserving maps $\,\cJ\,$ and $\,\cS\,$ induce by restriction a coGalois connexion between
the posets $\,\L_{\rm n}(\Gam)\,$ and $\,\L_\Gam(\fG)$. 
\end{pr}
   
\subsection{A nondegenerate pairing and the induced Galois connexion}
\label{Subsec:Pairing}
By Propositions~\ref{Pr:GalCon1}, \ref{Pr:coGCon}, the lattice $\,\L(\Gam)\,$ is related to the lattices $\,\L(\fZ)\,$
and $\,\L(\fG)\,$ through canonical maps defining a Galois connexion and a coGalois connection respectively. It is 
also natural to consider the relation between the lattices $\,\L(\fZ)\,$ and $\,\L(\fG)\,$.

First, with data and notation from \ref{Subsec:AssProfMod}, we obtain
\begin{pr}
\label{Pr:Pairing}
There is a canonical nondegenerate pairing 
\[
\la\,,\,\ra: \fG \times \fZ \lra A,\,(\varphi, \alpha) \mapsto \la \varphi, \alpha \ra:= \varphi(\alpha)
\]
satisfying the identities 
\[
\la \si \varphi, \alpha \ra = \si \la \varphi, \alpha \ra = \la {}^\si \varphi, \si \alpha \ra
\]
for all $\,\si \in \Gam, \varphi \in \fG, \alpha \in \fZ\,$.

The pairing above induces a canonical isomorphism of discrete torsion $\,\Gam$-modules \\
$\,\lam: \fZ \lra {\rm Hom}_\Gam(\fG, A)\,$, defined by $\,\lam(\alpha)(\varphi):= \varphi(\alpha)\,$, whose 
inverse \\ $\,\mu: {\rm Hom}_\Gam(\fG, A) \lra \fZ\,$ is defined by $\,\mu(\psi):= \psi \circ \eta\,$; here 
$\,{\rm Hom}_\Gam(\fG, A)\,$ denotes the discrete torsion abelian group of all continuous homomorphisms 
$\,\psi: \fG \lra A\,$ satisfying $\,\psi(\si \varphi) = \si \psi(\varphi)\,$ for all 
$\,\si \in \Gam, \varphi \in \fG\,$, together with the continuous action of $\,\Gam\,$ defined by 
$\,({}^\si \psi)(\varphi):= \si \cdot \psi({}^{\si^{- 1}} \varphi) = \psi(\si \cdot ({}^{\si^{- 1}} \varphi))\,$.
\end{pr}
\bp
The statement follows easily from the definitions given in \ref{Subsec:AssProfMod} of $\,\fG\,$ and of 
the actions $\,\Gam \times \fG \lra \fG, (\si, \varphi) \mapsto \si \varphi,\,{}^\si \varphi\,$. To prove that 
$\,\lam\,$ and $\,\mu\,$ are isomorphisms inverse to one another we have to show that $\,\mu \circ \lam = 1_{\fZ}\,$ 
and $\,\mu\,$ is injective. The first condition is satisfied since 
$\,(\mu \circ \lam)(\alpha) = \lam(\alpha) \circ \eta\,$ for all $\,\alpha \in \fZ\,$, and 
$\,(\lam(\alpha) \circ \eta)(\gam) = \lam(\alpha)(\eta_\gam) = \eta_\gam(\alpha) = \alpha(\gam)\,$ for all 
$\,\gam \in \Gam\,$. To check the injectivity of $\,\mu\,$, let $\,\psi \in \text{Hom}_\Gam(\fG, A)\,$ be such 
that $\,\mu(\psi) = \psi \circ \eta = 0\,$. As $\,\fG\,$ is topologically generated by the set $\,\eta(\Gam)\,$, 
it follows that $\,\psi = 0\,$ as desired.
\ep
Next we note that the posets $\,\L(\fG)\,$ and $\,\L(\fZ)\,$ are related through the canonical order-reversing maps
\[
\L(\fG) \lra \L(\fZ),\,{\bf a}\,\mapsto {\bf a}_\p:= \bigcap_{\varphi \in {\bf a}} \varphi_\p,
\]
where $\,\varphi_\p:= \Ker(\varphi)\,$ for $\,\varphi \in \fG \sse \text{Hom}(\fZ, A)\,$, and
\[
\L(\fZ) \lra \L(\fG),\,G \mapsto G_\p:= \bigcap_{\alpha \in G} \alpha_\p,
\]
where $\,\alpha_\p:= \Ker(\lam(\alpha)) = \{\varphi \in \fG\,|\,\varphi(\alpha) = 0\}$.

The following result is an analogue of Proposition~\ref{Pr:GalCon1}.
\begin{pr}
\label{Pr:GalCon2}
The following assertions hold.

{\rm (1)} The pair of order-reversing maps
\[
\L(\fG) \lra \L(\fZ),\,{\bf a} \mapsto {\bf a}_\p,\,\L(\fZ) \lra \L(\fG),\,G \mapsto G_\p,
\]
establishes a Galois connexion, i.e., $\,X \sse X_{\p\,\p}\,$ for any element $\,X\,$ of $\,\L(\fG)\,$ or $\,\L(\fZ)\,$.

{\rm (2)} The coherent map $\,{\bf a} \mapsto {\bf a}_\p\,$ is a morphism of complete semilattices \\  
$\,(\L(\fG),\, +) \lra (\L(\fZ),\,\cap)\,$, satisfying $\,({}^\si {\bf a})_\p = \si \cdot ({\bf a}_\p)\,$ for all 
$\,\si \in \Gam, {\rm a} \in \L(\fG)$. 

{\rm (3)} The continuous (not necessarily coherent) map $\,G \mapsto G_\p\,$ is a morphism of complete semilattices 
$\,(\L(\fZ),\,+) \lra (\L(\fG),\,\cap)\,$, satisfying $\,(\si G)_\p = {}^\si (G_\p)\,$ for all $\,\si \in \Gam\,$, 
$\,G \in \L(\fZ)$.

{\rm (4)} The maps above induce by restriction a Galois connexion between the posets $\,\L_\Gam(\fG)\,$ and
$\,\L_\Gam(\fZ)$.
\end{pr}
\begin{co}
\label{Cor:Composing}
The Galois connexion between the lattices $\,\L(\Gam)\,$ and $\,\L(\fZ)\,$ ({\em Proposition~\ref{Pr:GalCon1}})
is obtained by composing the coGalois connexion between the lattices $\,\L(\Gam)\,$ and $\,\L(\fG)\,$ 
({\em Proposition~\ref{Pr:coGCon}}) with the Galois connexion between the lattices $\,\L(\fG)\,$ and $\,\L(\fZ)\,$
({\em Proposition~\ref{Pr:GalCon2}}): $\,\Lam^\p = \cJ(\Lam)_\p,\,G^\p = \cS(G_\p)\,$ for 
$\,\Lam \in \L(\Gam),\,G \in \L(\fZ)$.
\end{co}

\subsection{Examples}
\label{Subsec:Examples} 

\begin{ex}
\label{Ex:1} \em
({\em The cyclotomic abstract coGalois theory} \cite{A-B, Serdica, JPAA}) Let $\,\Gam\,$ be a profinite group, $\,A\,$
a discrete quasi-cyclic group identified with a subgroup of $\,\Q/\Z\,$, $\,A^\vee:= \text{Hom}(A, \Q/\Z)\,$, its
Pontryagin dual, and $\,\Gam \times A \lra A,\,(\si, a) \mapsto \si a:= \chi(\si) a\,$, a continuous action given by
the continuous {\em cyclotomic character} $\,\chi: \Gam \lra (A^\vee)^\times\,$. With the notation from 
\ref{Subsec:AssProfMod}, we obtain $\,\fG = \fH = \text{Hom}(\fZ, A) = \text{Hom}(\fZ, \Q/\Z) = \fZ^\vee\,$, the
Pontryagin dual of the discrete torsion abelian group $\,\fZ\,$. $\,\fG\,$ is a profinite $\,\Gam$-module with respect
to the continuous action $\,\Gam \times \fG \lra \fG,\,(\si, \varphi) \mapsto \chi(\si) \varphi\,$, and the canonical map
$\,\eta: \Gam \lra \fG,\,\gam \mapsto (\eta_\gam: \fZ \lra A)\,$, defined by $\,\eta_\gam(\alpha):= \alpha(\gam)\,$,
is a continuous g-cocycle. Assuming without loss that $\eta\,$ is injective, i.e., the closed normal
subgroup $\,\fZ^\p\,$ is trivial, it follows that $\,\Gam\,$ is metabelian as an extension of the abelian profinite
quotient group $\,\Gam/\text{Fix}_\Gam(A) \cong \chi(\Gam) \sse (A^\vee)^\times\,$ by the abelian profinite group
$\,B^1(\Gam, A)^\p = \text{Fix}_\Gam(A)\,$. It follows also that any closed subgroup of $\,\fG\,$ is invariant under
the action of $\,\Gam\,$, any continuous homomorphism $\,\psi: \fG \lra A\,$ is $\,\Gam$-equivariant, and, by the
Pontryagin duality, the Galois connexion between the lattices $\,\L(\fZ)\,$ and $\,\L(\fG)\,$ is perfect, i.e.,
the maps $\,\L(\fZ) \lra \L(\fG),\,G \mapsto G_\p\,$, and $\,\L(\fG) \lra \L(\fZ),\,{\bf a} \mapsto {\bf a}_\p\,$, as
defined in \ref{Subsec:Pairing}, are lattice anti-isomorphisms inverse to one another. The last-mentioned fact and 
Corollary~\ref{Cor:Composing} are key ingredients used in \cite{A-B, Serdica, JPAA} in the investigation of 
the Galois connexion between the lattices $\,\L(\Gam)\,$ and $\,\L(\fZ)$.

In particular, if the action of $\,\Gam\,$ on $\,A\,$ is trivial then $\,\fZ = \text{Hom}(\Gam, A)\,$ and 
$\,\eta: \Gam \lra \fG = \text{Hom}(\Gam, A)^\vee\,$ is an isomorphism, therefore the Galois connexion between 
the lattices $\,\L(\Gam)\,$ and $\,\L(\text{Hom}(\Gam, A))\,$ is perfect. However, according to 
\cite[Proposition 5.6., Remarks 5.5., 5.8.]{JPAA} there exist nontrivial actions $\,\Gam \times A \lra A\,$ called
{\em strongly coGalois} such that the Galois connexion between the lattices $\,\L(\Gam)\,$ and $\,\L(\fZ)\,$ is
perfect and $\,\Gam \cong \fG\,$ (non-canonically) provided $\,\Gam\,$ is abelian.
\end{ex}
\begin{ex}
\label{Ex:2} \em
({\em The coGalois theory of separable radical extensions}) Let $\,L/K$ be a Galois extension with 
$\,\Gam:= \text{Gal}(L/K)\,$. In addition we assume that the extension $\,L/K\,$ is {\em radical}, i.e., 
$\,L = K(\text{Rad}(L/K))\,$, where
\[
\text{Rad}(L/K) = \{x \in L^\times\,|\,x^n \in K\,\mbox{for\,some}\,n \in \N \sm \{0\}\},
\]
so $\,\text{Rad}(L/K)^\Gam = (L^\times)^\Gam = K^\times\,$; as $\,L/K\,$ is separable, the exponents $\,n\,$ above 
may be assumed to be prime with the characteristic exponent of $\,K\,$. Note that any element $\,x \in \text{Rad}(L/K)$
is an $\,n$-th radical $\,\sqrt[n]{a}\,$ of an element $\,a \in K^\times$ for some $\,n \geq 1\,$; thus, 
$\,\text{Rad}(L/K)\,$ is precisely the set of all radicals belonging to $\,L\,$ of elements of $\,K^\times\,$.

With the notation from \ref{Subsec:coG} and \ref{Ex:1}, we consider the discrete $\,\Gam$-module $\,E:= L^\times\,$, 
the multiplicative group of the field $\,L\,$ with the Galois action, and we denote $\,A:= t(L^\times) = t(E) = \mu_L\,$,
with $\,A^\Gam = t(K^\times) = \mu_K\,$, the quasi-cyclic multiplicative group of the roots of unity contained in 
$\,L\,$, identified through a non-canonical monomorphism $\,\mu_L \lra \Q/\Z\,$ with a subgroup of $\,\Q/\Z\,$. 
Thus $\,\text{coG}(E) = \text{coG}(L/K):= t(L^\times/K^\times) = \text{Rad}(L/K)/K^\times\,$ is a discrete torsion 
abelian group with the induced Galois action. Since $\,H^1(\Gam, L^\times) = 0\,$ by Hilbert's Theorem 90, it 
follows by Lemma~\ref{Lem:H1} that the canonical map $\theta: \text{Rad}(L/K) \lra \fZ\,$ defined by
$\,\theta(x)(\gam) = \frac{\gam x}{x}\,$, for $\,x \in \text{Rad}(L/K), \gam \in \Gam\,$, induces an isomorphism
of the $\,\Gam$-modules $\,\text{coG}(E) = \text{coG}(L/K)\,$ and $\,\fZ\,$. 
 
Composing the perfect standard Galois connexion between the lattice $\,\L(L/K)$ of all intermediate fields 
$\,K \sse F \sse L\,$ and the lattice $\,\L(\Gam)\,$ with the Galois connexion between the lattices $\,\L(\Gam)\,$
and $\,\L(\fZ) \cong \L(\text{coG}(L/K))\,$ described in \ref{Subsec:GalCon}, \ref{Ex:1}, we obtain the coGalois 
connexion between the lattices $\,\L(L/K)\,$ and $\,\L(\text{coG}(L/K))\,$, the main object of investigation of the 
{\em coGalois theory of separable radical extensions\/}; consequently, all the results of this theory, in particular, 
the {\em Kummer theory}, could be obtained easily by transfering the corresponding results from its abstract version 
\ref{Ex:1}, \cite{A-B, Serdica, JPAA}. See \cite[6. Two examples]{JPAA} for a detailed analysis of the abstract version 
of the coGalois theory of the Galois extensions $\,L/K\,$, where $\,K\,$ is a finite field, $\,L = K^{\rm s}\,$,
respectively $\,K\,$ is a local field and $\,L\,$ its maximal tamely ramified extension.  
\end{ex}
\begin{ex}
\label{Ex:3} \em ({\em The additive analogue of} Example~\ref{Ex:2}) Using {\em Witt calculus} and 
{\em higher Artin-Schreier theory} \cite{A-S, Witt, FL}, we present an additive analogue of the multiplicative 
framework discussed in \ref{Ex:2}. 

Let $\,K\,$ be an arbitrary field of characteristic $\,p > 0\,$. Letting $\,V\,$ be the {\em shift operator} on 
the {\em Witt ring} $\,W(K)\,$, we consider the quotient rings $\,W_n(K):= W(K)/V^n W(K)\,$ consisting of 
{\em Witt vectors of length} $\,n \geq 1\,$ over $\,K\,$; in particular, $\,W_1(K) = K\,$. Since the shift operator 
$\,V: W(K) \lra W(K)\,$ is additive, the family of abelian groups $\,W_n(K)^+:= (W_n(K), +)\,$, indexed by the 
totally ordered set $\,\N^\ast:= \N \sm \{0\}\,$, forms a direct system with the canonical connecting 
monomorphisms $\,W_{n}(K)^+ \lra W_{n + 1}(K)^+\,$ induced by $\,V\,$. We denote 
$\,W_\infty(K):= \bigcup_{n \geq 1} W_n(K)^+\,$, the direct limit, identifying $\,W_n(K)^+\,(n \geq 1)\,$ with an
increasing sequence of subgroups of $\,W_\infty(K)\,$. On the other hand, the {\em Frobenius operator} $\,F\,$ is 
an endomorphism of the ring $\,W(K)\,$, and $F V = V F\,$, whence $\,F\,$ induces an endomorphism of the ring 
$\,W_n(K)\,$ for all $\,n \geq 1\,$ and so gives rise to an endomorphism, denoted also by $\,F\,$, of the abelian group 
$\,W_\infty(K)\,$ such that $\,F W_n(K)^+ \sse W_n(K)^+\,$ for all $\,n \geq 1\,$. Setting 
$\,e_n:= (1, 0, \dots, 0)\,$, the unit element of the ring $\,W_n(K)\,$, and using the identity $\,p x = F V x\,$ 
in the ring $\,W(K)\,$, we deduce that $\,e_n\,$ has order $\,p^n\,$. Consequently, $\,W_n(K)^+$ is a torsion abelian 
group of exponent $\,p^n\,$ containing $\,W_n(\F_p)^+ \cong p^{- n} \Z/\Z\,$, so $\,W_\infty(K)\,$ is a discrete 
$\,\Z_p$-module containing $\,W_\infty(\F_p) \cong (\Q/\Z)(p) = \Q_p/\Z_p$.

Now consider the family $\,(\fP_m)\,_{m \geq 1}\,$ of endomorphisms of the discrete $\,\Z_p$-module $\,W_\infty(K)\,$,
defined by $\,\fP_m(x) = F^m x - x\,$, and note that $\,\fP_m W_n(K)^+ \sse W_n(K)^+\,$ for $n, m \geq 1\,$. For 
$n = m = 1\,$ we obtain $\,\fP_1(x) = \fP(x) = x^p - x\,$ in $\,K^+ = W_1(K)^+\,$, the 
well known operator from the {\em Artin-Schreier theory} of abelian extensions of exponent $\,p\,$. It follows that
$\,\Ker(\fP_m) = W_\infty(K \cap \F_{p^m})\,$, in particular, $\,\Ker(\fP_1) = W_\infty(\F_p)$. 
 
Letting $\,L/K\,$ be a Galois extension with $\,\Gam:= {\rm Gal}(L/K)$, we denote by $\,k = K \cap \wt{\F_p}\,$,
$\,l:= L \cap \wt{\F_p}\,$ the relative algebraic closure of the prime field $\,\F_p\,$ in $\,K\,$ and $\,L\,$ 
respectively. Let $\,\overline{\Gam}\,$ be the procyclic group $\,{\rm Gal}(l/k)\,$, and 
$\,\Gam':= {\rm Gal}(L/(K \cdot l))\,$, the kernel of the epimorphism 
$\,\Gam \lra \overline{\Gam},\,\gam \mapsto\overline{\gam}:= \gam\,|_l\,$. The afore defined additive operators 
$\,\fP_m\,$ are seen as endomorphisms of $\,W_\infty(L)\,$ inducing by restriction endomorphisms of 
$\,W_\infty(E)\,$ for any subfield $\,E\,$ of $\,L\,$. The Galois group $\,\Gam\,$ acts continuously on the discrete
$\,\Z_p$-module $\,W_\infty(L),\,\gam \cdot W_n(L) = W_n(L)\,$ for all $\,\gam \in \Gam, n \geq 1\,$, 
$\,W_\infty(L)^\Gam = W_\infty(K),\,W_\infty(L)^{\Gam'} = W_\infty(K \cdot l)\,$, and $\gam\, \fP_m = \fP_m \gam\,$ for
all $\,\gam \in \Gam, m \geq 1\,$. For any $\,a = (a_0, \cdots, a_{n - 1}) \in W_n(L)\,$, we set 
$\,K(a):= K(a_0, \cdots, a_{n - 1})\,$, a finite subextension of $\,L/K\,$; more generally, for any subset $\,S\,$
of $\,W_\infty(L)\,$, we denote by $\,K(S)\,$ the composite of all $\,K(a)\,$ with $\,a \in S\,$. 

The additive version of the multiplicative group of radicals $\,{\rm Rad}(L/K)\,$ defined in \ref{Ex:2} is  
the $\,\Z_p[[\Gam]]$-submodule of $\,W_\infty(L)\,$ lying over $\,W_\infty(K) + W_\infty(l)\,$
\[
\mbox{Rad}_+(L/K):= \{x \in W_\infty(L)\,|\,\fP_m(x) \in W_\infty(K)\,\mbox{for\, some}\,m \geq 1\,\},
\]
the direct limit of the family of submodules $\,(\fP_m^{- 1}(W_\infty(K)))_{m \geq 1}\,$ of $\,W_\infty(L)\,$, with \\
$\,\fP_n^{- 1}(W_\infty(K)) \sse \fP_m^{- 1}(W_\infty(K))\,$ for $\,n\,|\,m\,$. Note that the $\,\Z_p[[\Gam]]$-submodule 
$\,W_\infty(l) =\,$ $\bigcup_{m \geq 1} \Ker(\fP_m)\,$ plays the role of the multiplicative group 
$\,\mu_L = t(L^\times)\,$ of the roots of unity contained in $\,L\,$.

The subextension $\,K({\rm Rad}_+(L/K))\,$ of the Galois extension $\,L/K\,$ is Galois; it is the the maximal 
subextension $\,N/K\,$ satisfying $\,N = K({\rm Rad}_+(N/K))\,$. By analogy with the radical extensions considered 
in \ref{Ex:2}, we say that $\,L/K\,$ is an {\em additive-radical extension} (for short {\em a-radical extension}) if 
$\,L = K({\rm Rad}_+(L/K))\,$. In particular, the usually called {\em Artin-Schreier extensions},
i.e., cyclic extensions of degree $\,p\,$ over a field of characteristic $\,p\,$, are minimal a-radical extensions. 

Note that for any $\,n \geq 1\,$, $\,M_n:= K(\fP_1^{- 1}(W_n(K)))\,$ is the maximal abelian subextension of 
$\,L/K\,$ whose exponent is a divisor of $\,p^n\,$, and hence $\,M_\infty:= \bigcup_{n \geq 1} M_n = 
K(\fP_1^{- 1}(W_\infty(K)))\,$ is the maximal abelian $\,p$-extension of $\,K\,$ contained in $\,L\,$; in particular, 
$\,M_\infty \cap (K \cdot l) = K \cdot l(p)\,$, where $\,l(p)\,$ is the maximal (procyclic) $\,p$-subextension of 
$\,l/k\,$. Using \cite[Propositions 1.1, 1.2]{G-S} as the first step of induction on $\,n\,$, we deduce that 
$\,{\rm Rad}_+(M_n/K) = \bigcup_{\F_{p^m} \sse k} \fP_m^{- 1}(W_n(K))\,$ for all $\,n \geq 1\,$, whence
$\,{\rm Rad}_+(M_\infty/K) = \bigcup_{\F_{p^m} \sse k} \fP_m^{- 1}(W_\infty(K))\,$.  

As an additive analogue of the coGalois group of a radical extension defined in \ref{Ex:2}, we take the quotient 
$\,\Z_p[[\Gam]]$-module $\,{\rm coG}_+(L/K):= {\rm Rad}_+(L/K)/W_\infty(K)\,$, the ``torsion'' 
of $\,W_\infty(L)/W_\infty(K)\,$ with respect to the endomorphisms $\,\fP_m,\,m \geq 1\,$. Considering the canonical 
action of $\,\Gam\,$ on $\,Z^1(\Gam, W_\infty(l))\,$ as defined in Remark~\ref{Rem:ActOnZ}, the map \\
$\,\theta: {\rm Rad}_+(L/K) \lra Z^1(\Gam, W_\infty(l))\,$, defined by $\,\theta(x)(\gam):= \gam x - x\,$, induces an
injective morphism of discrete $\,\Z_p[[\Gam]]$-modules $\,\wt{\theta}: {\rm coG}_+(L/K) \hra Z^1(\Gam, W_\infty(l))\,$,

Moreover, using the exact sequence of cohomology groups in low dimensions associated to the exact sequence of discrete 
$\,\Gam$-modules $0 \lra W_\infty(l) \lra W_\infty(L) \lra W_\infty(L)/W_\infty(l) \lra 0$, we deduce that 
$\,\wt{\theta}\,$ is an isomorphism since $\,H^1(\Gam, W_\infty(L)) = 0\,$ follows by induction from the 
additive Hilbert's Theorem 90. Thus, as in Example~\ref{Ex:2}, we have arrived to the abstract framework described in 
\ref{Subsec:GalCon}--\ref{Subsec:Pairing}: given a continuous action of a profinite group $\,\Gam\,$ on a discrete 
torsion $\,\Z_p$-module of the form $\,A:= W_\infty(l) = Q(W(l))/W(l) = \bigcup_{n \geq 1} p^{- n} W(l)/W(l)\,$, 
where $\,l\,$ is an algebraic extension of $\,\F_p\,$ and the fraction field $\,Q(W(l))\,$ of $\,W(l)\,$ is the 
unique, up to isomorphism, complete unramified discrete valued field of characteristic $\,0\,$ having $\,l\,$ as 
residue field, investigate the Galois connexion between the lattices $\L(\Gam)\,$ and $\,\L(Z^1(\Gam, A)\,$. Note
that we can replace the Witt ring $\,W(l)\,$ by the integral closure of $\,\Z_p\,$ in $\,W(l)\,$, the discrete 
valuation ring of the unramified field extension of $\,\Q_p\,$ having $\,l\,$ as residue field, provided the 
residue field extension $\,l/\F_p\,$ is infinite.
\end{ex}
\begin{ex} 
\label{Ex:4} \em ({\em The coGalois theory of Galois algebras over fields}) Let $\,K\,$ be a field, and $\,\fL\,$ 
be a discrete commutative $\,K$-algebra together with a continuous action of a profinite group $\,\Gam\,$. Thus 
$\,\fL\,$ is the union of its subalgebras $\,\fL^\Lam = \{f \in \fL\,|\,\gam f = f\,{\rm for\,all}\,\gam \in \Lam\}\,$,
where $\,\Lam\,$ ranges over the open normal subgroups of $\,\Gam\,$. In addition we assume that $\,\fL\,$ is a
{\em Galois $\Gam$-algebra}, i.e., for every open normal subgroup $\,\Lam\,$ of $\,\Gam\,$, the subalgebra 
$\,\fL^\Lam\,$ is {\em semisimple} ( a finite product of field extensions of $\,K\,$), and also a free 
$\,K[\Gam/\Lam]$-module of rank $1$ (there exists a {\em normal basis} for $\,\fL^\Lam\,|\,K\,$). With slight 
adaptations of the definitions and the arguments from \cite[5. Appendix: On Galois algebras]{L-R}, we obtain
\begin{lem}
\label{Lem:GalAlg}
$\,\fL\,|\,K\,$ is a Galois $\,\Gam$-algebra if and only if the following three conditions are satisfied.
\ben
\item[\rm (1)] The commutative ring $\,\fL\,$ is {\em regular} (in the sense of von Neumann), so $\,\fL\,$ is canonically
represented as a subdirect product of the family of the residue field extensions $\,(\fL/{\bf m})\,|\,K\,$,
where $\,{\bf m}\,$ ranges over all prime ($=$ maximal) ideals of $\,\fL\,$,
\item[\rm (2)] $\,\fL^\Gam = K\,$, and 
\item[\rm (3)] the continuous action of $\,\Gam\,$ on the profinite space $\,{\rm Spec}(\fL) = {\rm Max}(\fL)\,$,
canonically identified with the {\em Stone dual} of the boolean algebra $\,B(\fL)\,$ of all idempotents of
$\,\fL\,$, is transitive, and for some (for all) maximal ideal(s) $\,\bf{m}\,$, the stabilizer $\,\Gam_{\bf m}\,$ 
acts faithfully on the residue field $\,\fL/{\bf m}$. 
\een
\end{lem}
For a maximal ideal $\,{\bf m}\,$ of $\,\fL\,$ (uniquely determined by a primitive idempotent provided $\,\Gam\,$ 
is finite), the residue field $\,L = \fL/{\bf m}\,$ is a Galois extension of $\,K\,$, whose Galois group 
$\,{\rm Gal}(L/K)\,$ is isomorphic with the stabilizer $\,\Gam_{\bf m}\,$, the {\em decomposition group of 
$\,\fL\,|\,K\,$ associated to the maximal ideal} $\,{\bf m}\,$. Thus the Galois $\,\Gam$-algebra $\,\fL\,|\,K\,$ 
determines a conjugation class of closed subgroups of $\,\Gam\,$ as its decomposition groups. As stressed in 
\cite[A.7]{L-R}, the theory of Galois $\,\Gam$-algebras over a field $\,K\,$ is essentially the same as the theory
of Galois field extensions $\,L/K\,$, together with an embedding $\,{\rm Gal}(L/K) \hra \Gam\,$. Identifying 
$\,{\rm Gal}(L/K)\,$ with a closed subgroup $\,\Gam'\,$ of $\,\Gam\,$, the corresponding Galois $\,\Gam$-algebra 
$\,\fL\,$ is up to isomorphism the {\em induced algebra} $\,M_\Gam^{\Gam'}(L)\,$ of all continuous (locally constant) 
functions $\,f : \Gam \lra L\,$ satisfying $\,f(\si \gam) = \si f(\gam)\,$ for all $\,\si  \in \Gam', \gam \in \Gam\,$,
with the ring operations pointwise induced from the field $\,L\,$, the base field $\,K\,$ identified with the subfield 
of all constant functions, and the faithful action of $\,\Gam\,$ defined by $\,(\si f)(\gam) = f(\gam \si)\,$ for 
$\,f \in \fL\,$ and $\,\si, \gam \in \Gam\,$. 

For every $\,\gam \in \Gam\,$, let $\,\varphi_\gam \in {\rm Hom}_K(\fL, L)\,$ denote the epimorphism defined by 
$\,\varphi_\gam(f) = f(\gam)\,$ for $\,f \in \fL\,$. The map 
$\varphi : \Gam \lra {\rm Hom}_K(\fL, L),\,\gam \mapsto \varphi_\gam\,$, is a homeomorphism, considering on  
$\,{\rm Hom}_K(\fL, L)\,$ the Zariski ($=$ Stone) topology whose closed sets have the form 
$\,Z(I):= \{\,\psi \in {\rm Hom}_K(\fL, L)\,|\,\psi\,|_I = 0\,\}$, where $\,I\,$ is an ideal of $\,\fL\,$; note 
that $\,Z(I) = Z(I \cap B(\fL))\,$. The canonical action from the right of $\,\Gam\,$ on itself induced by 
multiplication determines via the homeomorphism $\,\varphi\,$ a continuous simple transitive action from the 
right of $\,\Gam\,$ on the profinite space $\,{\rm Hom}_K(\fL, L)$, defined by 
$\,(\psi \gam)(f) = \psi(\gam f)\,$ for $\,\gam \in \Gam, \psi \in {\rm Hom}\,_K(\fL, L), f \in \fL\,$.
The canonical continuous projection $\,\Gam \lra {\rm Max}(\fL),\,\gam \mapsto \Ker(\varphi_{\gam^{- 1}}) = 
\gam\,{\bf m}\,$, with $\,{\bf m} = \Ker(\varphi_1)\,$, induces an isomorphism of profinite $\,\Gam$-spaces
$\,\Gam/\Gam_{\bf m} \lra {\rm Max}(\fL)\,$, with $\,\Gam_{\bf m} = \Gam' = {\rm Gal}(L/K)$.   

By analogy with \ref{Ex:2}, call {\em coGalois group of the Galois $\,\Gam$-algebra} $\,\fL\,|\,K\,$ the 
coGalois group $\,{\rm coG}(\fL^\times):= t(\fL^\times/K^\times)\,$ of the discrete $\,\Gam$-module 
$\,\fL^\times = M_\Gam^{\Gam'}(L^\times)\,$ of units of the $\,K$-algebra $\,\fL\,$. Since $\,H^1(\Gam, \fL^\times) = 
H^1(\Gam, M_\Gam^{\Gam'}(L^\times)) \cong H^1(\Gam', L^\times) = 0\,$ by Shapiro's lemma \cite[Proposition 1.6.3]{NSW}
and Hilbert's Theorem 90, it follows by Lemma~\ref{Lem:H1} that the discrete torsion 
$\,\Gam$-module $\,{\rm coG}(\fL\,|\,K):= {\rm coG}(\fL^\times)\,$ is canonically isomorphic with 
$\,Z^1(\Gam, \mu_\fL)\,$, where $\,\mu_\fL:= t(\fL^\times) = M_\Gam^{\Gam'}(\mu_L)\,$. 

In abstract terms, we have arrived to a framework extending \ref{Ex:1} as follows. Let $\,\Gam\,$ be a profinite group, 
$\,\Gam'\,$ a closed subgroup of $\,\Gam\,$, $\,A\,$ a discrete quasi-cyclic group identified with a subgroup of 
$\,\Q/\Z\,$, $\,A^\vee:= \text{Hom}(A, \Q/\Z)\,$ its Pontryagin dual, and $\,\Gam' \times A \lra A,\,(\si, a) \mapsto 
\si a:= \chi(\si) a\,$, a continuous action given by a continuous homomorphism $\,\chi: \Gam' \lra (A^\vee)^\times\,$. 
Let $\,\fA:= M_\Gam^{\Gam'}(A)\,$ be the discrete torsion $\,\Gam$-module of all continuous functions 
$\,f: \Gam \lra A\,$ satisfying $\,f(\si \gam) = \chi(\si) \cdot f(\gam)\,$ for $\,\si \in \Gam', \gam \in \Gam\,$,
with the action $\,\Gam \times \fA \lra \fA\,$ defined by $\,(\si f)(\gam):= f(\gam \si)\,$ for $\,f \in \fA, \si, \gam
\in \Gam\,$; in particular, the map $\,\fA \lra A,\,f \mapsto f(1)\,$ is an isomorphism provided $\,\Gam = \Gam'\,$.
The object of investigation is the Galois connexion between the lattices $\,\L(\Gam)\,$ and $\,\L(Z^1(\Gam, \fA))\,$ 
using the framework described in \ref{Subsec:GalCon}-\ref{Subsec:Pairing}.
\end{ex} 

\section{Kneser and minimal non-Kneser triples}
\label{Sec:Kneser}

In this section we investigate a remarkable class of g-cocycles and associated structures, extending to 
the more general framework introduced in Section~\ref{Sec:GenFrame} the main results on {\em Kneser groups of cocycles} 
from \cite{A-B}. 
\subsection{Surjective cocycles}
\label{Subsec:SurCocycles}
Let $\,\Gam\,$ be a profinite group together with a continuous action on a profinite group $\,\fG\,$ and a continuous
$\,1$-cocycle $\,\eta: \Gam \lra \fG\,$. Setting $\,\De := \Ker(\eta),\,\fE := \fG \rtimes \Gam,\,p: \fE \lra \Gam\,$ 
the natural projection, $\,s_1 : \Gam \lra \fE\,$ the canonical section, $\,s_2 : \Gam \lra \fE\,$ the homomorphic
section induced by the cocycle $\,\eta\,$, and $\,\Gam_i :=
s_i(\Gam) \cong \Gam,\,i = 1, 2\,$, we obtain
\begin{lem} The following assertions are equivalent.
\label{Lem:SurCocycle} 

{\rm (1)} The cocycle $\,\eta : \Gam \lra \fG\,$ is surjective.

{\rm (2)} For any open ideal $\,{\bf a}\,$ of $\,\fG\,$, the induced cocycle $\,\eta_{{\bf a}}: \Gam \lra \fG/{\bf a}\,$
is surjective.

{\rm (3)} The map $\,\Gam/\De \lra \fG\,$ induced by $\,\eta\,$ is a homeomorphism.

{\rm (4)} The map $\,\Gam_1/(\Gam_1 \cap \Gam_2) \lra \fE/\Gam_2\,$ induced by the inclusion $\,\Gam_1 \sse \fE\,$ 
is a homeomorphism.

{\rm (5)} The map $\,\Gam_2/(\Gam_1 \cap \Gam_2) \lra \fE/\Gam_1\,$ induced by the inclusion $\,\Gam_2 \sse \fE\,$ 
is a homeomorphism.

{\rm (6)} $\,\fE = \Gam_1 \cdot \Gam_2 = \{\gam_1 \gam_2\,|\,\gam_i \in \Gam_i, i = \,1, 2\}$.

{\rm (7)} $\,\fE = \Gam_2 \cdot \Gam_1$.

{\rm (8)} The continuous map $\,\Gam \times \fG \lra \fG,\,(\gam, g) \mapsto \eta(\gam) \cdot \gam g\,$, is a 
transitive action on the underlying profinite space of $\,\fG\,$ having the closed subgroup $\,\De\,$ of $\,\Gam\,$ 
as stabilizer of the neutral element of $\,\fG$.
\end{lem} 
\bp
It suffices to note that the isomorphisms of profinite groups $\,\Gam \lra \Gam_i,\,i = 1, 2\,$,
induce the homeomorphisms of profinite spaces $\,\Gam/\De \lra \Gam_i/(\Gam_1 \cap \Gam_2),\,i = 1, 2\,$,
while the inclusion $\,\fG \sse \fE\,$ induces the homeomorphisms $\,\fG \lra \fE/\Gam_i,\,i = 1, 2\,$. On the other 
hand, the continuous maps from (3), (4) and (5) are obviously injective.
\ep
\begin{co}
\label{Cor:SurCocycle}
With notation above, assume that the $\,\Gam$-group $\,\fG\,$ is finite, so $\,\De\,$ is open. The following 
assertions are equivalent.

{\rm (1)} The cocycle $\,\eta: \Gam \lra \fG$ is surjective.

{\rm (2)} $\,(\Gam : \De) = (\Gam_1 : \Gam_1 \cap \Gam_2) = (\Gam_2 : \Gam_1 \cap \Gam_2) \geq\,|\fG|\,= 
(\fE : \Gam_1) = (\fE : \Gam_2)$.

{\rm (3)} $\,(\Gam : \De) = (\Gam_1 : \Gam_1 \cap \Gam_2) = (\Gam_2 : \Gam_1 \cap \Gam_2) =\,|\fG|\,= 
(\fE : \Gam_1) = (\fE : \Gam_2)$.
\end{co}
\begin{de}
\label{Def:KneserTriple}
With notation above, we call $\,(\Gam, \fG, \eta)\,$ a {\em Kneser triple} if the continuous cocycle 
$\,\eta: \Gam \lra \fG\,$ is surjective. 
\end{de}
Denote by $\,\cK\cZ^1\,$ the full subcategory of the category $\,\cZ^1\,$ defined in \ref{Subsec:GenCocycles}, having 
the Kneser triples as objects. 

\begin{rem} \em
\label{Rem:KneserStr}
According to Lemma~\ref{Lem:SurCocycle}, (8), the category $\,\cK\cZ^1\,$ is equivalent 
with the category of the systems $\,(\Gam, \De, \bullet)\,$, termed {\em Kneser structures}, consisting of a profinite 
group $\,\Gam\,$, a closed subgroup $\,\De\,$, and a continuous group operation $\,\bullet\,$ on the profinite space 
$\,X:= \Gam/\De = \{ \wh{\gam} := \gam \De\,|\,\gam \in \Gam \}$, with $\,\wh{1}\,$ as neutral element, such that the 
canonical transitive action of $\,\Gam\,$ on $\,X\,$ and the group operation $\,\bullet\,$ are related 
through the following condition
\[
\gam \cdot (x \bullet y) = (\gam \cdot x) \bullet I(\wh{\gam}) \bullet (\gam \cdot y)\,\,{\rm for}\,\gam \in \Gam,
x, y \in X,
\]
where $\,I(x)$ denotes the inverse of $\,x \in X\,$ with respect to the group operation $\,\bullet\,$; in particular,
$\,\gam \cdot I(x) = \wh{\gam} \bullet I(\gam \cdot x) \bullet \wh{\gam}\,$ for $\,\gam \in \Gam,\,x \in X\,$, and 
the restricted action of $\,\De\,$ on $\,X\,$ is compatible with the group operation $\,\bullet$. A morphism 
$\,(\Gam, \De, \bullet) \lra (\Gam', \De', \bullet)$ is a continuous morphism $\,\varphi : \Gam \lra \Gam'\,$ 
inducing a morphism of profinite groups $\,\wt{\varphi} : (\Gam/\De, \bullet) \lra (\Gam'/\De', \bullet')\,$. 
\end{rem}

\subsubsection{Surjectivity criteria for cocycles}
The next result, extending \cite[Proposition 1.14, Corollary 1.16, Corollary 1.17]{A-B}, provides a criterion for 
cocycles taking values in pronilpotent groups to be surjective.
\begin{pr} 
\label{Pr:Pronil}
Let $\,\eta \in Z^1(\Gam, \fG)$, where $\,\fG \cong \dprod_{p} \fG_p\,$ is a pronilpotent $\,\Gam\,$-group, and 
$\,\fG_p\,$ denotes the maximal pro-$p$-quotient of $\,\fG\,$ for any prime $\,p\,$, with the induced action of 
$\,\Gam\,$. Then the cocycle $\,\eta: \Gam \lra \fG\,$ is surjective if and only if the induced cocycle 
$\,\eta_p: \Gam \lra \fG_p\,$ is surjective for all prime numbers $\,p\,$.
\end{pr}
The criterion above is a consequence of Lemma~\ref{Lem:SurCocycle} (2) and the next lemma.
\begin{lem}
Let $\,\eta \in Z^1(\Gam, \fG)$, where $\,\fG = \dprod_{i = 1}^n \fG_i\,$ is a direct product of finite $\,\Gam$-groups 
such that $\,{\rm gcd}\,(|\fG_i|,\,|\fG_j|) = 1\,\text{for}\,i \neq j\,$. Then $\,\eta: \Gam \lra \fG$ is surjective if 
and only if the induced cocycle $\,\eta_i: \Gam \lra \fG_i\,$ is surjective for $\,i = 1, \dots, n$. 
\end{lem}
\bp
An implication is obvious. Conversely, assuming that $\,\eta_i: \Gam \lra \fG_i\,$ is surjective for 
$\,i = 1, \dots, n\,$, we obtain $\,|\fG| = \prod_{i = 1}^{n} |\fG_i| = \prod_{i = 1}^{n} (\Gam : \Ker(\eta_i))\,$
by Corollary~\ref{Cor:SurCocycle},(3). Since $\,\Ker(\eta) \sse \Ker(\eta_i)\,$, it follows that
$\,(\Gam : \Ker(\eta_i)) \,|\,(\Gam : \Ker(\eta))\,$, for $\,i = 1, \dots, n$. As $\,(\Gam : \Ker(\eta_i)) = |\fG_i|\,$, 
$\,i = 1, \dots, n\,$, are mutually relatively prime by hypothesis, it follows that 
$\,\prod_{i = 1}^{n} (\Gam : \Ker(\eta_i))\,|\,(\Gam : \Ker(\eta))$, whence $|\fG|\,|\,(\Gam : \Ker(\eta))\,$. 
Consequently, the cocycle $\,\eta$ is surjective by Corollary~\ref{Cor:SurCocycle},(2).
\ep
Another surjectivity criterion for g-cocycles is given by the next lemma.
\begin{lem}
\label{Lem:SurCrit}
Given a g-cocycle $\,\eta: \Gam \lra \fG\,$, the following assertions are equivalent.

{\rm (1)} $\,\eta: \Gam \lra \fG\,$ is surjective.

{\rm (2)} There exists a closed subgroup $\,\Lam \sse \Gam\,$ such that $\,{\bf a}:= \eta(\Lam)\,$ is an ideal 
of $\,\fG\,$, and the induced cocycle $\,\eta_{\bf a}: \Gam \lra \fG/{\bf a}\,$ is surjective. 
\end{lem}
\bp
The implication (1)$\Lra$(2) is obvious. To prove the converse, assume $\,\Lam \in \L(\Gam)\,$ satisfies (2), with
$\,{\bf a} = \eta(\Lam) \in \L(\fG)\,$, and let $\,g \in \fG\,$. As $\,\eta_{\bf a}: \Gam \lra \fG/{\bf a}$ is 
surjective by assumption, there exists $\,\si \in \Gam\,$ such that $\,\eta(\si)^{- 1} g \in {\bf a}\,$, and
hence $\,\eta(\si^{- 1}) \cdot (\si^{- 1} g) = \si^{- 1} (\eta(\si)^{- 1} g) \in {\bf a}\,$ since $\,{\bf a}\,$ is 
$\,\Gam$-invariant. On the other hand, as $\,{\bf a} = \eta(\Lam)\,$, we obtain 
$\,\eta(\si^{- 1}) \cdot (\si^{- 1} g) = \eta(\tau)\,$ for some $\,\tau \in \Lam\,$. Consequently, acting with $\,\si\,$,
it follows that $\,g = \eta(\si) \cdot (\si \eta(\tau)) = \eta(\si \tau)\,$, so $\,\eta(\Gam) = \fG\,$ as desired.
\ep
Now, given a g-cocycle $\,\eta: \Gam \lra \fG\,$, let $\,\overline{\De}:= \De' \cap \De''\,$, where
$\,\De':= \text{Fix}_\Gam(\fG)\,$, $\,\De'':= \{\si \in \Gam\,|\,\forall \gam \in \Gam,\,\eta(\gam \si \gam^{- 1}) = 
\gam \eta(\si)\}\,$. According to Lemma~\ref{Lem:PropGenCocycles}, $\,\overline{\De}\,$ is a closed normal subgroup
of $\,\Gam\,$, and $\,\eta(\overline{\De})\,$ is a central ideal of $\,\fG\,$. Applying Lemma~\ref{Lem:SurCrit} to
$\,{\bf a}:= \eta(\overline{\De})\,$, we obtain
\begin{co}
\label{Cor:SurCrit}
With data above, the g-cocycle $\,\eta: \Gam \lra \fG\,$ is surjective provided the induced cocycle
$\,\eta_{\bf a}: \Gam \lra \fG/{\bf a}\,$ is surjective. 
\end{co}

\begin{rem} \em
Though not stated explicitly, a particular form of Corrolary~\ref{Cor:SurCrit} plays implicitly a key role in 
the proof of \cite[Lemma 1.18, Theorem 1.20]{A-B}.
\end{rem}
\subsubsection{Bijective cocycles induced by self-actions}
\label{Subsubsec:BijectCocycles}
We describe a procedure for obtaining Kneser triples $\,(\Gam, \fG, \eta)\,$ where the cocycle 
$\,\eta\,$ is bijective. Let $\,\Gam\,$ be a profinite group, and $\,\om: \Gam \lra {\rm Aut}(\Gam)\,$
a continuous action of $\,\Gam\,$ on itself. Setting $\,\wt{\Gam}:= \Gam \rtimes \om(\Gam)\,$, $\,\wt{\Gam}\,$ becomes
a profinite $\,\Gam$-group via the canonical continuous action
\[
\Gam \times \wt{\Gam} \lra \wt{\Gam},\,(\gam, (\del, \theta)) \mapsto (\om(\gam)(\del),\,\om(\gam) \circ \theta \circ 
\om(\gam)^{- 1}), 
\]
and the continuous map $\,\eta_\om: \Gam \lra \wt{\Gam},\,\gam \mapsto (\gam, \om(\gam)^{- 1})\,$ is an 
injective cocycle.
\begin{de}
\label{Def:AdeqAction}
The continuous self-action $\,\om: \Gam \lra {\rm Aut}(\Gam)\,$ is said to be {\em adequate} if the image 
$\,\eta_\om(\Gam)\,$ is a (closed) subgroup of $\,\wt{\Gam}\,$, i.e., 
$\,\om(\theta(\gam)) = \theta \circ \om(\gam) \circ \theta^{- 1}\,$ for all $\,\gam \in \Gam, \theta \in \om(\Gam)\,$; we  
say that the profinite $\,\Gam$-group $\,\Gam_\om:= \eta_\om(\Gam)\,$ is a {\em deformation of} $\,\Gam\,$ via the
adequate self-action $\,\om$.
\end{de}
In other words, the self-action $\,\om: \Gam \lra {\rm Aut}(\Gam)\,$ is adequate if and only if the continuous binary 
operation on $\,\Gam\,$, defined by $\,\gam \bullet \del = \gam \cdot \om(\gam)^{- 1}(\del)\,$, is associative; in this
case, $\,(\Gam, \bullet)\,$ is a profinite $\,\Gam$-group isomorphic to $\,\Gam_\om\,$, and the identity 
$\,(\Gam, \cdot) \lra (\Gam, \bullet)\,$ is a cocycle, i.e. $\,\gam \cdot \del = \gam \bullet \om(\gam)(\del)\,$ for
$\,\gam, \del \in \Gam$. 

For any adequate self-action $\,\om: \Gam \lra {\rm Aut}(\Gam)\,$, $\,(\Gam, \Gam_\om, \eta_\om: \Gam \lra \Gam_\om)\,$ 
is a Kneser triple with a bijective cocycle $\,\eta_\om$. Note that $\,\Ker(\om)\,$ is identified with a common closed
normal subgroup of $\,\Gam\,$ and $\,\Gam_\om\,$ on which the cocycle $\,\eta_\om\,$ is the identity. Moreover $\,\Gam\,$
becomes a profinite $\,\Gam_\om$-group via the continuous action 
\[
\Gam_\om \times \Gam \lra \Gam,\,((\gam, \om(\gam)^{- 1}), \del) \mapsto \om(\gam)^{- 1}(\del), 
\]
whose kernel is identified with $\,\Ker(\om)\,$, and $(\Gam_\om, \Gam, \eta_\om^{- 1}: \Gam_\om \lra \Gam)\,$ is also 
a Kneser triple with a bijective cocycle $\eta_\om^{- 1}$. We obtain in this way a symmetric relation on profinite 
groups, weaker than the isomorphism-relation and finer than the cardinality equivalence, defined as follows.
\begin{de}
\label{Def:DefPair}
For profinite groups $\,\Gam\,$ and $\,\Lam\,$, $\,(\Gam, \Lam)\,$ is a {\em deformation pair} if the following 
equivalent conditions are satisfied.

{\rm (1)} There exists an adequate self-action $\,\om: \Gam \lra {\rm Aut}(\Gam)\,$ such that $\,\Lam \cong \Gam_\om$.

{\rm (2)} There exist continuous actions by automorphisms $\,\Gam \times \Lam \lra \Lam\,$,
$\Lam \times \Gam \lra \Gam\,$, and a continuous bijective cocycle $\,\eta \in Z^1(\Gam, \Lam)\,$ such that 
$\,\eta^{- 1} \in Z^1(\Lam, \Gam)$.
\end{de}
\begin{prob} 
Is the deformation relation transitive, i.e., an equivalence relation ?
\end{prob}
A profinite group $\,\Gam\,$ is termed {\em rigid}, if $\,\Gam \cong \Lam\,$ provided $\,(\Gam, \Lam)\,$ is a 
deformation pair. Among the rigid profinite groups, we mention the pro-$p$-cyclic groups, where $\,p\,$ is an odd prime
number, and the finite simple groups. By contrast, the pro-$2$-cyclic groups of order $\,\geq 4\,$ are not rigid; indeed,
for $\,\Gam = \Z/2^n \Z\,$, with $\,2 \leq n \in \N\,$, $(\Gam, \D_{2^n} = \Z/2^{n - 1} \Z \rtimes \Z/2 \Z)\,$ is a 
deformation pair, while $\,\D_{2^\infty} = \Z_2 \rtimes \Z/2 \Z\,$ is, up to isomorphism, the
unique profinite group $\,\Lam \not \cong \Z_2\,$ such that $\,(\Z_2, \Lam)\,$ is a deformation pair (see also
Example~\ref{Ex:SelfAct1}). 

Note that not all bijective cocycles are induced by adequate self-actions. 

\subsubsection{Examples}

\begin{ex}
\label{Ex:Local1} \em
Let $\,K\,$ be a local field with discrete valuation $\,v\,$, valuation ring $\,\cO\,$, maximal ideal $\,\fp\,$, 
finite residue field $\,k := \cO/\fp\,$, and $\,U^{(n)}:= 1 + \fp^n\,$, the group of $\,n$-units of $\,K\,$ for some 
natural number $\,n \geq 1\,$. Fix a local uniformizer $\,\pi\,$, and consider the faithful continuous action 
by multiplication of $\,U^{(n)}\,$ on $\,\cO^+\,$. The map $\,Z^1(U^{(n)}, \cO^+) \lra \cO^+,\,\eta
\mapsto \eta(1 + \pi^n)\,$ is an isomorphism whose inverse sends an element $\,a \in \cO\,$ to the 
cocycle $\,\eta_a\,$ defined by $\,\eta_a(x):= \pi^{- n} (x - 1) a\,$. The isomorphism above maps the subgroup 
$\,B^1(U^{(n)}, \cO^+)\,$ onto $\,\fp^n\,$ inducing an isomorphism $\,H^1(U^{(n)}, \cO^+) \cong \cO/\fp^n\,$. For 
$\,a \neq 0,\,\eta_a\,$ is injective, while $\,\eta_a\,$ is a g-cocycle if and only if it is bijective, i.e., 
$\,a \in \cO^\times\,$; one checks easily that for all $\,a \in O^\times\,$, the bijective cocycle $\,\eta_a\,$ is not
induced by an adequate self-action. Note that two cocycles $\,\eta_a\,$ and $\,\eta_b\,$ are cohomologous if and only 
if $\,a \equiv b\,{\rm mod}\,\fp^n\,$, while they are isomorphic as objects in the category $\,\cZ^1\,$ if and only if 
$\,a \cO = b \cO\,$, i.e., $\,v(a) = v(b)\,$. In particular, $\,\eta_a \cong \eta_1\,$ for all $\,a \in \cO^\times\,$, 
with the associated Kneser structure $\,(U^{(n)}, 1, \bullet)\,$, with $\,x \bullet y = x + y - 1\,$, neutral element 
$\,1\,$, and inverse $\,I(x) = 2 - x$. For any natural number $\,m \geq 1\,$, the bijective cocycle $\,\eta:= \eta_1$ 
induces a bijective cocycle $\,\eta^{(m)} \in Z^1(U^{(n)}/U^{(n + m)}, \cO/\fp^m)\,$ and $\,\eta = \dlim \eta^{(m)}$.
\end{ex}
\begin{ex}
\label{Ex:Local2} \em
This example plays a key role in the proof of Proposition~\ref{Pr:Ab2}. With the data from Example~\ref{Ex:Local1}, 
let $\,m > n \geq 1\,$, and consider the continuous action of $\,U^{(n)}\,$ on the quotient abelian group 
$\,\cO_m:= (\cO/\fp^m)^+\,$ with $\,{\rm Fix}_{U^{(n)}}(\cO_m) = U^{(m)}$ and 
$\,H^0(U^{(n)}, \cO_m) = \fp^{m - n} \cO_m\,$. For any $\,a \in \cO_m\,$, put 
$\,v_m(a):= {\rm max}\{ i \in \{ 0, 1, \dots, m \}\,|\,a \in \fp^i \cO_m\}\,$, so 
$\,\fp^i \cO_m = \pi^i \cO_m = \{ a \in \cO_m\,|\,v_m(a) \geq i \}$ for $\,i = 0, 1, \dots, m\,$.
Denote by $\,V\,$ the closed subgroup of $\,U^{(n)}\,$ generated by $\,1 + \pi^n\,$. We obtain an isomorphism 
\[
Z^1(U^{(n)}, \cO_m) \lra \cO_m \times {\rm Hom}(U^{(n)}/V, \fp^{m - n} \cO_m)
\] 
sending a cocycle $\,\eta\,$ to the pair $\,(\eta(1 + \pi^n), \wt{\eta})\,$, where the continuous homomorphism 
$\,\wt{\eta}\,$ is defined by $\,\wt{\eta}(x) = \eta(x) -\pi^{- n} (x - 1) \eta(1 + \pi^n)\,$ for all 
$\,x \in U^{(n)}\,$. The isomorphism above maps $\,B^1(U^{(n)}, \cO_m)\,$ onto $\,\fp^n \cO_m\,$ inducing 
an isomorphism  
\[
H^1(U^{(n)}, \cO_m) \cong \cO_n \times {\rm Hom}(U^{(n)}/V, \fp^{m - n} \cO_m) = H^1(U^{(n)}, \cO_n) \oplus 
{\rm Hom}(U^{(n)}/V, H^0(U^{(n)}, \cO_m))\,.
\]
For $\,a \in \cO_m, \alpha \in {\rm Hom}(U^{(n)}/V, \fp^{m - n} \cO_m)\,$, let $\,\eta_{a, \alpha}\,$
denote the unique cocycle $\,\eta\,$ satisfying $\,\eta(1 + \pi^n) = a, \wt{\eta} = \alpha\,$,
i.e., $\,\eta_{a, \alpha}(x) = \alpha(x) + \pi^{- n} (x - 1) a\,$ for all $\,x \in U^{(n)}\,$.
It follows that 
\[
\Ker(\eta_{a, \alpha}) = {\rm Eq}(\alpha, \eta_{- a, 0}) = \{ x \in U^{(n)}\,|
\,\alpha(x) = \pi^{-n} (1 - x) a \} \sse U^{(r(a))}\,,
\] 
where $\,r(a) := {\rm max}(n, m - v_m(a))\,$. Note that $\,U^{(r(a))} = \eta_{a, \alpha}^{- 1}(\fp^{m - n} \cO_m)\,$ 
is the maximal open subgroup $\,H\,$ of $\,U^{(n)}\,$ for which the restriction map $\,\eta_{a, \alpha}\,|_H\,$ is a 
homomorphism, so a cocycle $\,\eta\,$ is a homomorphism if and only if $\,\eta(1 + \pi^n) \in \fp^{m - n} \cO_m\,$, and
$\,\Ker(\eta_{a, \alpha}) = \Ker(\eta_{a, \alpha}\,|_{U^{(r(a))}})$. 
\end{ex}
\begin{ex}
\label{Ex:SelfAct1} \em
Let $\,\Gam = (\Z/n \Z,+)\,$ be a finite cyclic group of order $\,n \geq 2\,$. As $\,{\rm Aut}(\Gam) \cong 
(\Z/n \Z)^\times\,$, the self-actions of $\,\Gam\,$ on itself are identified via the map 
$\,u \mapsto \om_u: \Gam \lra {\rm Aut}(\Gam)\,$, where $\,\om_u(x)(y) = u^x \cdot y\,$, with the multiplicative subgroup
$\,\cU_n\,$ of $\,(\Z/n \Z)^\times\,$ consisting of those units $\,u \in (\Z/n \Z)^\times\,$ satisfying 
$\,u^n \equiv 1\,{\rm mod}\,n\,$. The adequate self-actions $\,\om_u\,$ are in 1-1 correspondence to those 
$\,u \in \cU_n\,$ satisfying the supplementary condition $\,u^{u - 1} \equiv 1\,{\rm mod}\,n\,$. Thus the adequate 
self-actions of $\,\Z/n \Z\,$ are identified with the set \\$\,\cU_n^{{\rm ad}} = 
\{u \in (\Z/n \Z)^\times\,|\,u^{r_u} \equiv\,1\,{\rm mod}\,n\}\,$, where $\,r_u:= (n, u - 1)\,$ for 
$\,u \in (\Z/n \Z)^\times\,$. We obtain a more explicit description of the set $\,\cU_n^{{\rm ad}}$ as follows.
\begin{lem}
\label{Lem:AdUnits}
Let $\,\cP_n\,$ be the set of odd prime divisors of $\,n\,$, to which we add the prime number $\,2\,$ provided 
$\,4\,|\,n\,$. For $\,u \in (\Z/n \Z)^\times\,$ and $\,p \in \cP_n\,$, put $\,u_p:= u\,{\rm mod}\,p^{v_p(n)}\,$, and 
denote by $\,o(u_p)\,$ the order of $\,u_p\,$. Then $\,u \in \cU_n^{{\rm ad}}\,$ if, and only if, the following 
conditions are satisfied.

{\rm (i)} $\,r_u \neq 1\,$, i.e., $\,u - 1 \notin (\Z/n \Z)^\times$.

{\rm (ii)} $\,p \in \cP_{r_u} \Lra v_p(n) \leq 2 v_p(r_u)$.

{\rm (iii)} $\,2 \neq p \in \cP_n \sm \cP_{r_u} \Lra 2 \leq o(u_p)\,|\,(r_u, p - 1)$.

{\rm (iv)} $\,2 \in \cP_n \sm \cP_{r_u} \Lra u \equiv - 1\,{\rm mod}\,2^{v_2(n) - 1}\,$, in particular, 
$\,o(u_2) = 2$.
\end{lem}
Lemma~\ref{Lem:AdUnits} remains valid for arbitrary procyclic groups of order $\,n\,$-a supernatural number. In
particular, $\,\cU_{p^\infty}^{{\rm ad}} = \cU^{{\rm ad}}((\Z_p)^\times) = \{1\}\,$ for $\,p \neq 2\,$, while 
$\,\cU_{2^\infty}^{{\rm ad}} = \cU^{{\rm ad}}((\Z_2)^\times) = \{1, - 1\}$; in the latter case, for $\,\Gam = \Z_2, u = - 1\,$,
we obtain $\,\Gam_{- 1} \cong \Z_2 \rtimes \Z/2 \Z \cong \Z/2 \Z \wh{*}_2 \Z/2 \Z\,$(the free product of two copies of
the cyclic group of order $\,2\,$ in the category of pro-$2$-groups) \footnote[2]{The procedure described in
\ref{Subsubsec:BijectCocycles} works for arbitrary topological groups, in particular, for discrete ones. 
For $\,\Gam = (\Z, +)\,$, we obtain $\,\cU^{{\rm ad}} = \Z^\times = \{\pm 1\}\,$, and $\,\Gam_{- 1} = \D_\infty = 
\Z \rtimes \Z/2 \Z \cong (\Z/2 \Z) * (\Z/2 \Z)$.}, the coGalois triple \footnote[3]{see Definition~
\ref{Def:coGTriple}.} $\,(\Gam, \Gam_{- 1}, \eta_{- 1})\,$, and the Kneser, but not coGalois, triple 
$\,(\Gam_{- 1}, \Gam, \eta_{- 1}^{- 1})\,$.

For $\,n = 4, \Gam = \Z/4 \Z\,$, $\,\cU_4^{{\rm ad}} = (\Z/4 \Z)^\times \cong \Z/2 \Z\,$. Setting 
$\,\wh{i}:= i\,{\rm mod}\,4\,$, we obtain $\,\Gam_{\wh{1}} \cong \Gam,\,\Gam_{\wh{3}} \cong \Z/2 \Z \times \Z/2 \Z\,$,
and $\,(\Gam, \Gam_{\wh{3}}, \eta_{\wh{3}})\,$ is a coGalois \footnotemark[3] triple, while 
$\,(\Gam_{\wh{3}}, \Gam, \eta_{\wh{3}}^{- 1})\,$ is a Kneser minimal non-coGalois triple 
\footnote[4]{see Definition~\ref{De:mncG} and Example~\ref{Ex:SelfAct1}\,(i).}.

For $\,n = 8, \Gam = \Z/8 \Z\,$, $\,\cU_8^{{\rm ad}} = (\Z/8 \Z)^\times \cong \Z/2 \Z \times \Z/2 \Z\,$. Setting 
$\,\wh{i}:= i\,{\rm mod}\,8\,$, it follows that $\,\Gam = \Gam_{\wh{1}} \cong \Gam_{\wh{5}}\,$, while 
$\,\Gam_{\wh{7}} \cong \D_8 \cong \Z/4 \Z \rtimes \Z/2 \Z\,$, and $\,\Gam_{\wh{3}} \cong Q\,$ (the quaternion group). 
The triples $\,(\Gam, \Gam_{\wh{7}}, \eta_{\wh{7}})\,$ and $\,(\Gam, \Gam_{\wh{3}}, \eta_{\wh{3}})\,$ are coGalois 
\footnotemark[3], but the Kneser triples $\,(\Gam_{\wh{7}}, \Gam, \eta_{\wh{7}}^{- 1})\,$ and 
$\,(\Gam_{\wh{3}}, \Gam, \eta_{\wh{3}}^{- 1})\,$ are not. Note also that $\,(\D_8, Q)\,$ is a deformation pair (see
Remarks~\ref{Rems:mncG}(3)).     
\end{ex}

\subsection{Kneser ideals}
\label{Subsec:KneserIdeals}
Let $\,\Gam\,$ be a profinite group, $\,\fG\,$ a profinite $\,\Gam\,$-group, and $\,\eta: \Gam \lra \fG\,$ a continuous
g-cocycle, with $\,\De:= \Ker(\eta),\,\De':= \text{Fix}_\Gam(\fG)\,$. We may assume without loss that the
triple $\,(\Gam, \fG, \eta)\,$ is normalized, i.e., $\De \cap \De' = \{1\}\,$. Let $\,\L(\Gam\,|\,\De),\,\L(\fG)\,$ 
be the lattices connected through the operators 
$\,\cJ: \L(\Gam\,|\,\De) \lra \L(\fG),\,\cS: \L(\fG) \lra \L(\Gam\,|\,\De)\,$ defined in \ref{Subsec:GeneralFrame}.

\begin{de}
\label{Def:KneserIdeal}
$\,{\bf a} \in \L(\fG)\,$ is called a {\em Kneser ideal} (with respect to $\,\eta\,$) if the cocycle \\
$\,\eta_{\bf a}: \Gam \lra \fG/{\bf a}\,$ induced by $\,\eta\,$ is surjective, while $\,\fG/{\bf a}\,$ is called a 
{\em Kneser quotient} (with respect to $\,\eta$\,) of the profinite $\,\Gam$-group $\,\fG$.
\end{de}
We denote by $\,\K(\fG)\,$ the set of all Kneser ideals
of $\,\fG\,$, partially ordered by
inclusion, with $\,\fG\,$ as the last element, canonically anti-isomorphic to the
poset of all Kneser quotients of $\,\fG\,$.
Note that $\,\K(\fG)$ is an upper
subset of the lattice $\,\L(\fG)$, i.e., $\,\forall {\bf a} \in \K(\fG),\,\forall {\bf b} \in \L(\fG), {\bf a} \sse
{\bf b} \Lra {\bf b} \in \K(\fG)$. 

\begin{lem}
\label{Lem:Section}
If $\,(\Gam, \fG, \eta)\,$ is a Kneser triple then $\,\K(\fG) = \L(\fG)$, and the continuous map 
$\,\cS : \L(\fG) \lra \L(\Gam\,|\,\De)\,$ is a section of the coherent map $\,\cJ : \L(\Gam\,|\,\De) \lra \L(\fG)$.
\end{lem}
\bp
By assumption $\,\eta: \Gam \lra \fG\,$ is surjective, so $\,\{ 1 \} \in \K(\fG)$, therefore $\,\K(\fG) = \L(\fG)$ 
since $\,\K(\fG)\,$ is an upper subset of $\,\L(\fG)\,$. Let $\,{\bf a} \in \L(\fG)\,$ and put $\,\Lam := 
\cS({\bf a}) \in \L(\Gam\,|\,\De),\,{\bf b} := \cJ(\Lam) \sse {\bf a}\,$. As $\,\cS \circ \cJ \circ \cS = \cS\,$, 
we obtain $\,\cS({\bf b}) = \Lam\,$. On the other hand, the canonical maps 
$\,\eta'_{\bf a}: \Gam/\Lam \lra \fG/{\bf a}\,$ and $\,\eta'_{\bf b} : \Gam/\Lam \lra \fG/{\bf b}\,$ induced by 
$\,\eta\,$ are bijective. Since $\,\eta'_{\bf a}\,$ is the composition of $\,\eta'_{\bf b}\,$ with the natural 
projection $\,\fG/{\bf b} \lra \fG/{\bf a}$, it follows that $\,{\bf b} = {\bf a}\,$, hence $\cJ \circ \cS = 
1_{\L(\fG)}\,$ as desired.
\ep  
\begin{pr}
\label{Pr:KneserSpace}
Given a triple $\,(\Gam, \fG, \eta)\,$, where $\,\eta: \Gam \lra \fG\,$ is a g-cocycle, the following 
assertions hold.

{\rm (1)} Let $\,{\bf a} \in \L(\fG)\,$. Then $\,{\bf a} \in \K(\fG)\,$ if and only if $\,(\Gam : \cS({\bf b})) =
(\fG : {\bf b}) ( = |\fG/{\bf b}| )\,$ for all {\em open} ideals $\,{\bf b} \in \L(\fG\,|\,{\bf a})$.

{\rm (2)} $\K(\fG)\,$ is a closed subspace of the spectral space $\,\L(\fG)$.

{\rm (3)} For every Kneser ideal $\,{\bf a}\,$, the set of minimal Kneser ideals contained in $\,{\bf a}\,$
is nonempty.

{\rm (4)} The subspace $\,\K(\fG)_{\rm min}$ of the spectral space $\,\K(\fG)\,$ consisting of 
all minimal Kneser ideals of $\,\fG\,$ is Hausdorff. 
\end{pr} 
\bp
(1) follows by Lemma~\ref{Lem:SurCocycle},(2) and  Corollary~\ref{Cor:SurCocycle} applied to the surjective cocycle
$\,\eta_{\bf a}: \Gam \lra \fG/{\bf a}\,$ for any $\,{\bf a} \in \K(\fG)$. 

(2) Let $\,{\bf a} \in \L(\fG) \sm \K(\fG)\,$. By (1) there exists an open ideal $\,{\bf b} \in \L(\fG\,|\,{\bf a})\,$  
such that $\,{\bf b} \notin \K(\fG)\,$. Thus $\,\L({\bf b})\,$ is an open neighbourhood of $\,{\bf a}\,$ and 
$\,\L({\bf b}) \cap \K(\fG) = \es\,$, so $\,\K(\fG)\,$ is closed in the spectral space $\,\L(\fG)\,$; in particular, 
$\,\K(\fG)\,$ is a spectral space with respect to the induced topology.

(3) follows by (1) and Zorn's lemma.

(4) Let $\,{\bf a}_i \in \K(\fG)_{\rm min},\, i = 1, \,2\,$, be such that $\,{\bf a}_1 \neq {\bf a}_2\,$. Consequently, 
there exist open ideals $\,{\bf b}_i \in \L(\fG\,|\,{\bf a}_i),\, i = 1,\, 2\,$, such that 
$\,{\bf b}_1 \cap {\bf b}_2 \notin \K(\fG)\,$, since otherwise it would follow by compactness that each open ideal 
lying over $\,{\bf a}_1 \cap {\bf a}_2\,$ is Kneser, whence $\,{\bf a}_1 \cap {\bf a}_2 \in \cK(G)\,$ by (1), contrary 
to the minimality condition satisfied by the distinct Kneser ideals $\,{\bf a}_1\,$ and $\,{\bf a}_2\,$. For such a 
pair $\,({\bf b}_1, {\bf b}_2),\,\L({\bf b}_i)\,$ is an open neighbourhood of $\,{\bf a}_i\,$ for $\,i = 1,\,2\,$, and 
$\,\L({\bf b}_1) \cap \L({\bf b}_2) \cap \K(\fG) = \es\,$, so $\,\K(\fG)_{\rm min}\,$ is a Hausdorff space with 
respect to the topology induced from the spectral space $\,\K(\fG)\,$. 
\ep

As shown in Remark~\ref{Re:CoherentS}, the continuous map $\,\cS: \L(\fG) \lra \L(\Gam\,|\,\De)$ is not necessarily
a coherent map. Related to this fact, the following question arises naturally.
\begin{prob}
\label{Prob:CoherentSresK}
Is the continuous restriction map $\,\cS\,|_{\K(\fG)}: \K(\fG) \lra \L(\Gam\,|\,\De),\,{\bf a} \mapsto 
\eta^{- 1}({\bf a})\,$, necessarily coherent ? 
\end{prob}
A positive answer to the question above is given in \cite[Theorem 2.13]{Serdica} in the particular framework of
cyclotomic abstract coGalois theory presented in Example~\ref{Ex:1}.

\begin{res}
\label{Rems:KneserSbgps} \em
(1) Let $\,\K(\Gam)\,$ denote the image of the restriction map $\,\cS\,|_{\K(\fG)}: \K(\fG) \lra \L(\Gam\,|\,\De)\,$, 
and call its members {\em Kneser subgroups} of the profinite group $\,\Gam\,$. The map $\,\cS\,|_{\K(\fG)}\,$ is
not necessarily injective; moreover $\,\cS({\bf a}) = \cS({\bf b})\,$, with $\,{\bf a}, {\bf b} \in \K(\fG)\,$, does 
not imply $\,{\bf a} \cong {\bf b}\,$ as profinite groups (see \cite[Remarks 2.6 (1)]{Serdica} for a simple example).
Note also that, in general, $\,\K(\Gam)\,$ is not an upper subset of $\,\L(\Gam\,|\,\De)\,$, and hence not necessarily a 
closed subset of the spectral space $\,\L(\Gam\,|\,\De)\,$ (for an example see \cite[Remarks 2.6 (2)]{Serdica}); however
$\,\K(\Gam)\,$ is closed with respect to the profinite topology as image through a continuous map of $\,\K(\fG)\,$ with
its profinite topology.
 
(2) Let $\,\H\K(\Gam)\,$ denote the subset of $\,\K(\Gam)\,$ consisting of those $\,\Lam \in \L(\Gam\,|\,\De)\,$ for 
which $\,\L(\Gam\,|\,\Lam) \sse \K(\Gam)\,$, termed {\em hereditarily Kneser subgroups} of $\,\Gam\,$. In the
particular framework of cyclotomic abstract coGalois theory, $\,\H\K(\Gam)\,$ is a closed subspace of the spectral
space $\,\L(\Gam\,|\,\De)\,$ \cite[Corollary 2.15]{Serdica}, and an explicit {\em hereditarily-Kneser criterion} for 
closed subgroups of $\,\Gam\,$ is provided by \cite[Lemma 3.1, Theorem 3.2]{Serdica}; field theoretic interpretations 
are contained in \cite[Remark 2.16]{Serdica}.  
\end{res}
 
\subsection{A general Kneser criterion}
\label{Subsec:KneserCrit}
We end this section with a general ommitting-type criterion for an ideal to be Kneser and with an open problem concerning 
the classification of certain finite algebraic structures arising from this criterion.

Given a g-cocycle $\,\eta \in Z^1(\Gam, \fG)$, we denote by $\,\N\K(\fG)\,$ the complementary set of 
$\,\K(\fG)\,$ in the set $\,\L(\fG)\,$ of all ideals of the profinite $\,\Gam$-group $\,\fG\,$. By 
Proposition~\ref{Pr:KneserSpace},(2), $\,\N\K(\fG)$ is an open and hence a lower subset of the spectral space 
$\L(\fG)$. Note that $\,\N\K(\fG) \neq \es \Llra \eta(\Gam) \neq \fG\,$. We denote by $\,\N\K(\fG)_{\rm max}\,$ 
the set of all maximal members of the poset $\,\N\K(\fG)\,$. By Proposition~\ref{Pr:KneserSpace},(1), $\,{\bf m}\,$ is 
an open ideal of $\,\fG\,$ provided $\,{\bf m} \in \N\K(\fG)_{\rm max}\,$, and $\,\N\K(\fG)\,$ is the union of the 
basic open compact sets $\,\L({\bf m})\,$ for $\,{\bf m}\,$ ranging over $\,\N\K(\fG)_{\rm max}$. Consequently, we
obtain
\begin{pr} 
\label{Pr:AbsKneserCrit}
{\em (Abstract Kneser Criterion)} Let $\,\eta: \Gam \lra \fG\,$ be a g-cocycle. Then the following
assertions are equivalent for any ideal $\,{\bf a}\,$ of the profinite $\,\Gam$-group $\fG$. 

{\rm (1)} $\,{\bf a} \in \K(\fG)$, i.e., the induced cocycle $\,\eta_{\bf a}: \Gam \lra \fG/{\bf a}\,$ is
surjective.

{\rm (2)} $\,{\bf a} \not \sse {\bf m}\,$ for all $\,{\bf m} \in \N\K(\fG)_{\rm max}\,$, i.e.,
$\,\L(\fG\,|\,{\bf a}) \cap \N\K(\fG)_{\rm max} = \es$.
\end{pr}
The following class of finite algebraic structures arises naturally from the abstract Kneser criterion above, 
considering the quotient $\,\Gam$-groups $\,\fG/{\bf m}\,$ with $\,{\bf m} \in \N\K(\fG)_{\rm max}$.
\begin{de}
\label{De:mnK}
A triple $\,(\Gam, \fG, \eta)\,$ consisting of a finite group $\Gam\,$, a finite $\,\Gam$-group $\fG$ and a 
g-cocycle $\,\eta: \Gam \lra \fG\,$, with $\,\De := \Ker(\eta)\,$, is called a {\em minimal non-Kneser triple} 
(for short {\em mnK triple}) if the following conditions are satisfied.

{\rm (1)} $\,\eta(\Gam) \neq \fG\,$, i.e., the cocycle $\,\eta\,$ is not surjective,

{\rm (2)} for every ideal $\,{\bf a} \neq \{ 1 \}\,$ of the $\,\Gam$-group 
$\,\fG\,,(\Gam : \eta^{- 1}({\bf a})) = (\fG : {\bf a})$, i.e., the induced cocycle 
$\,\eta_{\bf a}: \Gam \lra \fG/{\bf a}$ is surjective, and  

{\rm (3)} $\,\bigcap_{\gam \in \Gam} \gam \De \gam^{- 1} = {\rm Fix}_\Gam(G) \cap \De = \{ 1 \}$, i.e., 
the triple $\,(\Gam, \fG, \eta)\,$ is normalized.
\end{de}  
Note that the conjunction of conditions (1) and (2) above is equivalent with the sentence 
$\,\N\K(\fG) = \N\K(\fG)_{\bf max} = \{1\}$. 

In its full generality, the following problem is far from being trivial. 
\begin{prob}
\label{Prob:mnK}
Classify (up to isomorphism in $\,\cZ^1\,$) the minimal non-Kneser triples.
\end{prob}
Partial answers are given in Section~\ref{Sec:mnK}.

\section{CoGalois and minimal non-coGalois triples}

In this section we introduce two remarkable classes of surjective cocycles extending to a more 
general framework some notions and results from \cite{A-B} on {\em coGalois groups of cocycles}. 
\begin{de}
\label{Def:coGTriple}
A triple $\,(\Gam, \fG, \eta)\,$ is called {\em coGalois} if the cocycle $\,\eta: \Gam \lra \fG\,$ is surjective and
the maps $\,\cJ : \L(\Gam\,|\,\Ker(\eta)) \lra \L(\fG)\,$ and $\,\cS : \L(G) \lra \L(\Gam\,|\,\Ker(\eta))\,$ are 
lattice isomorphisms inverse to one another, i.e., the coGalois connexion between the lattices 
$\,\L(\Gam\,|\,\Ker(\eta))\,$ and $\,\L(\fG)\,$ is perfect.
\end{de}
Characterizations of coGalois triples are given by the next lemma whose proof is straightforward.
\begin{lem}
\label{Lem:coGtriple}
Let $\,\eta \in Z^1(\Gam, \fG)\,$ be a g-cocycle with $\,\De:= \Ker(\eta)\,$. The following assertions are equivalent.

{\rm (1)} $\,(\Gam, \fG, \eta)\,$ is a coGalois triple.

{\rm (2)} $\,\eta\,$ is surjective and $\,\cS \circ \cJ = 1_{\L(\Gam\,|\,\De)}$.

{\rm (3)} $\,\eta\,$ is surjective and $\,\cJ : \L(\Gam\,|\,\De) \lra \L(\fG)\,$ is injective. 

{\rm (4)} $\,\eta\,$ and $\,\cS : \L(\fG) \lra \L(\Gam\,|\,\De)\,$ are surjective.

{\rm (5)} $\,\eta(\Lam) = \cJ(\Lam)\,$ for all $\,\Lam \in \L(\Gam\,|\,\De)$.

{\rm (6)} $\,\eta(\Lam) \in \L(\fG)\,$ for all $\,\Lam \in \L(\Gam\,|\,\De)$.
\end{lem}

\begin{de}
Let $\,\eta \in Z^1(\Gam, \fG)\,$ be a g-cocycle. An ideal $\,{\bf a} \in \L(\fG)\,$ is called {\em coGalois} if 
the induced triple $\,(\Gam,\,\fG/{\bf a},\,\eta_{\bf a}: \Gam \lra \fG/{\bf a})\,$ is coGalois.
\end{de}
We denote by $\,\mathbb{C G}(\fG)\,$ the poset with the last element $\,\fG\,$ of all coGalois ideals of $\,\fG\,$. Some
properties of coGalois ideals are collected together in the next result.
\begin{pr}
\label{Pr:coGaloisSpace}
Let $\,\eta \in Z^1(\Gam, \fG)\,$ be a g-cocycle. The following assertions hold.

{\rm (1)} $\,\mathbb{C G}(\fG)\,$ is an upper subset of $\L(\fG)\,$, contained in $\,\K(\fG)$. 

{\rm (2)} For any ideal $\,{\bf a} \in \L(\fG)\,$, $\,{\bf a} \in \mathbb{C G}(\fG)\,$ if and only if 
$\,{\bf b} \in \mathbb{C G}(\fG)\,$ for all open ideals $\,{\bf b}\,$ containing $\,{\bf a}$.

{\rm (3)} $\,\mathbb{C G}(\fG)\,$ is a closed subspace of the spectral space $\,\K(\fG)$.

{\rm (4)} For every $\,{\bf a} \in \mathbb{C G}(\fG)\,$ there exists at least one minimal coGalois ideal
$\,{\bf b} \sse {\bf a}$.

{\rm (5)} The space $\,\mathbb{C G}(\fG)_{\rm min}\,$ of all minimal members of $\,\mathbb{C G}(\fG)\,$ is 
Hausdorff with respect to the topology induced from the spectral space $\,\mathbb{C G}(\fG)$.
\end{pr}
\bp
The proof is similar with the proof of Proposition~\ref{Pr:KneserSpace}.
\ep

\begin{rem} \em 
The restriction map $\,\cS\,|_{\C\mathbb{G}(\fG)}: \C\mathbb{G}(\fG) \lra \L(\Gam\,|\,\De)\,$ is injective and coherent,
inducing a homeomorphism of the spectral space $\,\C\mathbb{G}(\fG)\,$ onto the closed subspace 
$\,\C\mathbb{G}(\Gam):= \cS(\C\mathbb{G}(\fG))\,$ of the spectral space $\,\L(\Gam\,|\,\De)\,$, contained in the
subspace $\,\H\K(\Gam)\,$ defined in Remarks~\ref{Rems:KneserSbgps}; call its members {\em coGalois subgroups} of 
the profinite group $\,\Gam\,$. An explicit {\em coGalois criterion for hereditarily Kneser subgroups} of $\,\Gam\,$
is provided by \cite[Theorem 3.2]{JPAA} in the particular framework of cyclotomic abstract coGalois theory.  
\end{rem}
Given a surjective cocycle $\,\eta \in Z^1(\Gam, \fG)\,$, we denote by $\,\mathbb{N C G}(\fG)\,$ the complementary set 
of $\,\mathbb{C G}(\fG)\,$ in the set $\,\L(\fG) = \K(\fG)\,$ of all ideals of the profinite $\,\Gam$-group $\,\fG\,$. 
By Proposition~\ref{Pr:coGaloisSpace},(3), $\,\mathbb{N C G}(\fG)\,$ is an open and hence a lower subset of the spectral space 
$\L(\fG)$. By Lemma~\ref{Lem:coGtriple},(6), $\,\mathbb{N C G}(\fG) \neq \es \Llra \eta(\Lam) \notin \L(\fG)\,$ for
some $\,\Lam \in \L(\Gam\,|\,\Ker(\eta)) \sm \{\Ker(\eta), \Gam\}\,$. We denote by $\,\mathbb{N C G}(\fG)_{\rm max}\,$ 
the set of all maximal members of the poset $\,\mathbb{N C G}(\fG)\,$. By Proposition~\ref{Pr:coGaloisSpace},(4), 
$\,{\bf m}\,$ is an open ideal of $\,\fG\,$ provided $\,{\bf m} \in \mathbb{N C G}(\fG)_{\rm max}\,$, and 
$\,\mathbb{N C G}(\fG)\,$ is the union of the basic open compact sets $\,\L({\bf m})\,$ for $\,{\bf m}\,$ ranging 
over $\,\mathbb{N C G}(\fG)_{\rm max}$. Consequently, we obtain the following analogue of 
Proposition~\ref{Pr:AbsKneserCrit}.
\begin{pr}
\label{Pr:AbscoGaloisCrit}
{\rm (Abstract coGalois Criterion)} Let $\,\eta \in Z^1(\Gam, \fG)\,$ be a surjective cocycle. The following assertions 
are equivalent for any $\,{\bf a} \in \L(\fG) = \K(\fG)$.

{\rm (1)} $\,{\bf a} \in \mathbb{C G}(\fG)\,$, i.e., the induced triple 
$\,(\Gam,\,\fG/{\bf a},\,\eta_{\bf a}: \Gam \lra \fG/{\bf a})\,$ is coGalois.

{\rm (2)} $\,{\bf a} \not \sse {\bf m}\,$ for all $\,{\bf m} \in \mathbb{N C G}(\fG)_{\rm max}\,$, i.e.,
$\,\L(\fG\,|\,{\bf a}) \cap \mathbb{N C G}(\fG)_{\rm max} = \es$.
\end{pr} 
Considering the quotient $\,\Gam$-groups $\,\fG/{\bf m}\,$ for $\,{\bf m} \in \mathbb{N C G}(\fG)_{\rm max}$,
we obtain the following class of finite algebraic structures.
\begin{de}
\label{De:mncG}
A triple $\,(\Gam, \fG, \eta)\,$ consisting of a finite group $\Gam\,$, a finite $\,\Gam$-group $\fG$ and a 
surjective cocycle $\,\eta: \Gam \lra \fG\,$, with $\,\De := \Ker(\eta)\,$, is called a (Kneser) 
{\em minimal non-coGalois triple} 
(for short {\em mncG triple}) if the following conditions are satisfied.

{\rm (1)} $\,\eta(\Lam)\,$ is not an ideal of $\,\fG\,$ for some $\,\Lam \in \L(\Gam\,|\,\De) \sm \{\De, \Gam\}\,$,
i.e., $\,(\Gam, \fG, \eta)\,$ is not a coGalois triple,

{\rm (2)} for every ideal $\,{\bf a} \neq \{ 1 \}\,$ of $\,\fG\,$, the induced triple 
$\,(\Gam,\,\fG/{\bf a},\,\eta_{\bf a}: \Gam \lra \fG/{\bf a})\,$ is coGalois, and  

{\rm (3)} $\,\bigcap_{\gam \in \Gam} \gam \De \gam^{- 1} = {\rm Fix}_\Gam(G) \cap \De = \{ 1 \}$, i.e., 
the triple $\,(\Gam, \fG, \eta)\,$ is normalized.
\end{de}  
Note that the conjunction of conditions (1) and (2) above is equivalent with the sentence 
$\,\mathbb{N C G}(\fG) = \mathbb{N C G}(\fG)_{\bf max} = \{1\}$. 

\begin{prob}
\label{Prob:mncG}
Classify (up to isomorphism in $\,\cZ^1\,$) the (Kneser) minimal non-coGalois triples.
\end{prob}

\begin{res} 
\label{Rems:mncG} \em
(1) Problem~\ref{Prob:mncG} is solved in \cite{A-B} in the particular framework described in Example~\ref{Ex:1}. 
Assume $\,\fG\,$ is a finite $\,\Gam$-module of exponent $\,k\,$, and the action of $\,\Gam\,$ is given by
a character $\,\chi: \Gam \lra (\Z/k \Z)^\times\,$. Let $\,\eta \in Z^1(\Gam, \fG)\,$. According to \cite[Lemma 2.17,
Corollary 2.18]{A-B}, $\,(\Gam, \fG, \eta)\,$ is a Kneser mncG triple if and only if the triple $\,(\Gam, \fG, \eta)\,$ 
is, up to isomorphism, of one of the following three types.

(i) $\,k = 4\,$, $\,\Gam = \la \si, \tau\,|\,\si^2 = \tau^2 = (\si \tau)^2 = 1\ra \cong 
\Z/2 \Z \times \Z/2 \Z\,$, $\,\fG = \Z/4 \Z\,$, $\,\chi(\si) = - 1\,{\rm mod}\,4$, $\,\chi(\tau) = 1\,{\rm mod}\,4\,$, 
$\,\eta(\si) = 1\,{\rm mod}\,4, \eta(\tau) = 2\,{\rm mod}\,4$.

(ii) $\,k = 4\,$, $\,\Gam = \D_8 = \la \si, \tau\,|\,\si^2 = \tau^4 = (\si \tau)^2 = 1 \ra \cong 
\Z/4 \Z \rtimes \Z/2 \Z\,$, $\,\fG = (\Z/2 \Z) e_1 \oplus (\Z/4 \Z) e_2$, $\,\chi(\si) = - 1\,{\rm mod}\,4\,$,
$\,\chi(\tau) = 1\,{\rm mod}\,4\,$, $\,\eta(\si) = e_2$, $\,\eta(\tau) = e_1 + e_2$.

(iii) $\,k = p r\,$ with $\,p \neq 2\,$ prime, $\,1 \neq r\,|\,(p - 1)\,$, $\,\Gam = \la \si, \tau\,|\,
\si^r = \tau^p = \si \tau \si^{- 1} \tau^{- u} = 1 \ra \cong \Z/p \Z \rtimes_u \Z/r \Z\,$, where 
$\,u \in (\Z/p r \Z)^\times\,$ such that the order of $\,u\,{\rm mod}\,p \in (\Z/p \Z)^\times\,$ is $\,r\,$ and
$\,l\,|\,(u - 1)\,$ for all $\,l\,|\,r\,$ with $\,l \neq 2\,$ prime or $\,l = 4\,$, $\,\fG = \Z/p r \Z\,$,
$\,\chi(\si) = u, \chi(\tau) = 1\,{\rm mod}\,p r\,$, $\eta(\si) = p\,{\rm mod}\,p r,\,\eta(\tau) = r\,{\rm mod}\,pr$.

Based on the classification of mncG triples above, \cite[Theorem 2.19]{A-B} provides an explicit form of 
Proposition~\ref{Pr:AbscoGaloisCrit} in the framework of cyclotomic abstract coGalois theory.

(2) The bijective cocycle $\,\Z/2 \Z \times \Z/2 \Z \lra \Z/4 \Z\,$ from (1)(i) is induced by an action of
$\,\Z/2 \Z \times \Z/2 \Z\,$ on itself, while the bijective cocycle $\,\D_8 \lra \Z/2 \Z \oplus \Z/4 \Z\,$ from 
(1)(ii) is not since $\,(\D_8, \Z/2 \Z \oplus \Z/4 \Z)\,$ is not a deformation pair.

(3) The dihedral group $\,\D_8 =\Z/4 \Z \rtimes \Z/2 \Z \cong \la \si, \tau\,|\,\si^2 = \tau^4 = (\si \tau)^2 = 1 \ra\,$
and the quaternion group
$\,Q \cong \la \rho, \theta\,|\,\rho^4 = 1, \rho^2 = \theta^2, \rho \theta \rho^{- 1} = \theta^{- 1} \ra\,$ form a 
deformation pair: consider the actions $\,\D_8 \times Q \lra Q\,$, $\,Q \times \D_8 \lra D_8\,$, defined by
\[
^\si \rho = \rho^{- 1},\, ^\si \theta = \theta,\, ^\tau \rho = \rho,\, ^\tau \theta = \theta;\,^\rho \si = \tau^2 \si,\,
^\rho \tau = \tau,\,^\theta \si = \si,\,^\theta \tau = \tau.
\]
Setting $\,\eta(\si) = \rho,\,\eta(\tau) = \theta\,$, we obtain a bijective cocycle $\,\eta \in Z^1(\D_8, Q)\,$ such
that $\,\eta^{- 1} \in Z^1(Q, \D_8)\,$ as desired. Since $\,\D_8/C(\D_8) \cong Q/C(Q) \cong \Z/2 \Z \times \Z/2 \Z\,$,
it follows easily that the Kneser triples $\,(\D_8, Q, \eta)\,$ and $\,(Q, \D_8, \eta^{- 1})\,$ are both mncG, providing
examples of Kneser mncG triples in a purely noncommutative framework. 

Similarly, for any $\,m \geq 3\,$, 
the dihedral group $\,\D_{2^m} \cong \la \si, \tau\,|\,\si^2 = \tau^{2^{m - 1}} = (\si \tau)^2 = 1 \ra\,$ and the
generalized quaternion group $\,Q_{2^m} \cong \la \rho, \theta\,|\,\rho^4 = 1, \rho^2 = \theta^{2^{m - 2}},
\rho \theta \rho^{- 1} = \theta^{- 1} \ra\,$ form a deformation pair inducing Kneser non-coGalois triples, but the
minimality condition is satisfied only for $\,m = 3$.

(4) Other remarkable Kneser mncG triples in a purely noncommutative framework are obtained as follows. Let $\,K = \F_q\,$
be a finite field, $\,q = p^f, f \geq 1, p \neq 2\,$, $\,L = \F_{q^2}\,$, $\,{\rm Gal}(L/K) = \la \alpha \ra \cong
\Z/2 \Z\,$, with $\,\alpha(x) = x^q\,$, and $\,\chi: L^\times/(L^\times)^2 \lra \Z/2 \Z\,$ the unique
isomorphism. Setting $\,\Gam:= PGL(2, q^2) = PGL(2, L)\,$, extend $\,\alpha\,$ to an automorphism of $\,\Gam\,$, and
consider the self-action $\,\om: \Gam \lra {\rm Aut}(\Gam)\,$, defined by $\,\om(A\,{\rm mod}\,L^\times):= 
\alpha^{\chi({\rm det}(A) \cdot (L^\times)^2)}\,$ for $\,A \in GL(2, L)\,$. It follows that 
$\,\om(\Gam) = \la \alpha \ra \cong \Z/2 \Z\,$, and $\,\Ker(\om) = PSL(2, L) = PSL(2, q^2)\,$. The self-action $\,\om\,$
is adequate, and the deformation of $\,\Gam\,$ via $\,\om\,$ is the {\em Zassenhaus group} $\,\Gam_\om = M(q^2)\,$. 
The non-isomorphic finite groups $\,\Gam\,$ and $\,\Gam_\om\,$ of the same order $\,(q^2 - 1) q^2 (q^2 + 1)\,$, having 
the {\em simple} group $\,\Ker(\om) = PSL(2, q^2)\,$ as a common (normal) subgroup of index $\,2\,$, act both 
{\em faithfully} and {\em sharply $\,3$-transitive} on the projective line $\,\P^1(\F_{q^2})\,$. The induced Kneser
triples $\,(\Gam, \Gam_\om, \eta_\om)\,$ and $\,(\Gam_\om, \Gam, \eta_\om^{- 1})\,$ are both mncG.                                                                                                                                                                                                                                                                                        
\end{res}

\section{Partial answers to classification Problem~\ref{Prob:mnK}}
\label{Sec:mnK}

This last section, devoted to partial answers to Problem~\ref{Prob:mnK}, provides some classes of mnK triples 
including the very simple ones which occur in \cite[Lemma 1.18, Theorem 1.20]{A-B} in the framework of cyclotomic 
abstract coGalois theory described in Example~\ref{Ex:1}.

First some useful lemmas. A proper subclass of mnK triples is provided by the next obvious lemma.
\begin{lem}
\label{Lem:Simple}
Let $\,\eta \in Z^1(\Gam, \fG)\,$ be a normalized g-cocycle such that $\,\eta(\Gam) \neq \fG\,$. Then 
$\,(\Gam, \fG, \eta)\,$ is a mnK triple whenever the $\,\Gam$-group $\,\fG\,$ is {\em simple}, i.e.,
$\,\{ 1 \}\,$ and $\,\fG\,$ are its only ideals.
\end{lem}

\begin{re} \em
Not all mnK triples are of the type above. The simplest example of a mnK triple  $\,(\Gam, \fG, \eta)\,)\,$, where
the $\,\Gam$-group $\,\fG\,$ is not simple, is obtained by taking $\,\Gam := \Z/2\Z,\,\fG := \Z/4\Z\,$ with the 
unique non-trivial action of $\,\Gam\,$ on $\,\fG\,$. Then $\,Z^1(\Gam, \fG) \cong \Z/4\Z$ is generated by the 
injective g-cocycle $\,\eta\,$ defined by $\,\eta(1\,{\rm mod}\,2) = 1\,{\rm mod}\,4\,$. Note that 
$\,\eta(\Gam) = \{ 0\,{\rm mod}\,4,\,1\,{\rm mod}\,4 \} \neq \fG$. As $\,0,\,\fG$ and 
$\,{\bf a}:= \{ 0\,{\rm mod}\,4,\,2\,{\rm mod}\,4 \} \cong \Z/2\Z$ are the only submodules of $\,\fG$, and the 
induced cocycle $\,\eta_{\bf a}: \Gam \lra \fG/{\bf a} \cong \Z/2\Z$ is an isomorphism, it follows that 
$\,(\Gam, \fG, \eta)\,$ is a mnK triple but the $\,\Gam$-module $\,\fG\,$ is not simple. 
\end{re}
\begin{co}
Assume $\,\Gam = \fG\,$ is a non-abelian finite simple group acting on itself by inner automorphisms 
$\,(\gam, g) \mapsto \gam g \gam^{- 1}$. The following assertions hold.

{\rm (1)} Every non-trivial cocycle $\,\eta \in Z^1(\Gam, G)\,$ is a normalized g-cocycle. 

{\rm (2)} $\,(\Gam, \fG, \eta_g)\,$ is a mnK triple, where $\,\eta_g : \Gam \lra \fG\,$ is the coboundary 
$\,\gam \mapsto [g, \gam] := g \gam g^{- 1} \gam^{- 1}\,$ for any $\,1 \neq g \in \fG\,$, with 
$\,\Ker(\eta_g) = C_\Gam(g)\,$, the centralizer of $\,g$.

{\rm (3)} For $\,g,\,h \in \fG \sm \{ 1 \}\,$, the mnK triples $\,(\Gam, \fG, \eta_g)\,$ and 
$\,(\Gam, \fG, \eta_h)\,$ are isomorphic if and only if $\,\varphi(g) = h\,$ for some automorphism $\,\varphi\,$ 
of $\,\fG = \Gam$. 
\end{co}
\bp
(1) For a non-trivial cocycle $\,\eta \in Z^1(\Gam, G)\,$, let $\,\fH\,$ denote the subgroup of $\,\fG\,$ generated 
by $\,\eta(\Gam)\,$. As the cocycle $\,\eta\,$ is non-trivial, $\,\fH\,$ is a proper subgroup of $\,\fG\,$. Since 
$\,\si \eta(\tau) \si^{- 1} = \eta(\si)^{- 1} \eta(\si \tau)\,\in \fH\,$ for all $\,\si, \tau \in \Gam\,$, it follows 
that $\,\fH\,$ is a normal subgroup of $\,\fG\,$, and hence $\,\fH = \fG = \Gam\,$. Consequently, $\,\eta\,$ is a 
g-cocycle. It is also normalized since $\,\Ker(\eta) \neq \Gam\,$ implies 
$\,\dbigcap_{\gam \in \Gam} \gam \Ker(\eta) \gam^{- 1} = \{ 1 \}$ as the group $\,\Gam\,$ is simple.

(2) By (1), $\,\eta_g\,$ is a normalized g-cocycle. On the other hand, $\,1 \neq g \in C_\Gam(g) = \Ker(\eta_g)\,$,
therefore $\,\eta_g(\Gam) \neq \fG\,$ since $\,\Gam = \fG\,$ is finite. The conclusion is immediate by 
Lemma~\ref{Lem:Simple}.

(3) follows from the definition of isomorphic cocycles.
\ep
\smallskip

\begin{lem}
\label{Lem:Nilpotent}
Let $\,(\Gam, \fG, \eta)\,$ be a mnK triple, $\,\De:= \Ker(\eta)\,$. Then the following assertions hold.

{\rm (1)} $\,{\rm Fix}\,_\Gam(\fG) \cap \eta^{- 1}(C(\fG)) = \{ 1 \}$.

{\rm (2)} If $\,\fG\,$ is abelian then $\,{\rm Fix}\,_\Gam(\fG) = \{ 1 \}\,$, i.e., $\,\Gam\,$ acts faithfully 
on $\,\fG$.

{\rm (3)} If $\,\fG\,$ is nilpotent then it is a $\,p$-group for some prime number $\,p\,$. If
$\,p \notdiv (\Gam : \De)\,$, then $\,\fG\,$ is a simple $\,\F_p[\Gam]$-module, the action of $\,\Gam\,$ is
faithful, and $\,\eta\,$ is a coboundary.
\end{lem}
\bp
(1) Letting $\,\overline{\De}:= {\rm Fix}\,_\Gam(\fG) \cap \eta^{- 1}(C(\fG))\,$, it suffices to show that 
$\,\overline{\De} \sse \De\,$ since $\,\De \cap {\rm Fix}\,_\Gam(\fG) = \{ 1 \}\,$ by assumption. By 
Lemma~\ref{Lem:PropGenCocycles},(3, 4), $\,{\bf a} := \eta(\overline{\De})\,$ is an ideal of the $\,\Gam$-group 
$\,\fG\,$. On the other hand, since $\,\eta\,$ is not surjective by assumption, it follows by 
Corollary~\ref{Cor:SurCrit} that the induced cocycle $\,\eta_{\bf a}: \Gam \lra \fG/{\bf a}\,$ is not surjective, 
so $\,{\bf a} = \{ 1 \}\,$ by the minimality property of $\,\eta$, and hence $\,\overline{\De} \sse \De\,$ as desired.

(2) is an immediate consequence of (1).

(3) Since the cocycle $\,\eta : \Gam \lra \fG\,$ is not surjective and $\,\fG\,$ is nilpotent by assumption, it
follows by Proposition~\ref{Pr:Pronil} that there exists a prime number $\,p\,$ such that the induced cocycle
$\,\eta_p: \Gam \lra \fG_p\,$ is not surjective, therefore the kernel of the natural projection $\,\fG \lra \fG_p\,$ 
is trivial by the minimality property of $\,\eta\,$, and hence $\,\fG \cong \fG_p\,$ as required.

Assuming $\,p \notdiv (\Gam : \De)$, it follows that the $\,\Gam$-group $\,\fG\,$ is simple. In particular, since 
the center $\,C(\fG)\,$ of the $\,p$-group $\,\fG\,$ is a nontrivial ideal of the $\,\Gam$-group $\,\fG\,$, we obtain 
$\,\fG = C(\fG)$, so $\,\fG\,$ is an abelian $\,p$-group. On the other hand, $\,p\fG = 0\,$ since $\,p\fG \neq \fG\,$. 
Thus $\,\fG\,$ is a simple $\,\F_p[\Gam]$-module as desired. The faithfulness of the action of $\,\Gam\,$ follows by 
(2), while $\,\eta(\gam) = (\gam - 1)g\,$ for $\,\gam \in \Gam\,$, with 
$\,g = - \frac{1}{(\Gam : \De)} \dsum_{\si \in \Gam/\De} \eta(\si)\,$.
\ep
The proof of the next lemma is straightforward.
\begin{lem}
Let $\,p\,$ be a prime number, $\,\Gam\,$ a non-trivial finite group, and $\,\fG\,$ a simple $\,\F_p[\Gam]$-module 
such that the action of $\,\Gam\,$ on $\,\fG\,$ is faithful. Let $\,K \cong \F_q,\,q = p^f,\,f \geq 1\,$, denote 
the field of the endomorphisms of the simple $\,\F_p[\Gam]$-module $\,\fG\,$. The following
assertions hold.

{\rm (1)} $\,\fG\,$ is a simple $\,K[\Gam]$-module, $\,{\rm End}\,_{K[\Gam]}(\fG) = K\,$, the group $\,\Gam\,$, 
identified with a subgroup of the linear group $\,GL_K(\fG)\,$, is not a $\,p$-group, and 
$\,C(\Gam) = \Gam \cap K^\times \cong \Z/r\Z, \,r\,|\,(q - 1)\,$. 

{\rm (2)} Every non-trivial cocycle $\,\eta: \Gam \lra \fG\,$ is a normalized g-cocycle. 

{\rm (3)} $\,(\Gam, \fG, \eta)\,$ is a mnK triple provided $\,\eta\,$ is a non-trivial coboundary.

{\rm (4)} Every cocycle $\,\eta: \Gam \lra \fG\,$ for which $\,p \not |\,(\Gam : \Ker(\eta))\,$ is a coboundary.
\end{lem}

\subsection{The abelian case of Problem~\ref{Prob:mnK}}
\label{Subsec:Abelian}
Now we approach the classification Problem~\ref{Prob:mnK} in the special case when $\,\Gam\,$ and $\,\fG\,$ are abelian.
\begin{lem}
\label{Lem:Ab1}
Let $\,(\Gam, \fG, \eta)\,$ be a mnK triple, with abelian $\,\Gam\,$ and $\,\fG\,$. Then the following assertions hold.

{\rm (1)} $\,\fG\,$ is an abelian $\,p$-group for some prime number $\,p,\,\Gam\,$ acts faithfully on $\,\fG\,$, 
and the g-cocycle $\,\eta\,$ is injective. Let $\,p^n, n \geq 1\,$, be the exponent of $\,\fG$.

{\rm (2)} $\,\Gam\,$ is a $\,p$-group $\,\Llra \fG^\Gam \neq 0 \Llra \fG^\Gam \cong \Z/p\Z \Llra |\fG| =
p |\Gam|$. 

{\rm (3)} There exists a unique minimal non-zero $\,\Gam$-submodule of $\,\fG$.

{\rm (4)} The image $\,R\,$ of the canonical ring homomorphism $\,(\Z/p^n\Z)\,[\Gam] \lra {\rm End}(\fG)\,$ is a 
finite commutative local ring of characteristic $\,p^n\,$, and $\,\Gam\,$, identified with a subgroup of $\,R^\times\,$,
generates $\,R\,$ as $\,(\Z/p^n\Z)$-module, in particular as ring.

{\rm (5)} If $\,\Gam\,$ is a $\,p$-group then the residue field $\,k := R/\fm \cong \F_p\,$, $\,\Gam = 1 + \fm\,$, 
and the $\,R$-module $\,\fG\,$ is (non-canonically) isomorphic to the Pontryagin dual
$\,R^\vee = {\rm Hom}(R^+, 1/p^n\,\Z/\Z)\,$ with the induced structure of $\,R$-module. In particular, if the local 
ring $\,R\,$ is principal then the $\,R$-module $\,\fG\,$ is free of rank $\,1$.      
\end{lem}
\bp
The assertion (1) is immediate from Lemma~\ref{Lem:Nilpotent}.

(2) We have to show that $\,|\fG| = p |\Gam|$ and $\,\fG^\Gam \cong \Z/p\Z\,$ provided $\,\fG^\Gam \neq 0\,$.
Put $\,\Lam:= \eta^{- 1}(\fG^\Gam)\,$. Since $\,(\Gam, \fG, \eta)$ is a mnK triple, it follows 
$\,|\fG| = (\Gam : \Lam) |\fG^\Gam|\,$. Assuming $\,\Lam \neq 1\,$, the injective cocycle $\,\eta\,$ induces 
by restriction a nontrivial monomorphism $\,\eta|_\Lam: \Lam \lra \fG^\Gam\,$ whose image 
$\,\fM:= \eta(\Lam) \cong \Lam\,$ is a non-zero $\,\Gam$-submodule of $\,\fG\,$, therefore the induced
cocycle $\,\eta_\fM \in Z^1(\Gam, \fG/\fM)\,$ is surjective since $\,(\Gam, \fG, \eta)\,$ is a mnK triple.
By Lemma~\ref{Lem:SurCrit} we deduce that the cocycle $\,\eta\,$ is surjective, whence a contradiction. Thus 
$\,\Lam = 1\,$, so it remains to show that $\,\fG^\Gam \cong \Z/p\Z\,$.

Choose an element $\,h \in \fG^\Gam\,$ of order $\,p\,$ and denote by $\,\fH \cong \Z/p\Z\,$ the 
$\,\Gam$-submodule generated by $\,h\,$. As $\,\fH \sse \fG^\Gam\,$, we obtain 
$\,\eta^{- 1}(\fH) \sse \eta^{- 1}(\fG^\Gam) = 1\,$, and hence
\[
(\fG : \fG^\Gam) = (\Gam : \eta^{- 1}(\fG^\Gam))
= |\Gam| = (\Gam : \eta^{- 1}(\fH)) = (\fG : \fH),
\]
so $\,\fG^\Gam = \fH \cong \Z/p\Z\,$ as desired. 

(3) If $\,\Gam\,$ is a $\,p$-group then $\,\fG^\Gam \cong \Z/p \Z\,$ (by (2)) is obviously the
unique minimal non-zero $\,\Gam$-submodule of $\,\fG\,$. Assume $\,\Gam\,$ is not a $\,p$-group and 
$\,\fH_i, i = 1, 2\,$, are distinct minimal non-zero $\,\Gam$-submodules of $\,\fG\,$, so $\,\fH_1 \cap \fH_2 = 0\,$. 
Put $\,\Lam_i := \eta^{- 1}(\fH_i), i = 1, 2\,$. Since $\,(\Gam, \fG, \eta)\,$ is mnK, it follows that 
$\,(\Gam : \Lam_i) = (\fG : \fH_i)\,$ is a $\,p$-th power and hence $\,\Lam_i \neq 1, i = 1, 2\,$, as 
$\,\Gam\,$ is not a $\,p$-group. On the other hand, 
$\,\Lam_1 \cap \Lam_2 = \eta^{- 1}(\fH_1 \cap \fH_2) = \eta^{- 1}(0) = \Ker(\eta) = 1\,$, therefore the
inclusion $\,\Lam_1 \hra \Gam\,$ induces a monomorphism $\,\Lam_1 \lra \Gam/\Lam_2\,$, so 
$\,\Lam_1\,$ is a $\,p$-group and hence $\,|\Gam| = |\Lam_1| (\Gam : \Lam_1)\,$ is a $\,p$-th
power, which is a contradiction.
  
(4) We have only to show that the finite commutative ring $\,R\,$ is local, i.e., the only idempotents of
$\,R\,$ are the elements $\,0\,$ and $\,1\,$. Let $\,e_i,\,i = 1, \dots, s\,$, denote the minimal non-zero
idempotents of $\,R\,$, the atoms of the boolean algebra $\,B(R) := \{e \in R\,|\,e^2 = e\}\,$ of  
idempotents of $\,R\,$ with respect to the partial order $\,e \leq f \Llra e f = e\,$. Put 
$\,R_i := R e_i\,$, $\,\fG_i := R_i \fG = e_i \fG, i = 1, \dots, s\,$. Then $\,\dsum_{i = 1}^{s}{e_i} = 1, e_i e_j = 0\,$ 
for $\,i \neq j\,$, the $\,R_i\,$'s are local rings, $\,B(R_i) = \{0, 1_{R_i}:= e_i\},\,R \cong \dprod_{i = 1}^{s} R_i\,$
is a semi-local ring, and $\,R = \bigoplus_{1 \leq i \leq s} R_i, \fG = \bigoplus_{1 \leq i \leq s} \fG_i\,$ are 
$\,R$-module direct sums. Note that $\,\fG_i \neq 0\,$ for $\,i = 1, \dots, s\,$, since 
$\,0 \neq e_i \in R \sse {\rm End}(\fG)\,$. By (3), we conclude that $\,s = 1\,$, i.e., $\,R\,$ is a local
ring as desired. 

(5) Let $\,h\,$ be a generator of $\,\fH:= \fG^\Gam \cong \Z/p\Z\,$, the unique minimal non-zero $\,R$-submodule of 
$\,\fG\,$. The kernel of the surjective homomorphism of $\,R$-modules $\,R \lra \fH\,$, $\lam \mapsto \lam h\,$, is 
the maximal ideal $\,\fm\,$ of $\,R\,$, so $\,k:= R/\fm \cong (\Z/p^n\Z)/p(\Z/p^n\Z) \cong \F_p\,$ and 
$\,\Gam \sse 1 + \fm\,$. 

Since $\,|\fG| = |\Gam| \cdot p\,\,|\,|1 + \fm| \cdot p = |\fm| \cdot|R/\fm| = |R|\,$,
we deduce that $|\fG| \leq |R|\,$. Consider the Pontryagin dual $\,\fG^\vee = {\rm Hom}(\fG, 1/p^n\,\Z/\Z)\,$ of the
finite abelian group $\,\fG\,$ of exponent $\,p^n\,$ with its natural structure of $\,R\,$-module given by 
$\,(\lam \psi)(g) = \psi(\lam g)\,$ for $\,\lam \in R, \psi \in \fG^\vee, g \in \fG\,$. Choose 
$\,\varphi \in \fG^\vee\,$ such that $\,\varphi\,|\,_\fH \neq 0\,$. The homomorphism of $\,R$-modules 
$\,R \lra \fG^\vee, \lam \mapsto \lam \varphi\,$, is injective: Indeed, assuming the contrary, let
$\,0 \neq \lam \in R\,$ be such that $\,\lam \varphi = 0\,$, i.e., $\,\varphi\,|\,_{\lam \fG} = 0\,$. As 
$\,R \sse {\rm End}(\fG)\,$, it follows that $\,\lam \fG\,$ is a non-zero $\,R$-submodule of $\,\fG\,$ and hence  
$\,\fH \sse \lam \fG\,$, so $\,\varphi\,|\,_\fH = 0\,$, which is a contradiction. As we already know that 
$\,|\fG^\vee| = |\fG| \leq |R|\,$, we deduce that $\,\Gam = 1 + \fm\,$ and the injective map above is an isomorphism, 
inducing by duality an isomorphism of $\,R$-modules (depending on the choice of $\,\varphi\,$)
\[ 
\fG \lra R^\vee = {\rm Hom}(R^+, 1/p^n\,\Z/\Z),\,g \mapsto \psi_g\,\mbox{with}\,\psi_g(\lam) = \varphi(\lam g).
\] 
The isomorphism $\,\fG \lra R^\vee\,$ maps $\,\fH = \fG^\Gam\,$ onto $\,(R^\vee)^\Gam = (R/\fm)^\vee \cong k^+\,$, 
inducing an isomorphism $\,\fG/\fH \cong \fm^\vee = \{\psi\,|\,_\fm\,:\,\psi \in R^\vee\}\,$. Composing succesively 
the bijective map $\,\fm \lra \Gam, \lam \mapsto 1 + \lam$, (the inverse of the canonical cocycle 
$\,\gam \mapsto \gam - 1$), the cocycle $\eta : \Gam \lra \fG$, the isomorphism above $\,\fG \lra R^\vee\,$ 
and the natural projection $\,R^\vee \lra \fm^\vee\,$, we obtain an isomorphism of $\,R$-modules 
$\,\fm \lra \fm^\vee, \lam \mapsto \lam^\vee\,$, with 
$\,\lam^\vee(\mu) = \varphi(\mu \eta(1 + \lam)) =$ $\,\varphi(\lam \eta(1 + \mu))\,$, inducing a nondegenerate pairing
$\,\fm \times \fm \lra \Q/\Z,\,(\lam, \mu) \mapsto \lam^\vee(\mu) = \mu^\vee(\lam)\,$, 
which is compatible with the canonical action of $\,R\,$ on $\,\fm$.

Assuming that $\,R\,$ is principal, $\,\fm = R \theta\,$ for any $\,\theta \in \fm \sm \fm^2\,$, and every ideal of 
$\,R\,$ is of the form $\,\fm^i = R \theta^i\,$, whence for each $\,a \in R\,$, there exists a unique nonnegative 
integer $\,i \leq m\,$ such that $\,a = u \theta^i\,$, where $\,u \in R^\times\,$ and $\,m\,$ is the nilpotency 
index of $\,\fm\,$. Putting $\,S := \Z/p^n\Z,\,R = S \oplus S \theta \oplus \cdots \oplus S \theta^{e - 1}\,$ 
is an $\,S$-module direct sum where $\,e$ is the greatest integer $\,i \leq m$ such that $\,p \in \fm^i\,$. 
It follows that $\,\theta\,$ satisfies an {\em Eisenstein polynomial}
\[
f(x) := x^e - p(a_{e - 1}x^{e - 1} + \cdots + a_0),
\]
where $\,a_i \in S, a_0 \in S^\times\,$. Thus $\,R \cong S[x]/(f(x), p^{n - 1}x^t)\,$, where
$\,1 \leq t:= m - (n - 1)e \leq e\,$. As $\,S$-module, $\,R\,$ is free of rank $\,e\,$ if and only if $\,t = e\,$,
i.e., $\,m = n e$.

As we already know that $\,\fG\,$ and $\,R^\vee\,$ are isomorphic $\,R$-modules, we deduce that $\,\fG\,$ is
free of rank $\,1\,$ since the finite commutative ring $\,R\,\,$ is principal by assumption. 
\ep
The next two results add new informations to those contained in Lemma~\ref{Lem:Ab1}.

\begin{pr}
\label{Pr:Ab1}
Let $\,p\,$ be a prime number, and $\,\eta \in Z^1(\Gam, \fG)\,$ be such that $\,\fG\,$ is a finite abelian $\,p$-group,
while the finite abelian group $\,\Gam\,$ is not a $\,p$-group. With notation from {\em Lemma~\ref{Lem:Ab1}}, the 
following assertions are equivalent.

{\rm (1)} $\,(\Gam, \fG, \eta)\,$ is a mnK triple.

{\rm (2)} $\,R \cong \F_q\,$ is a finite field with $\,q = p^f, f \geq 1\,$ for $\,p \neq 2,\,f \geq 2\,$ for 
$\,p = 2\,$, $\,\fG\,$ is a one-dimensional $\,R\,$-vector space identified with $\,R^+\,$, the group $\,\Gam\,$, 
identified with a subgroup of the multiplicative group $\,R^\times\,$, is cyclic of order $\,1\neq r\,|\,(q - 1)\,$ 
such that $\,f\,$ is the order of $\,p\,{\rm mod}\,r \in (\Z/r\Z)^\times\,$, and 
$\eta \in Z^1(\Gam, R^+) = B^1(\Gam, R^+) \cong (\Z/p\Z)^f\,$ is, up to multiplication by elements in $\,R^\times\,$, 
the coboundary $\,u \in \Gam \mapsto u - 1 \in R^+$.   
\end{pr}
\bp
$(1) \Lra (2)$. By Lemma~\ref{Lem:Ab1},(4), $\,R\,$ is a finite local ring; $\,k:= R/\fm\,$ its residue field.
By Lemma~\ref{Lem:Ab1},(3) again, there exists a unique non-zero minimal $\,R$-submodule $\,\fH$ of $\,\fG$, so  
$\,\fH = Rh\,$ for any $\,0 \neq h \in \fH\,$, and the surjective morphism of $\,R$-modules 
$\,R \lra \fH,\lam \mapsto \lam h\,$ induces an isomorphism $\,\fH \cong k^+\,$. Put $\,\Lam:= \eta^{- 1}(\fH)\,$.
As $\,(\Gam, \fG, \eta)\,$ is a mnK triple, it follows that $\,(\Gam : \Lam) = (\fG : \fH)\,$ is a $\,p$-th power, 
and hence $\,\Lam\,$ is not a $\,p$-group. In particular, $\,\Lam \neq 1\,$, and $\,\fH = R \eta(\si)\,$ for any 
$\,1 \neq \si \in \Lam\,$.

As $\,\Gam \sse R^\times\,$ is not a $\,p$-group, its $\,p$-primary component 
$\,\Gam(p) := \Gam \cap (1 + \fm) = {\rm Fix}\,_\Gam(\fH)\,$ is a proper subgroup of $\,\Gam\,$, therefore its complement
$\,\Gam(p') \cong \Gam/\Gam(p)\,$ is identified with a nontrivial subgroup of the multiplicative group $\,k^\times\,$; 
in particular, $\,|k|:= q = p^f > 2$, so $\,f \geq 2\,$ for $\,p = 2\,$. Thus $\,1 \neq r := |\Gam(p')|\,|\,(q - 1)\,$, 
and $\,\Gam(p') \cong \Z/r\Z\,$. Moreover, since $\,\Gam\,$ generates the local ring $\,R\,$, it follows that 
$\,k = \F_p(\Gam(p'))\,$, and hence $\,f\,$ is the order of $\,p\,{\rm mod}\,r \in (\Z/r\Z)^\times\,$. 

Let us show that $\,\Gam(p) \cap \Lam = 1\,$, and hence $\,\Lam = \Gam(p') \cong \Z/r\Z,\,|\Gam(p)| = (\fG : \fH)$. 
Assuming the contrary, let $\,h := \eta(\si)\,$ for some $\,1 \neq \si \in \Gam(p) \cap \Lam\,$, so $\,\fH = R h\,$. 
For any $\,\tau \in \Lam\,$, we obtain 
$\,\tau h = \eta(\tau \si) - \eta(\tau) = \eta(\si \tau) - \eta(\tau) = \eta(\si) + (\si - 1) \eta(\tau) = 
\eta(\si) = h\,$,
i.e., $\,\Lam\,$ acts trivially on $\,\fH\,$. Consequently, the restriction map $\,\eta|_{\Lam} : \Lam \lra \fH\,$ 
is a monomorphism, therefore $\,\Lam \cong \eta(\Lam)\,$ is a $\,p$-group, which is a contradiction.

It remains to show that $\,\Gam(p) = 1\,$, so $\,\Gam = \Lam,\,R \cong k, \,\fG = \fH \cong k^+$ as required. Assuming 
the contrary, let $\,\fG'\,$ be the subgroup of $\,\fG\,$ generated by $\,\eta(\Gam(p)) \neq 0\,$. $\,\fG'$ is obviously 
stable under the action of $\,\Gam(p)\,$ and also fixed by $\,\Gam(p') = \Lam\,$ since, acting on generators, we obtain 
$\,\tau \eta(\si) = \eta(\si) + (\si - 1)\eta(\tau) = \eta(\si)\,\,\mbox{for}\,\,\si \in \Gam(p) = 
{\rm Fix}_\Gam(\fH), \tau \in \Lam = \eta^{- 1}(\fH)\,$.
Thus $\,\fG'\,$ is a non-zero $\,\Gam$-submodule of $\,\fG\,$, and hence $\,\fH \sse \fG'\,$. On the other hand, 
the cocycle $\,\wt{\eta} : \Gam(p) \lra \fG/\fH\,$ induced by $\,\eta\,$ is bijective since 
$\,\Ker(\wt{\eta}) = \Gam(p) \cap \eta^{- 1}(\fH) = 1\,$ and $\,|\Gam(p)| = (\fG : \fH)\,$, therefore $\,\fG' = \fG\,$ 
and $\,1 \neq \Lam \sse {\rm Fix}_\Gam(\fG') = {\rm Fix}_\Gam(\fG) = 1\,$, which is a contradiction.

$(2) \Lra (1)$. For an arbitrary prime number $\,p\,$ and an arbitrary integer $\,r \geq 2\,$ such that $\,(p, r) = 1\,$ 
and $\,f \geq 2\,$ for $\,p = 2\,$, where $\,f\,$ is the order of $\,p\,{\rm mod}\,r \in (\Z/r\Z)^\times\,$, it follows 
easily that $\,(\Gam, k^+, \eta)\,$ is a mnK triple, where $\,k = \F_q,\,q := p^f,\,\Gam \cong \Z/r\Z\,$ is the unique 
subgroup of order $\,r\,$ of the multiplicative group $\,k^\times\,$, acting canonically on $\,k^+\,$, and 
$\,0 \neq \eta \in Z^1(\Gam, k^+) = B^1(\Gam, k^+) \cong (\Z/p\Z)^f\,$.
Consequently, the pairs $\,(p, r\,)$ as above classify, up to isomorphism, the mnK triples $\,(\Gam, G, \eta)\,$ with 
abelian $\,\Gam\,$ and $\,\fG\,$ such that $\,(|\Gam|,\,|\fG|) = 1$.  
\ep

\begin{pr}
\label{Pr:Ab2}
Let $\,\eta \in Z^1(\Gam, \fG)\,$, where $\,\Gam\,$ and $\,\fG\,$ are finite abelian $\,p$-groups for some
prime number $\,p\,$. With the notation from {\em Lemma~\ref{Lem:Ab1}}, the following assertions are equivalent.

{\rm (1)} $\,(\Gam, \fG, \eta\,)$ is a mnK triple, and the finite local ring $\,R\,$ of characteristic 
$\,p^n\,$ is principal.

{\rm (2)} One of the following conditions is satisfied.

{\rm (i)} $\,n = 1\,$, i.e., the characteristic of $\,R\,$ is $\,p\,$, and $\,R \cong \F_p[x]/(x^m)\,$
with $\,m \geq 2\,$ and $\,(m, p) = 1$.

{\rm (ii)} $\,p = n = 2\,$, i.e., the characteristic of $\,R\,$ is $\,4\,$, $\,R \cong (\Z/4\Z)[x]/(f(x), 2x^e)\,$,
where $\,f \in (\Z/4\Z)[x]\,$ is an Eisenstein polynomial of degree $\,e \geq 1\,$, the nilpotency index $\,m = 2 e\,$ 
is even, and $\,R\,$ is free of rank $\,e\,$ as $\,\Z/4\Z$-module.

{\rm (iii)} $\,p = n = 2\,$, i.e., the characteristic of $\,R\,$ is $\,4\,$,  
$\,R \cong (\Z/4\Z)[x]/(f(x), 2x^t)\,$, where $\,f \in (\Z/4\Z)[x]\,$ is an Eisenstein polynomial of degree 
$\,e \geq 2, 0 < t < e\,$, and the nilpotency index $\,m = e + t\,$ is odd.

In all three cases above, the $\,R$-module $\,\fG\,$ can be identified with $\,R^+,\,R/\fm \cong \F_p$,
$\,\Gam = 1 + \fm\,$, and the injective cocycle $\,\eta \in Z^1(\Gam, R^+)\,$ is unique up to multiplication with 
elements from $\,R^\times\,$ and summation with homomorphisms defined on $\,\Gam\,$ with values in 
$\,(R^+)^\Gam \cong \Z/p\Z$. 
\end{pr}
\bp
$(1) \Lra (2)$. By Lemma~\ref{Lem:Ab1},(5), we can identify the $\,R$-module $\,\fG\,$ with $\,R^+\,$, 
$\,\Gam = 1 + \fm\,$, $\,\eta \in Z^1(\Gam, R^+)\,$ is injective, and 
$\,R = (\Z/p^n\Z)[\theta] \cong \Z/p^n\Z [x]/(f(x), p^{n - 1}x^t)$, where $\,\theta R = \fm\,$ and $\,\theta\,$ 
is a root of the Eisenstein polynomial of degree $\,e\,$
\[
f(x) = x^e - p(a_{e - 1}x^{e - 1} + \cdots + a_0),\,a_i \in \Z/p^n\Z, a_0 \in (\Z/p^n\Z)^\times,
\]
and $\,1 \leq t = m - (n - 1)e \leq e\,$, where $\,m\,$ is the nilpotency index of the maximal ideal $\,\fm\,$. 
Choose $\,\wt{a_i} \in \Z_p\,$ such that $\,\wt{a_i} \equiv a_i\,{\rm mod}\,p^n, i = 0, \dots , e - 1\,$. 
The Eisenstein polynomial
\[
\wt{f}(x) = x^e - p(\wt{a_{e - 1}}x^{e - 1} + \cdots + \wt{a_0}) \in \Z_p[x]
\]
is irreducible over $\,\Q_p\,$. Let $\,\pi \in \wt{\Q_p}\,$ be a root of $\,\wt{f}\,$. Then $\,K = \Q_p(\pi)\,$ 
is a totally ramified extension of $\,\Q_p\,$ of degree $\,e\,$, with the valuation ring $\,\cO = \Z_p[\pi]\,$, the 
maximal ideal $\,\fp = \pi \cO\,$ and the residue field $\,k = \cO/\fp \cong \Z_p/p \Z_p \cong \F_p\,$.   
Sending $\,\pi\,$ to $\,\theta\,$, we obtain an epimorphism $\,\cO \lra R\,$, so 
$\,R \cong \cO_m := \cO/\fp^m,\,\fm = \fp \cO_m, \Gam \cong U^{(1)}/U^{(m)}\,$.
The cocycle $\,\eta : \Gam \lra \cO_m\,$ extends to a cocycle $\,\wt{\eta} : U^{(1)} \lra \cO_m\,$ with
$\,\Ker(\wt{\eta}) = U^{(m)}\,$. According to Example~\ref{Ex:Local2},
\[
\wt{\eta}(x) = \wt{\eta}_{u, \alpha}(x) \equiv \pi^{m - 1} \alpha(x) + \f{x - 1}{\pi} u\,{\rm mod}\,\fp^m,
\]
for some $\,u \in \cO_m^\times, \alpha \in T_{m + 1}^\vee = {\rm Hom}(T_{m + 1}, \F_p^+)\,$, with 
$\,T_{m + 1}:= \frac{U^{(1)}}{U^{(m + 1)} V (U^{(1)})^p}\,$, where $\,V\,$ denotes the closed subgroup 
of $\,U^{(1)}\,$ generated by $\,1 + \pi\,$. As $\,U^{(m + 1)} \sse \Ker(\wt{\eta}) \sse U^{(m)}\,$ and 
$\,m \geq 2\,$, the necessary and sufficient condition for the required equality $\,\Ker(\wt{\eta}) = U^{(m)}\,$
is that $\,U^{(m)} \cap (U^{(1)})^p \sse U^{(m + 1)}\,$ and the restriction $\,\alpha\,|\,_{U^{(m)}}\,$ is the 
nontrivial character $\,\chi_{\overline{u}} : U^{(m)}/U^{(m + 1)} \lra \F_p^+\,$, defined by
\[
\chi_{\overline{u}}(x) = \f{1 - x}{\pi^m} \overline{u}\,\,\mbox{for}\,x \in U^{(m)},\,\mbox{with}\,\overline{u}:= 
u\,{\rm mod}\,\fm.
\]

Assuming $\,m > e\,$, whence $\,n \geq 2\,$, put $\,a := (1 + \pi^{m - e})^p \in (U^{(1)})^p\,$. Since
\[
v(a - 1) \geq {\rm min}\,(v(p \pi^{m - e}), v(\pi^{(m - e)p})) = {\rm min}\,(m, (m -e)p),
\] 
the condition above implies the inequality $\,(m - e)p \leq m\,$, therefore $\,(n - 1)e < m \leq \f{pe}{p - 1}\,$. 
It follows that $\,2 \leq n \leq \f{p}{p - 1} = 1 + \f{1}{p - 1}\,$, so $\,n = p = 2\,$. Consequently, if 
$\,p \neq 2\,$ then $\, n = 1\,$, i.e., $\,{\rm char}\,R = p\,$, and $\,m = e = t \geq 2, R \cong \F_p[x]/(x^m)\,$, 
so we may take $\,\wt{f}(x) = x^m - p\,$, $\,K = \Q_p(p^{\f{1}{m}})\,$.

Next assume $\,p\,|\,m\,$, and put $\,b := (1 + \pi^{\f{m}{p}})^p \in U^{(1)}\,$. Since
\[
v(b - 1) \geq {\rm min}\,(v(p \pi^{\f{m}{p}}),
v(\pi^m)) = {\rm min}\,(e + \f{m}{p}, m),
\] 
we deduce that $\,e + \f{m}{p} \leq m\,$, i.e., $\, m \geq \f{p e}{p - 1}\,$, since otherwise
$\,b \in (U^{(m)} \cap (U^{(1)})^p) \sm U^{(m + 1)}\,$, i.e., a contradiction. Consequently, if $\,p \neq 2\,$ 
then $\,(m, p) = 1\,$ as stated by (i), i.e., the totally ramified extension $\,K\,|\,\Q_p\,$ is tame,
since otherwise we obtain $\,e = m \geq \f{pe}{p - 1}\,$, a contradiction. If $\,p = 2\,|\,m\,$ then 
$\,m \geq 2 e\,$, in particular, $\,m > e\,$, and hence $\,n = 2\,$ and $\,m \leq 2e\,$, so $\,m = 2e\,$, as stated
by (ii). If $\,p = 2 \nmid\, m\,$, in particular, $\,m \geq 3\,$, then either $\,n = 1\,$ (case (i)) or
$\,n = 2, e < m < 2e\,$ (case (iii)).

According to \cite[Theorem 2]{CL}, the principal local ring $\,R \cong \cO_m\,$ is determined up to isomorphism 
by its invariants $\,p, n, m\,$ only in the case (i), the cases (ii) and (iii) with $\,e\,$ odd, and the case (iii) 
with $\,t = 1\,$, i.e., $\,e\,$ even and $\,m = e + 1\,$. As shown above, the cocycles 
$\,\eta: \Gam = 1 + \fm \lra \fG = R^+\,$ for which $\,(\Gam, \fG, \eta)\,$ is a mnK triple, i.e., $\,\eta\,$ is 
injective, are in 1-1 correspondence with the pairs $\,(u, \alpha) \in R^\times \times T_{m + 1}^\vee\,$ satisfying 
$\,\alpha|_{U^{(m)}} = \chi_{\overline{u}}\,$. The condition $\,U^{(m)} \cap (U^{(1)})^p \sse U^{(m + 1)}\,$
is equivalent with the fact that $\,U^{(m)}/U^{(m + 1)} \cong \Z/p\Z\,$ is the kernel of the canonical projection 
$\,T_{m + 1} \lra T_{m}\,$, so $\,T_{m + 1}^\vee \cong T_m^\vee \oplus (U^{(m)}/U^{(m + 1)})^\vee\,$. Thus the 
cocycles above are parametrized by the elements of the direct product 
$\,R^\times \times T_m^\vee =\,$ $\,R^\times \times (\Gam/\la 1 + \theta \ra \Gam^p)^\vee\,$, and hence they are, up to 
isomorphism in $\,\cZ^1$, in 1-1 correspondence to the elements of the elementary $\,p$-group 
$\,(\Gam/\la 1 + \theta \ra \Gam^p)^\vee\,$. Consequently, the cocycle $\,\eta\,$ is unique up to multiplication by 
elements of $\,R^\times\,$ and summation with homomorphisms from $\,\Gam\,$ to $\,\fG^\Gam \cong \Z/p\Z\,$ as desired.

$(2) \Lra (1)$. Assume the finite principal local ring $\,R\,$ satisfies one of the conditions (i)-(iii). 
As in the first part of the proof, we choose a suitable finite extension $\,K\,$ of $\,\Q_p\,$ with valuation ring 
$\,\cO\,$, maximal ideal $\,\fp = \pi \cO\,$, and residue field $\,k \cong \F_p\,$, such that 
$\,R \cong \cO_m = \cO/\fp^m\,$, where $\,m\,$ is the nilpotency index of the maximal ideal $\,\fm\,$ of $\,R\,$. 
Thanks to the arguments from the first part of the proof, it suffices to check the condition 
$\,U^{(m)} \cap (U^{(1)})^p \sse U^{(m + 1)}\,$. Assuming the contrary, let $\,x \in \fp, c := (1 + x)^p - 1\,$ 
be such that $\,v(c) = m\,$. We consider separately the cases (i)-(iii).

Case (i). It follows that
\[ 
m = v(c) \geq\,{\rm min}\,(v(px), v(x^p)) = {\rm min}\,(m + v(x), p v(x)) \geq {\rm min}\,(m + 1, p v(x)),
\]
therefore $\,m = p v(x)\,$, contrary to the hypothesis $\,(m, p) = 1$.

Case (ii). We obtain $\,c = (1 + x)^2 - 1 = x(x + 2)\,$, so
\[ 
2 e = m = v(c) = v(x) + v(x + 2)\geq v(x) +\,{\rm min}\,(v(x), e),
\]
hence $\,v(x) \leq e\,$. Assuming $\,v(x) < v(2) = e\,$, it follows that $\,2 e = v(c) = 2 v(x)\,$, which is a 
contradiction. Thus $\,v(x) = e\,$, therefore $\,v(1 + \frac{2}{x}) = v(c) - 2 v(x) = 0\,$, contrary to the fact that
$\,1 + \frac{2}{x} \equiv 2 \equiv 0\,{\rm mod}\,\fp\,$ since $\,\cO/\fp \cong \F_2$. 

Case (iii). It follows that $\,m = v(c) \geq v(x) + {\rm min}\,(e, v(x))\,$, therefore
\[ 
e > t = m - e \geq v(x) +\,{\rm min}\,(0, v(x) - e),
\]
so $\,v(x) < e\,$. Consequently, $\,m = v(c) = 2 v(x)\,$, again a contradiction, as $\,m\,$ is odd by hypothesis. 
\ep

\begin{res} \em
(1) For large enough prime numbers $\,p\,$, an effective description of the corresponding mnK triples 
is easy. Indeed, let us consider the case (i) above with $\,p > m\,$. Then $\,\Gam = 1 + \fm \cong (\Z/p\Z)^{m - 1}\,$ 
with a base consisting of the elements $\,\gam_i := 1 + \theta^i,\,i = 1, \dots, m - 1\,$, and 
$\,R^\times \cong \Gam \times \F_p^\times\,$. The canonical injective cocycle $\,\eta : \Gam \lra R^+\,$
is completely determined by its values $\,\eta(\gam_i) = \theta^{i - 1}, i = 1, \dots, m - 1\,$, and every injective 
cocycle $\,\Gam \lra R^+\,$ has the form $\,u (\eta + \theta^{m - 1} \beta)\,$ with $\,u \in R^\times\,$, 
$\,\beta \in {\rm Hom}\,(\Gam/\la \gam_1 \ra, \F_p^+) \cong (\Z/p\Z)^{m - 2}$.

(2) The classification of the mnK triples for which the local ring $\,R\,$ is not principal seems to be a more 
difficult task. We give only a simple example of such triples. Let $\,W\,$ be a vector space over the prime field 
$\,\F_p, p \neq 2\,$, with the base $\,\{ \theta_i\,|\, i = 1, \dots, s \}, s \geq 2\,$, and let 
$\,R = \F_p \oplus W\,$ with the multiplication given by $\,(x \oplus y) \cdot (x^{\prime} \oplus y^{\prime}) =
xx^{\prime} \oplus (xy^{\prime} + x^{\prime}y)\,$. $\,R\,$ is a local ring with the maximal ideal $\,\fm = W\,$ whose 
nilpotency index is $\,2\,$. The canonical map $\,\fm \lra \Gam = 1 + \fm, x \mapsto 1 + x\,$ is an isomorphism, 
the elements $\,\gam_i := 1 + \theta_i, i = 1, \dots, s\,$, form a base of the elementary $\,p$-group $\,\Gam\,$, 
and $\,R^\times \cong \F_p^\times \times \Gam \cong \Z/(p - 1)\Z  \times (\Z/p\Z)^s\,$. Setting $\,\theta_0 = 1\,$, 
the Pontryagin dual $\,\fG := R^\vee = {\rm Hom}\,(R^+, \F_p)\,$ is a vector space over $\,\F_p\,$ with the dual base 
$\,\{ \theta_i^\vee\,|\, i = 0, \dots, s \}\,$. The canonical $\,R$-module structure on $\,\fG\,$ is defined by the 
relations $\,\theta_0 \theta_j^\vee = \theta_j^\vee, \theta_i \theta_j^\vee = \delta_{i\,j} \theta_0^\vee,
i = 1, \dots, s, j = 0, \dots, s\,$, where $\,\delta_{i\,j}$ is the Kronecker symbol. It follows that
$\,\fH:= \fG^\Gam = \F_p \theta_0^\vee \cong \Z/p\Z\,$, and the induced action of $\,\Gam\,$ on the quotient 
$\,\fG/\fH \cong (\Z/p\Z)^s\,$ is trivial. Every cocycle $\,\eta \in Z^1(\Gam, \fG)\,$ is completely determined by 
its values $\,\eta(\gam_i), i = 1, \dots, s$, which must satisfy the condition 
$\,\theta_i \eta(\gam_j) = \theta_j \eta(\gam_i)\,$ for $\,1 \leq i, j \leq s$. Writing 
$\,\eta(\gam_i) = \dsum_{j = 0}^{s} \lam_{i\,j} \theta_j^\vee\,$, the cocycles $\,\eta\,$ are in 1-1 correspondence 
to the pairs $\,(\lam_0, \Lam)\,$ consisting of a homomorphism $\,\lam_0: \Gam \lra \F_p, \gam_i \mapsto \lam_{i\,0}\,$,
and a symmetric matrix $\,\Lam = ( \lam_{i\,j})_{1 \leq i, j \leq s}\,$ with entries in $\,\F_p$.
 
$(\Gam, \fG, \eta)\,$ is a mnK triple $\,\Llra \eta\,$ is injective $\,\Llra\,$ the matrix $\,\Lam\,$ is invertible
$\,\Llra\,{\rm det}\,(\Lam) \neq 0 \Llra\,$ the quadratic form 
$\,Q(x) = \dsum_{1 \leq i, j \leq s} \lam_{i\,j} x_i x_j\,$ in $\,x_1, \dots, x_s$ is nondegenerate. Consider the 
group $\,\cG := {\rm Aut}\,(\Gam, \fG)\,$ consisting of the pairs $\,(\varphi, \psi)\,$, where 
$\,\varphi \in {\rm Aut}(\Gam)\,$ and $\,\psi\in {\rm Aut}(\fG)\,$ satisfy the condition
$\,\psi(\varphi(\gam) g) = \gam \psi(g)\,$ for $\,\gam \in \Gam, g \in \fG\,$, with the composition law 
$\,(\varphi, \psi) \circ (\varphi', \psi') = (\varphi' \circ \varphi, \psi \circ \psi')$. The elements of $\,\cG\,$ 
are in 1-1 correspondence with the triples 
$\,(A = (a_{i\,j})_{1 \leq i, j \leq s} \in GL_s(\F_p), b_0 \in \F_p^\times, {\bf b} = (b_1, \dots, b_s) \in \F_p^s)\,$,
assigning to such a triple $\,(A, b_0, {\bf b})\,$ the pair $\,(\varphi, \psi)\,$ defined by 
\[
\varphi(\gam_i) = \dprod_{j = 1}^{s} \gam_j^{a_{i\,j}}, \psi(\theta_0^\vee) = b_0 \theta_0^\vee, \psi(\theta_i^\vee) = 
b_0 (b_i \theta_0^\vee + \dsum_{j = 1}^{s} a_{j\,i} \theta_j^\vee), i = 1, \dots, s.
\]
The group $\,\cG\,$ acts on the injective cocycles sending a cocycle $\,\eta\,$ to the cocycle 
$\,\psi \circ \eta \circ \varphi\,$. The action is transitive whenever $\,s\,$  is odd and hence in this case the 
mnK triple $\,(\Gam, \fG, \eta)\,$ is unique up to isomorphism. Indeed let $\,\eta_0\,$ be the canonical injective 
cocycle defined by $\,\eta_0(\gam_i) = \theta_i^\vee, i = 1, \dots, s\,$, and $\,\eta\,$ be an arbitrary injective 
cocycle defined by the pair$\,(\lam_0, \Lam)\,$. Put 
\begin{center}
$\,b_0 = \left\{ \begin{array}{lcl}
1 & \mbox{if} & {\rm det}\,(\Lam) \in \F_p^2  \\
\,u \in \F_p \sm \F_p^2 & \mbox{if} & {\rm det}\,(\Lam) \not \in \F_p^2
\end{array}
\right.$
\end{center}
By \cite[1.7, Proposition 5]{Serre-A}, there exists $\,A \in GL_s(\F_p)\,$ such that $\,\Lam = A \cdot (b_0 I_s) \cdot\, ^t A\,$,
where $\,^t A\,$ denotes the transpose of $\,A\,$. Let $\,{\bf b} = (b_1, \dots, b_s) \in \F_p^s\,$ be the
unique solution of the linear system in $\,x_1, \dots, x_s$
\[
\dsum_{j = 1}^{s} a_{i\,j} x_j = b_0^{- 1} \lam_{i\,0} +
\dsum_{j = 1}^{s} \f{a_{i\,j} (1 - a_{i\,j})}{2},\,i = 1, \dots, s.
\] 
One cheks easily that $\,\eta = \psi \circ \eta_0 \circ \varphi\,$, where the pair $\,(\varphi, \psi) \in \cG\,$ is
determined by the triple $\,(A, b_0, {\bf b})\,$ as defined above. Similarly, it follows for $\,s\,$ even that 
$\,(\Gam, \fG, \eta_i)\,$, $\,i = 0, 1\,$, are up to isomorphism the only two mnK triples associated to the pair 
$\,(\Gam, \fG)\,$, where the injective cocycle $\,\eta_1\,$ is defined by
$\,\eta_1(\gam_i) = \theta_i^\vee, i = 1, \dots, s - 1, \eta_1(\gam_s) = u \theta_s^\vee\,$ with 
$\,u \in \F_p \sm \F_p^2$. 
\end{res}
\begin{co}
\label{Cor:Character}
Let $\,(\Gam, \fG, \eta)\,$ be a mnK triple, and assume $\,\fG\,$ is abelian of exponent $\,k\,$ and the action 
of $\,\Gam\,$ is induced by a character $\,\chi: \Gam \lra (\Z/k\Z)^\times\,$, i.e., $\,\gam g = \chi(\gam) g$ 
for $\,\gam \in \Gam, g \in \fG\,$. Then either $\,k = p\,$ is an odd prime number or $\,k = 4\,$, and 
$\,\fG \cong \Z/k \Z\,$.

If $\,k = p \neq 2\,$ then $\,(\Gam, \fG, \eta) \cong (U, \F_p^+, u \mapsto u - 1)\,$, where 
$\,U \cong \Gam \cong \Z/r\Z\,$, $2 \leq r\,|\,(p - 1)\,$, is the unique subgroup of order $\,r =\,|\Gam|\,$ of the
multiplicative group $\,\F_p^\times\,$ acting by multiplication on $\,\F_p^+$.

If $\,k = 4\,$ then $\,(\Gam, \fG, \eta) \cong (\Z/2\Z, \Z/4\Z,\,1\,{\rm mod}\,2 \mapsto 1\,{\rm mod}\,4)\,$ with
the unique non-trivial action of $\,\Z/2\Z\,$ on $\,\Z/4\Z$.
\end{co} 
\bp
Since $\,(\Gam, \fG, \eta)\,$ is a mnK triple and $\,\fG\,$ is abelian, it follows by Lemma~\ref{Lem:Nilpotent},(2) that
$\,{\rm Fix}_\Gam(\fG) = 1\,$, and hence $\,\Gam \cong \chi(\Gam) \sse (\Z/k \Z)^\times\,$ is abelian. It remains
to apply Lemma~\ref{Lem:Ab1}, Proposition~\ref{Pr:Ab1} and Proposition~\ref{Pr:Ab2}.
\ep
As a consequence of Corollary~\ref{Cor:Character} and Proposition~\ref{Pr:AbsKneserCrit}, we find again 
\cite[Lemma 1.18, Theorem 1.20]{A-B}, the abstract version in the framework described in Example~\ref{Ex:1} of the
classical Kneser criterion for separable radical extensions \cite{K}, \cite[Theorem 11.1.5]{CT}.

\end{document}